\documentclass[final,3p,times]{elsarticle}
\usepackage{amssymb}
\usepackage{amsmath}
\usepackage{amsthm}
\usepackage{booktabs}
\usepackage{xcolor}
\journal{Journal of Computational Physics}
\usepackage{subfig}
\usepackage[numbers]{natbib}
\usepackage{graphicx}
\usepackage{url}
\usepackage{caption}
\usepackage{multirow}

\usepackage{pifont}

\usepackage{xcolor}
\usepackage{multirow}

\definecolor{goodgreen}{RGB}{0,120,0}
\definecolor{badred}{RGB}{180,0,0}

\usepackage{algorithm}
\usepackage{algpseudocode}  % 算法伪代码样式

\usepackage{float}
\usepackage{zymacros}

\usepackage{tabularx}
\usepackage{booktabs}
\usepackage{array}
\usepackage{amsmath,amssymb,mathtools}
\usepackage{bm}

\usepackage{booktabs}
\usepackage{tabularx}

\usepackage{cleveref}

\usepackage{amsthm}

\theoremstyle{plain}
\newtheorem{theorem}{Theorem}[section]
\newtheorem{lemma}[theorem]{Lemma}
\newtheorem{proposition}[theorem]{Proposition}
\newtheorem{corollary}[theorem]{Corollary}

\theoremstyle{definition}
\newtheorem{assumption}[theorem]{Assumption}

\theoremstyle{remark}
\newtheorem{remark}[theorem]{Remark}

\usepackage[most]{tcolorbox} % 在导言区添加

\tcbset{colback=gray!4,colframe=gray!40,boxrule=0.4pt,arc=2pt,
        left=4pt,right=4pt,top=4pt,bottom=4pt}

\begin{document}
\begin{frontmatter}

\title{Prescribed-Basis Coefficient-to-Coefficient Neural Operator for Partial Differential Equations} 

\author[1]{Chuqi Chen\corref{cor1}}
\ead{cchenck@umich.edu}

\author[2,3]{Yang Xiang\corref{cor1}} 
\ead{maxiang@connect.ust.hk}

\author[2]{Weihong Zhang\corref{cor1}}
\ead{wzhangde@connect.ust.hk}

\cortext[cor1]{Corresponding author.}
%\fntext[fn1]{These authors contributed equally to this work.}

%% Author affiliation
\affiliation[1]{organization={Department of Mathematics},
            addressline={University of Michigan}, 
            city={Ann Arbor},
            state={MI},
            country={USA}}
\affiliation[2]{organization={Department of Mathematics},
            addressline={The Hong Kong University of Science and Technology}, 
            city={Clear Water Bay}, 
            country={Hong Kong SAR}}
\affiliation[3]{organization={Algorithms of Machine Learning and Autonomous Driving Research Lab},
            addressline={HKUST Shenzhen-Hong Kong Collaborative Innovation Research Institute}, 
            city={Shenzhen},
            country={China}}

\footnotetext[1]{Author names are ordered alphabetically to indicate equal contribution.}

%% Abstract

\begin{abstract}
Operator learning provides a data-driven approach to approximating solution operators of partial differential equations, but its effectiveness depends strongly on how input and output functions are represented. 
Point-value representations can make the trainable map mesh-dependent and high-dimensional; snapshot-based POD/PCA reductions require aligned data and basis construction, while learned representations introduce additional trainable encoders, decoders, or neural bases.
We propose the \emph{Fixed-Basis Coefficient-to-Coefficient Network} (FB-C2CNet), which learns PDE solution maps in fixed, data-independent approximation spaces using prescribed bases as function encoders and decoders. 
Input observations are encoded by regularized least-squares projection onto bases such as finite element, random-feature, or radial-basis-function bases. 
A neural network maps the resulting input coefficients to output coefficients, and the fixed decoder reconstructs the solution at arbitrary target locations. 
This separation of basis selection from network training avoids neural basis learning and snapshot-based basis extraction, reduces the dimension of the trainable map, and lowers training cost. 
With suitable prescribed bases, FB-C2CNet also accommodates scattered, non-aligned, and sample-dependent observations.
We analyze the stability--bias trade-off of regularized coefficient encoding and the intrinsic projection error determined by the output space. 
Experiments on elliptic, nonlinear time-dependent, weak-solution, high-dimensional, and inverse Stokes boundary-recovery problems demonstrate competitive accuracy with reduced trainable dimension and training cost, including for high-resolution and irregularly sampled data.
\end{abstract}

%% Keywords

\end{frontmatter}

%%%%%%%%%%%%%%%%%%%%%%%%%
%%%%%%%%%%%%%%%%%%%%%%%%%

\section{Introduction}

Deep learning has recently emerged as a powerful paradigm in scientific computing, with applications ranging from solving differential equations to accelerating large-scale simulations and data-driven modeling~\cite{lu2021deepxde,pathak2022fourcastnet,bi2023accurate,lee2025physics}. 
Among various methodologies, operator learning has attracted significant attention because it aims to approximate mappings between function spaces. 
Unlike neural PDE solvers such as Physics-Informed Neural Networks (PINNs)~\cite{pinn,karniadakis2021physicsinformed,cuomo2022scientific}, the Deep Ritz Method~\cite{deepritz}, Weak Adversarial Networks (WANs)~\cite{wan,bao2025pfwnn}, and random-feature-based solvers~\cite{chen2022bridging,dong2023method,sun2024local,chen2024quantifying}, which usually approximate one solution instance at a time, neural operators learn solution maps from input functions to output functions. 
Typical inputs include coefficient fields, source terms, boundary data, or initial conditions, while the outputs are the corresponding PDE solution fields. 
Representative neural operator architectures include DeepONet~\cite{lu2021deeponet}, Fourier Neural Operator (FNO)~\cite{li2021fourier}, and related variants~\cite{lu2022comprehensive,li2020neural,tripura2023wavelet,fanaskov2023spectral,son2025elm}. 
For example, one may consider the elliptic solution operator
\[
\mathcal{G}:(a,f)\mapsto u,
\]
where $u$ solves
\[
-\nabla\cdot(a\nabla u)=f
\]
under appropriate boundary conditions. 
Once trained, an operator-learning model can be repeatedly evaluated for new input functions and can therefore serve as a fast surrogate for parametric or data-driven PDE simulations.

Despite this potential, a central issue in neural operator design is the finite-dimensional representation of input and output functions. 
Since the underlying PDE solution operators act between infinite-dimensional function spaces, a practical neural architecture must first choose how functions are encoded and decoded. 
A common strategy is to represent functions by their pointwise values on a fixed grid or a fixed set of sensors. 
This pointwise discretization is simple and widely used, but the neural input and output dimensions are directly tied to the spatial or spatiotemporal resolution. 
For high-resolution two- and three-dimensional problems, this can lead to very large input and output vectors, increased network sizes, higher memory usage, and longer training times. 
Moreover, when observations are scattered, nonuniform, or sample-dependent, a pointwise representation on a common grid may require additional interpolation, mesh transfer, or resampling procedures.

Reduced-coordinate representations provide a natural way to alleviate this dimensionality issue~\cite{benner2015survey,hesthaven2016certified,swischuk2019projection}. 
Classical reduced-order and POD/PCA-type neural models represent solution fields in low-dimensional coordinates and learn maps between these coordinates~\cite{hesthaven2018non,bhattacharya2021model,lu2022comprehensive}. 
However, standard snapshot-based POD/PCA reductions construct data-dependent bases from training fields represented on a common discretization, typically through a global SVD or PCA step~\cite{audouze2009reduced,bhattacharya2021model}. 
This global basis construction can become expensive when the snapshot matrix is large, for example for large datasets, high-resolution fields, or high-dimensional domains, and the resulting representation is tied to the snapshot discretization. 
Another line of work learns basis functions, trunk networks, encoders, or decoders as part of the neural operator architecture~\cite{lu2021deeponet,hua2023basis,ingebrand2025basis,bahmani2025resolution}. 
Such learned-basis approaches can be flexible, but they introduce additional trainable components and increase the optimization cost.

In this work, we take a different route and propose the \emph{Fixed-Basis Coefficient-to-Coefficient Network} (FB-C2CNet). 
The central idea is to learn PDE solution operators in prescribed, data-independent approximation spaces. 
Instead of learning basis functions or extracting them from training snapshots, we choose the input and output bases before training, for example from finite element, random-feature, or radial-basis-function approximation spaces. 
The basis selection is therefore a numerical approximation design choice rather than a trainable part of the neural operator.

For notational simplicity, we first describe the single-input/single-output setting. 
Let
\[
G:X\to Y,\qquad f\mapsto u=G(f),
\]
be an operator mapping an input function \(f\in X\) to an output function \(u\in Y\). 
We approximate the input and output function spaces by two prescribed finite-dimensional spaces
\[
V_{\mathrm{in}}^{m_1}
=
\operatorname{span}\{\varphi_j\}_{j=1}^{m_1},
\qquad
V_{\mathrm{out}}^{m_2}
=
\operatorname{span}\{\psi_j\}_{j=1}^{m_2}.
\]
Given sampled observations of $f$, the input coefficient vector $\bm a\in\mathbb{R}^{m_1}$ is obtained by a regularized least-squares projection onto $V_{\mathrm{in}}$, so that
\[
f(x)\approx \sum_{j=1}^{m_1} a_j \varphi_j(x).
\]
A neural network \(N_\theta\) is then trained only to learn the finite-dimensional coefficient map
\[
N_\theta:\mathbb{R}^{m_1}\to\mathbb{R}^{m_2},
\qquad
\bm a\mapsto \bm b.
\]
The predicted output field is reconstructed by the prescribed output basis:
\[
\widehat{u}_\theta(y)
=
\sum_{j=1}^{m_2}
\big[N_\theta(\bm a)\big]_j\psi_j(y).
\]
Thus, FB-C2CNet does not train the input basis, output basis, encoder, or decoder. 
The only trainable component is the coefficient-to-coefficient map \(N_\theta\).

This prescribed-basis formulation has several practical advantages. 
First, FB-C2CNet avoids both neural basis learning and snapshot-based POD/PCA basis extraction. 
As a result, the dimension of the trainable map is determined by the chosen approximation spaces rather than by the number of grid points, leading to fewer trainable parameters and lower training cost. 
Second, because the encoder and decoder are defined through evaluations of prescribed basis functions, the method can encode scattered, non-aligned, and sample-dependent input observations, and it can reconstruct outputs at arbitrary target locations. 
This property is particularly useful for sparse sensing and inverse problems, where measurements may be available only at irregular sensor locations. 
Third, the basis spaces can be selected according to the geometry, boundary conditions, sampling pattern, and dimension of the PDE problem. 
For example, finite element bases provide a natural interface with existing numerical discretizations on meshed domains, random-feature bases are convenient for mesh-free and high-dimensional settings, and radial basis functions are useful for scattered-data and geometric-field representations.

The fixed-basis viewpoint also makes the main error sources transparent. 
The performance of FB-C2CNet depends on the approximation capability of the prescribed input and output spaces, the stability and bias of the regularized coefficient recovery, and the learning error of the coefficient map. 
In this paper, we analyze the stability--bias trade-off induced by regularized input encoding and identify the output projection error associated with the prescribed decoder as an intrinsic accuracy floor. 
This separation connects neural operator learning with classical numerical approximation and provides practical guidance for basis selection and regularization.

\begin{figure}[!ht]
\centering
\includegraphics[width=0.9\textwidth]{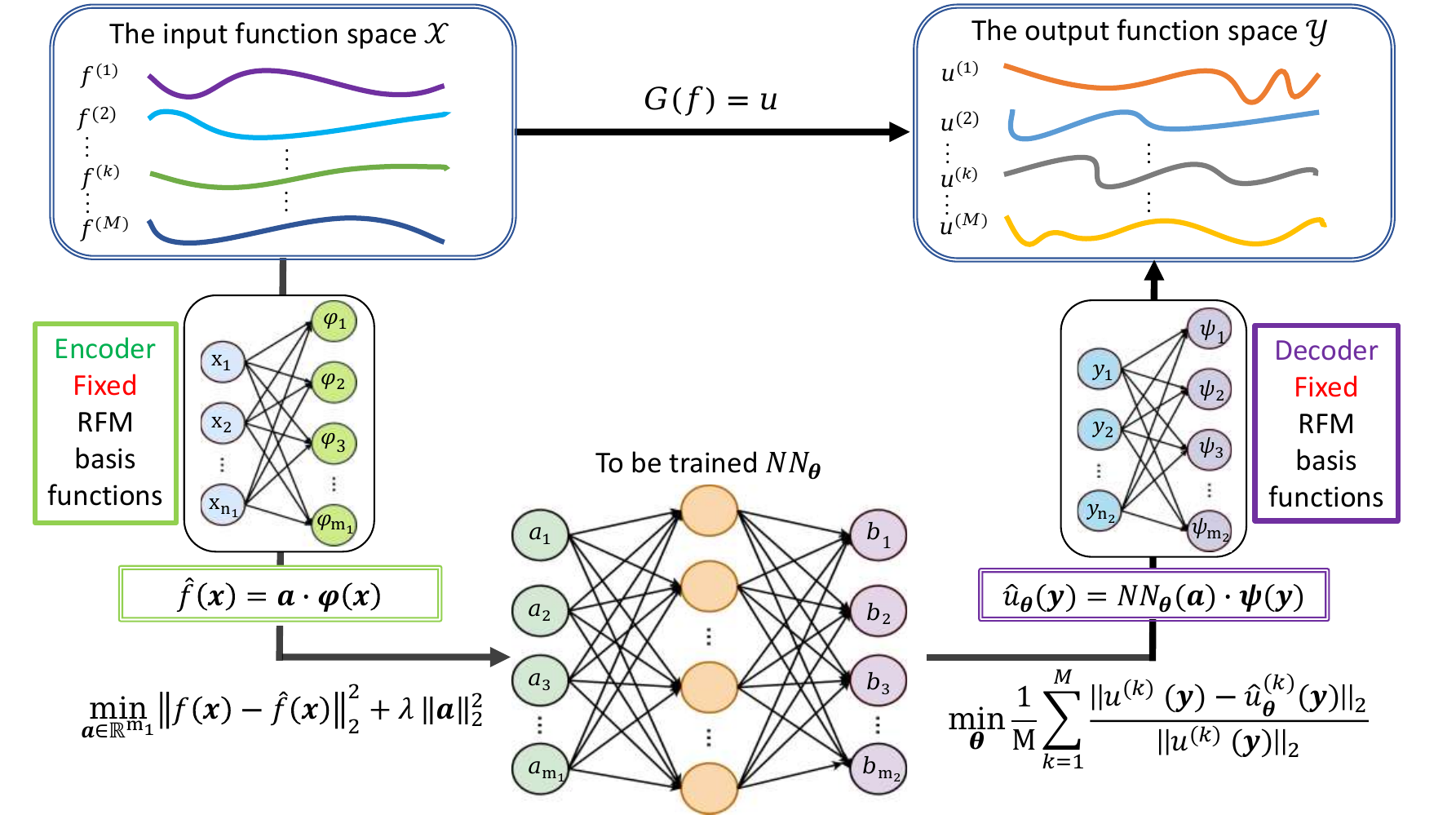}
\caption{
Workflow of the Fixed-Basis Coefficient-to-Coefficient Network (FB-C2CNet), illustrated with random-feature basis functions.
Input observations are encoded by a fixed prescribed-basis encoder into input coefficients; a neural network learns only the coefficient-to-coefficient map; and the output field is reconstructed by a fixed prescribed-basis decoder. 
    Both the encoder and decoder are non-trainable and are defined through evaluations of the prescribed basis functions.
}
\label{fig:workflow_FB-C2CNet}
\end{figure}

The main idea of FB-C2CNet is illustrated in Figure~\ref{fig:workflow_FB-C2CNet}. 
A detailed description of the proposed methodology is presented in Section~2. 
This section continues with a review of related work, followed by a summary of the main contributions of this paper.

\subsection{Related work}

Neural-network-based methods for partial differential equations can be broadly divided into two categories.
The first category aims to approximate individual PDE solutions.
Representative examples include Physics-Informed Neural Networks (PINNs)~\cite{pinn,karniadakis2021physicsinformed}, the Deep Ritz method~\cite{deepritz}, Weak Adversarial Networks (WANs)~\cite{wan,bao2025pfwnn}, and random-feature-based solvers~\cite{chen2022bridging,dong2023method,sun2024local}.
These methods usually train a neural network for a particular PDE instance by incorporating strong, weak, or variational forms of the governing equations into the loss function.
The second category is operator learning, which aims to approximate solution operators mapping input functions, such as coefficient fields, source terms, boundary data, or initial conditions, to output solution fields.
Representative neural operators include DeepONet~\cite{lu2021deeponet}, Fourier Neural Operators (FNOs)~\cite{li2021fourier}, Graph Neural Operators (GNOs)~\cite{li2020neural}, Wavelet Neural Operators (WNOs)~\cite{tripura2023wavelet}, and Spectral Neural Operators (SNOs)~\cite{fanaskov2023spectral}.
The present work belongs to the operator-learning framework.

A major design issue in neural operators is how input and output functions are represented.
DeepONet represents the input function by its values at sensor locations through a branch network and evaluates the output through a trunk network depending on the target coordinate~\cite{lu2021deeponet}.
POD-DeepONet replaces the learned trunk representation by POD modes extracted from training data~\cite{lu2022comprehensive}.
Other variants introduce alternative output representations, such as fixed bases from extreme learning machines~\cite{son2025elm}, finite element shape functions in mesh-informed neural networks~\cite{yamazaki2025finite}, or random-feature and partition-of-unity representations~\cite{chen2024quantifying}.
Although these methods differ in their output decoders, the input function is typically still provided to the branch network through pointwise values at sensors or grid points.
Consequently, the input dimension grows with the number of sensors, which can increase the parameter count and training cost in high-resolution or high-dimensional settings.

Reduced-coordinate and model-reduction-based approaches address this issue by representing functions in low-dimensional coordinates.
Classical reduced-order models and POD/PCA-type neural methods learn maps between reduced coefficients rather than pointwise fields~\cite{audouze2009reduced,benner2015survey,swischuk2019projection,hesthaven2016certified,hesthaven2018non,bhattacharya2021model}.
These approaches are closely related to our work in spirit, since they also reduce operator learning to a finite-dimensional map between coefficient vectors.
However, standard snapshot-based POD/PCA reductions construct data-dependent bases from training fields represented on a common discretization, typically through a global SVD or PCA step.
Thus, the resulting basis is tied to the snapshot grid, and scattered or non-aligned observations usually require additional interpolation, mesh transfer, or projection procedures.
For large snapshot matrices arising from high-resolution fields, large datasets, or high-dimensional domains, the basis construction itself can also become a non-negligible computational cost.

Another line of work learns basis functions, encoders, decoders, or coordinate-dependent representations as part of the neural operator architecture.
BasisONet represents input functions by coefficients with respect to a basis and trains a trunk network to construct neural basis functions for the output~\cite{hua2023basis}.
Basis-to-Basis (B2B) operator learning uses function encoders to obtain coefficient representations and learns maps between basis representations~\cite{ingebrand2025basis}.
BelNet~\cite{zhang2023belnet}, NOMAD~\cite{seidman2022nomad}, C2BNet~\cite{zhang2025coefficient}, and resolution-independent neural operators~\cite{bahmani2025resolution} also develop flexible representations for mesh-free or resolution-independent operator learning.
These methods provide powerful mechanisms for handling flexible discretizations and arbitrary target coordinates, but they often introduce additional trainable representation modules, such as learned encoders, learned decoders, trunk networks, or dictionary-learning components.

FB-C2CNet can be viewed as a coefficient-to-coefficient neural operator, but its main distinction lies in how the coefficient spaces are constructed.
The input and output bases are prescribed before training and are not extracted from training snapshots or learned as part of the neural architecture.
Input observations are encoded into coefficients by a fixed coefficient-recovery procedure, and the output field is decoded from the predicted coefficients using a fixed prescribed basis.
Therefore, only a single neural network, namely the coefficient-to-coefficient map, is trained.
This design separates numerical approximation from neural-network optimization, reduces the number of trainable components, and leads to faster training.
It also provides flexibility in choosing approximation spaces: finite element bases can be used when a mesh and boundary-conforming discretization are available, random-feature bases are suitable for mesh-free and high-dimensional settings, and radial basis functions are useful for scattered-data and geometric-field representations.
Table~\ref{tab:method-comparison} summarizes the position of FB-C2CNet relative to representative neural operator architectures in terms of basis construction, trainable representation modules, and sampling flexibility.

\begin{table}[!ht]
\centering
\caption{
Comparison of FB-C2CNet with representative neural operator architectures considered in this work. 
The comparison focuses on how input and output functions are represented, how the representation is constructed, and which components are trained.
}
\label{tab:method-comparison}
\renewcommand{\arraystretch}{1.18}
\scriptsize

\begin{tabularx}{\textwidth}{@{}
p{0.15\textwidth}
>{\raggedright\arraybackslash}X
>{\raggedright\arraybackslash}X
>{\raggedright\arraybackslash}X
>{\raggedright\arraybackslash}p{0.14\textwidth}
>{\raggedright\arraybackslash}X@{}}
\toprule
\textbf{Method}
& \textbf{Input representation}
& \textbf{Output representation}
& \textbf{Representation source}
& \textbf{Input/Output sampling}
& \textbf{Main trainable part} \\
\midrule

DeepONet~\cite{lu2021deeponet}
& Point values at sensors
& Coordinate-wise output
& Learned trunk representation
& Usually fixed sensors
& Branch and trunk networks \\

\midrule

FNO~\cite{li2021fourier}
& Grid values
& Grid values
& Fourier transform with learned lifting and projection
& Structured grids
& Neural operator layers \\

\midrule

PCA-Net~\cite{bhattacharya2021model}
& PCA coefficients
& PCA coefficients with PCA reconstruction
& Data-driven PCA basis from aligned snapshots
& Aligned discretization
& Coefficient map  (1 NN) \\

\midrule

BasisONet~\cite{hua2023basis}
& Basis coefficients
& Coordinate-wise output through neural basis
& Learned neural basis / trunk representation
& Flexible after coefficient encoding
& Branch network and neural basis (3 NNs) \\

\midrule

B2B~\cite{ingebrand2025basis}
& Function-encoder coefficients
& Basis coefficients
& Trained function-encoder/ decoder representation
& Flexible
& Function encoders and coefficient map (3 NNs) \\

\midrule

\textbf{FB-C2CNet}
& \textbf{Prescribed-basis coefficients}
& \textbf{Prescribed-basis coefficients}
& \textbf{FEM / RFM / RBF bases fixed before training}
& \textbf{Flexible via projection}
& \textbf{Coefficient map (1 NN)} \\

\bottomrule
\end{tabularx}
\end{table}

\subsection{Contributions}

In this work, we propose the \emph{Fixed-Basis Coefficient-to-Coefficient Operator Network} (FB-C2CNet), 
an operator-learning framework built on prescribed, data-independent approximation spaces. 
The main contributions are summarized as follows:

\begin{itemize}
    \item \textbf{A prescribed-basis coefficient-to-coefficient neural operator.}
    We propose FB-C2CNet, which learns PDE solution operators by mapping between coefficients in prescribed input and output approximation spaces. 
    Unlike POD/PCA-based methods that construct data-dependent bases from aligned snapshots, and unlike learned-basis neural operators that train representation modules, FB-C2CNet fixes the approximation spaces before training and learns only the coefficient-to-coefficient map.

    \item \textbf{Efficient training through fixed prescribed-basis encoding and decoding.}
By separating basis selection from neural-network optimization, FB-C2CNet reduces the dimension and complexity of the trainable operator. 
The framework can directly use a range of prescribed approximation spaces, including finite element bases determined by the computational geometry, mesh-free random-feature bases for scattered or high-dimensional samples, and radial-basis-function bases for geometric fields such as signed distance functions.

    \item \textbf{Stability and approximation analysis.}
    We analyze how coefficient recovery affects the stability and accuracy of FB-C2CNet. 
    In particular, we show how ridge regularization or SVD cut-off stabilizes the input coefficient representation while introducing an approximation bias. 
    We also identify the output projection error associated with the prescribed output space as an intrinsic accuracy floor for any fixed-basis prediction.

    \item \textbf{Numerical validation on regular, irregular, and inverse PDE problems.}
    We validate FB-C2CNet on a broad set of PDE operator-learning tasks, including elliptic equations, nonlinear time-dependent equations, weak-solution learning, high-dimensional problems, multi-input/output problems, and an inverse Stokes boundary-recovery problem. 
    The experiments demonstrate competitive accuracy with reduced training time and highlight the method's ability to handle non-aligned input observations and flexible output evaluation.
\end{itemize}

The remainder of the paper is organized as follows. 
Section~2 presents the proposed fixed-basis coefficient-to-coefficient operator-learning framework. 
Section~3 discusses basis selection, coefficient recovery, and the associated stability and approximation properties. 
Section~4 reports numerical experiments demonstrating the accuracy, efficiency, and flexibility of the proposed method. 
Section~5 concludes the paper and discusses future directions.

\section{Methodology}

We aim to approximate a PDE solution operator $G:X\to Y$, $f\mapsto u=G(f)$, where $X$ and $Y$ are function spaces defined on bounded domains.
Instead of representing input and output functions by dense pointwise values, FB-C2CNet represents them by coefficients in prescribed finite-dimensional approximation spaces.
The input function is first encoded into coefficients with respect to a fixed input basis.
A neural network then learns a map from input coefficients to output coefficients, and the output field is reconstructed by a fixed output basis.
Thus, the only trainable component is the coefficient-to-coefficient map, while the basis functions, encoder rule, and decoder are fixed before training.

\subsection{Problem formulation}

Suppose we are given $M$ input-output pairs $\{(f^{(k)},u^{(k)})\}_{k=1}^{M}\subset X\times Y$, with $u^{(k)}=G(f^{(k)})$.
The input function $f^{(k)}$ is observed at points $\{\bm x_i^{(k)}\}_{i=1}^{n_1^{(k)}}$, and the output function $u^{(k)}$ is evaluated at target points $\{\bm y_i^{(k)}\}_{i=1}^{n_2^{(k)}}$.
For simplicity, when all samples share the same input and output grids, we write these points as $\{\bm x_i\}_{i=1}^{n_1}$ and $\{\bm y_i\}_{i=1}^{n_2}$.
The corresponding sampled vectors are
\[
\bm f^{(k)}
=
\big(f^{(k)}(\bm x_1),\ldots,f^{(k)}(\bm x_{n_1})\big)^\top
\in\mathbb R^{n_1},
\qquad
\bm u^{(k)}
=
\big(u^{(k)}(\bm y_1),\ldots,u^{(k)}(\bm y_{n_2})\big)^\top
\in\mathbb R^{n_2}.
\]
The goal is to train a neural operator that predicts $\bm u^{(k)}$ for new input observations.

\subsection{Fixed-Basis Coefficient-to-Coefficient Operator Network}

FB-C2CNet starts by choosing two prescribed approximation spaces
\[
V_{\mathrm{in}}^{m_1}
=
\operatorname{span}\{\varphi_j\}_{j=1}^{m_1},
\qquad
V_{\mathrm{out}}^{m_2}
=
\operatorname{span}\{\psi_j\}_{j=1}^{m_2}.
\]
The input space $V_{\mathrm{in}}^{m_1}$ is used to encode input functions, while the output space $V_{\mathrm{out}}^{m_2}$ is used to decode predicted solution fields.
These bases are selected before training and are not updated during neural-network optimization.

Given an input function $f$, we approximate it by
\[
f(\bm x)\approx \sum_{j=1}^{m_1}a_j\varphi_j(\bm x),
\qquad
\bm a=(a_1,\ldots,a_{m_1})^\top\in\mathbb R^{m_1}.
\]
Let $\Phi\in\mathbb R^{n_1\times m_1}$ be the input basis evaluation matrix, with $\Phi_{ij}=\varphi_j(\bm x_i)$.
The input coefficient vector is obtained by regularized least-squares fitting:
\[
\bm a_\lambda(f)
=
\arg\min_{\bm c\in\mathbb R^{m_1}}
\frac{1}{n_1}\|\Phi\bm c-\bm f\|_2^2
+
\lambda\|\bm c\|_2^2,
\qquad
\bm a_\lambda(f)
=
(\Phi^\top\Phi+n_1\lambda I)^{-1}\Phi^\top\bm f.
\]
If the input observation points are sample-dependent, the same formula is applied with the sample-specific matrix $\Phi^{(k)}_{ij}=\varphi_j(\bm x_i^{(k)})$.
In this case, the basis functions and the coefficient-recovery rule are still fixed; only the basis evaluation matrix changes with the observed locations.

After encoding, FB-C2CNet trains a neural network $N_\theta:\mathbb R^{m_1}\to\mathbb R^{m_2}$, $\bm a\mapsto \bm b$, where $\bm b=(b_1,\ldots,b_{m_2})^\top$ is the predicted output coefficient vector.
The corresponding output field is reconstructed by the fixed decoder
\[
\widehat u_\theta(\bm y)
=
\sum_{j=1}^{m_2}
\big[N_\theta(\bm a_\lambda(f))\big]_j\psi_j(\bm y).
\]
At the target points $\{\bm y_i\}_{i=1}^{n_2}$, define the output basis evaluation matrix $\Psi\in\mathbb R^{n_2\times m_2}$, with $\Psi_{ij}=\psi_j(\bm y_i)$.
Then the predicted output vector is
\[
\widehat{\bm u}_\theta
=
\Psi N_\theta(\bm a_\lambda(f)).
\]
Thus, the learnable part of FB-C2CNet is only the finite-dimensional coefficient map $N_\theta$.
The input basis, output basis, coefficient-recovery procedure, and output decoder are non-trainable.

\subsection{Training objective}

In this work, we train FB-C2CNet using a physical-space relative $L^2$ loss.
For a training pair $(f^{(k)},u^{(k)})$, let $\bm a^{(k)}=\bm a_\lambda(f^{(k)})$ and $\widehat{\bm u}^{(k)}_\theta=\Psi N_\theta(\bm a^{(k)})$.
The empirical loss is
\[
\mathcal L(\theta)
=
\frac{1}{M}
\sum_{k=1}^{M}
\frac{
\left\|
\Psi N_\theta(\bm a^{(k)})-\bm u^{(k)}
\right\|_2
}{
\left\|\bm u^{(k)}\right\|_2
}.
\]
The relative loss normalizes the prediction error by the magnitude of each output sample and prevents large-amplitude solutions from dominating the optimization.
Although output coefficients can also be computed by projecting $u^{(k)}$ onto $V_{\mathrm{out}}^{m_2}$, we use the physical-space loss because it directly measures the error of the reconstructed solution field and avoids over-penalizing coefficient non-uniqueness when the output basis is non-orthogonal.

The same formulation also applies to sample-dependent output locations.
In that case, one replaces $\Psi$ by the sample-specific output evaluation matrix $\Psi^{(k)}_{ij}=\psi_j(\bm y_i^{(k)})$ and uses
\[
\widehat{\bm u}^{(k)}_\theta
=
\Psi^{(k)}N_\theta(\bm a^{(k)}).
\]
This is one reason why FB-C2CNet can naturally evaluate outputs at arbitrary target locations.

\subsection{Extension to multi-component inputs and outputs}

The above presentation focuses on a single-input/single-output operator for clarity.
For multi-component inputs, such as $f=(f_1,\ldots,f_p)$, one may either encode each component separately and concatenate the resulting coefficient vectors, or use a vector-valued input basis to encode all components jointly.
Similarly, for multi-component outputs, one may use separate scalar-valued output bases for each component or a shared vector-valued output basis.
After these coefficient representations are constructed, the learning problem has the same form $N_\theta:\mathbb R^{m_1}\to\mathbb R^{m_2}$, $\bm a\mapsto \bm b$.
This allows FB-C2CNet to handle scalar-valued, vector-valued, and multi-input/output PDE solution operators within the same coefficient-to-coefficient framework.

The overall workflow of FB-C2CNet is illustrated in Figure~\ref{fig:workflow_FB-C2CNet}.
A summary of the main notation is provided in Table~\ref{tab:notation}.

\begin{table}[!ht]
\centering
\caption{Summary of notation and key quantities used in FB-C2CNet.}
\label{tab:notation}
\renewcommand{\arraystretch}{1.12}
\small
\begin{tabular}{@{}clc@{}}
\toprule
\textbf{Category} & \textbf{Description} & \textbf{Notation} \\
\midrule
\multirow{4}{*}{Data}
& Number of training samples & $M$ or $N_{\mathrm{train}}$ \\
& Number of testing samples & $N_{\mathrm{test}}$ \\
& Number of input observation points & $n_1$ \\
& Number of output evaluation points & $n_2$ \\
\midrule
\multirow{4}{*}{Basis representation}
& Input approximation space & $V_{\mathrm{in}}^{m_1}$ \\
& Output approximation space & $V_{\mathrm{out}}^{m_2}$ \\
& Number of input basis functions & $m_1$ \\
& Number of output basis functions & $m_2$ \\
\midrule
\multirow{4}{*}{Matrices and maps}
& Input basis evaluation matrix & $\Phi$ \\
& Output basis evaluation matrix & $\Psi$ \\
& Encoder regularization parameter & $\lambda$ \\
& Coefficient neural network & $N_\theta$ \\
\bottomrule
\end{tabular}
\end{table}

\newpage
\begin{algorithm}[!ht]
\renewcommand{\thealgorithm}{}   % 去掉算法编号
\algrenewcommand\Require{\textbf{Input:}}
\algrenewcommand\Ensure{\textbf{Output:}}
\centering
\footnotesize
\caption{Fixed-Basis Coefficient-to-Coefficient Operator Learning}
\label{alg:coef2coef}

% -------------------- Phase I --------------------
\begin{algorithmic}
  \State {\color{blue}\textbf{Phase I}: \textbf{Choose fixed input/output basis systems:} 
    $V_{\mathrm{in}}=\text{span}\{\phi_j(\vx)\}_{j=1}^{m_1}$,\;
    $V_{\mathrm{out}}=\text{span}\{\psi_j(\vy)\}_{j=1}^{m_2}$.}
\end{algorithmic}

\vspace{3pt}

\noindent
\begin{minipage}{0.48\linewidth}
\begin{tcolorbox}
\textbf{Strategy 1: Random feature model}\\[2pt]
\textbf{Input:} $m_1,m_2$, $R_m^{\mathrm{in/out}}$, 
$M_p^{\mathrm{in/out}},J_n^{\mathrm{in/out}}$\\[3pt]
\textbf{Procedure:}
\begin{itemize}
  \item For $n=1,\dots,M_p^{\mathrm{in}}$ and $j=1,\dots,J_n^{\mathrm{in}}$:
  \[
  \phi_{nj}(\vx)=\omega_n(\vx)\sigma(\vk^{\mathrm{in}}_{nj}\!\cdot\!\vx+b^{\mathrm{in}}_{nj}),
  \]
  with $\vk^{\mathrm{in}}_{nj}\!\sim\!U([-R_m^{\mathrm{in}},R_m^{\mathrm{in}}]^d)$,
  $b^{\mathrm{in}}_{nj}\!\sim\!U([-R_m^{\mathrm{in}},R_m^{\mathrm{in}}])$.
  \item Similarly construct $\psi_{nj}(\vy)$ for output basis with
  $\vk^{\mathrm{out}}_{nj}$, $b^{\mathrm{out}}_{nj}$.
\end{itemize}
\textbf{Output:} 
$V_{\mathrm{in}}=\text{span}\{\phi_j(\vx)\}_{j=1}^{m_1}$,
$V_{\mathrm{out}}=\text{span}\{\psi_j(\vy)\}_{j=1}^{m_2}$.
\end{tcolorbox}
\end{minipage}\hfill
\begin{minipage}{0.48\linewidth}
\begin{tcolorbox}
\textbf{Strategy 2: Finite element model}\\[2pt]
\textbf{Input:} element types $E^{\mathrm{in/out}}$, polynomial degrees 
$p^{\mathrm{in/out}}$, evaluation points 
$\boldsymbol{\xi}^{\mathrm{in/out}}$\\[3pt]
\textbf{Procedure:}
\begin{itemize}
  \item Generate the mesh $\mathcal{T}_h$; then, for the chosen element family and degree, construct the global FE basis.
\end{itemize}
\textbf{Output:} 
$V_{\mathrm{in}}=\text{span}\{\phi_j(\vx)\}_{j=1}^{m_1}$,
$V_{\mathrm{out}}=\text{span}\{\psi_j(\vy)\}_{j=1}^{m_2}$.
\end{tcolorbox}
\end{minipage}

\vspace{6pt}
% -------------------- Phase II --------------------
\begin{algorithmic}
  \State {\color{blue}\textbf{Phase II}: \textbf{Encode input coefficients via least squares.}} \\
  \Require data sets $\{\mathbf{f}^{(k)}\}_{k=1}^M$, each $\mathbf{f}^{(k)}=\{f^{(k)}(\vx_i)\}_{i=1}^{n_1}$; penalty $\lambda$
  \For{$k=1$ \textbf{to} $M$}
    \State compute $\va^{(k)}\in\mathbb{R}^{m_1}$ by
    \Statex \hspace{\algorithmicindent}$\displaystyle
      \va^{(k)} = \underset{\va\in\mathbb{R}^{m_1}}{\arg\min}\;
      \frac{1}{n_1}\sum_{i=1}^{n_1}
      \left| f^{(k)}(\vx_i)-\sum_{j=1}^{m_1}a_j\phi_j(\vx_i)\right|^2
      +\lambda\|\va\|_2^2.$
  \EndFor  \\
  \Ensure $\{\va^{(k)}\}_{k=1}^M$

  \vspace{3pt}
  \State {\color{blue}\textbf{Phase III}: \textbf{Train coefficient-to-coefficient network $NN_\theta$.}} \\
  \Require targets $\{\mathbf{u}^{(k)}\}_{k=1}^M$, each $\mathbf{u}^{(k)}=\{u^{(k)}(\vy_i)\}_{i=1}^{n_2}$; step size $\eta$
  \State form $\Psi\in\mathbb{R}^{n_2\times m_2}$ with $\Psi_{ij}=\psi_j(\vy_i)$
  \While{not converged}
    \State for each $k$:
    \Statex \hspace{\algorithmicindent}
      $\vb^{(k)}=NN_\theta(\va^{(k)})$, 
      $\hat{\vu_{\bm \theta}}^{(k)}=\Psi\vb^{(k)}$, 
       $\vu^{(k)}=\big(u^{(k)}(\vy_i)\big)_{i=1}^{n_2}$
    \State loss:
    \Statex \hspace{\algorithmicindent}$\displaystyle
      L(\bm \theta)=\frac{1}{M}\sum_{k=1}^{M}\frac{\|\hat{\vu_{\bm \theta}}^{(k)}-\vu^{(k)}\|_2}
           {\|\vu^{(k)}\|_2}$
    \State update: $\bm \theta \leftarrow \bm\theta - \eta\nabla_{\bm \theta}L(\bm \theta)$
  \EndWhile
  \State \textbf{Output:} $N\!N_{\bm \theta}$
\end{algorithmic}
\end{algorithm}

\section{Prescribed Basis Spaces, Coefficient Recovery, and Error Analysis}
\label{sec:basis-coefficient-analysis}

This section analyzes three principal error sources in FB-C2CNet: the approximation errors of the prescribed input and output spaces, the stability and bias of coefficient recovery from sampled observations, and the learning error of the coefficient-to-coefficient map. 
We use the notation introduced in Section~2, with
$V_{\mathrm{in}}^{m_1}=\operatorname{span}\{\varphi_j\}_{j=1}^{m_1}$,
$V_{\mathrm{out}}^{m_2}=\operatorname{span}\{\psi_j\}_{j=1}^{m_2}$, and
$N_\theta:\mathbb R^{m_1}\to\mathbb R^{m_2}$.

\subsection{Choice of prescribed basis spaces}
\label{subsec:prescribed-basis-spaces}

A key feature of FB-C2CNet is that the input and output approximation spaces are selected before neural-network training. 
The basis construction is therefore a numerical approximation choice rather than part of the optimization of $N_\theta$. 
The basis spaces may be selected according to the domain geometry, boundary conditions, regularity, dimension, and sampling pattern of the problem.

In this work, we mainly use random feature model (RFM), finite element method (FEM), and radial basis function (RBF) bases. 
RFM bases provide mesh-free representations that are convenient for scattered observations and high-dimensional problems. 
FEM bases are naturally adapted to meshed domains, boundary conditions, and existing numerical discretizations. 
RBF bases are useful for scattered-data approximation and geometric fields such as signed distance functions.

\paragraph{Random feature model bases.}

The RFM constructs prescribed basis functions from randomly generated nonlinear features. 
We mainly use the localized partition-of-unity form introduced in~\cite{chen2022bridging}. 
Let $\{\omega_n\}_{n=1}^{M_p}$ be localized partition functions. 
The corresponding basis functions are
\begin{equation}
\varphi_{n,j}(\mathbf x)
=
\omega_n(\mathbf x)
\sigma\!\left(
\mathbf k_{n,j}\cdot\mathbf x+b_{n,j}
\right),
\qquad
n=1,\ldots,M_p,
\quad
j=1,\ldots,J_n,
\label{eq:pou-rfm-basis-function}
\end{equation}
where $\mathbf k_{n,j}\sim\mathcal U([-R_m,R_m]^d)$ and
$b_{n,j}\sim\mathcal U([-R_m,R_m])$ are sampled once and then fixed. 
The activation function $\sigma$ may be chosen as $\tanh$, $\sin$, or another smooth nonlinear function.

The approximation takes the form
\begin{equation}
u_m(\mathbf x)
=
\sum_{n=1}^{M_p}
\sum_{j=1}^{J_n}
a_{n,j}\varphi_{n,j}(\mathbf x).
\label{eq:pou-rfm-approximation}
\end{equation}
The global RFM is recovered by taking $M_p=1$ and $\omega_1\equiv1$. 
Once the random features are fixed, only the outer coefficients $\{a_{n,j}\}$ are estimated. 
The one-dimensional partition function used in the experiments is given in Appendix~\ref{app:pou-window}.

RFM bases can be evaluated directly at scattered or sample-dependent locations, which makes them convenient for irregular observations and high-dimensional problems~\cite{nelsen2024operator,chi2024random,chen2023random}. 
However, random features may be strongly correlated, and the resulting basis evaluation matrix may contain small singular values. 
Stable coefficient recovery may therefore require ridge regularization or a truncated-SVD procedure.

\paragraph{Finite element bases.}

Given a mesh $\mathcal T_h$ of a domain $\Omega\subset\mathbb R^d$, FEM basis functions are constructed locally on the mesh elements~\cite{brenner2008mathematical,ciarlet2002finite}. 
For linear $P_1$ elements, the nodal basis functions satisfy
\begin{equation}
\varphi_j(\mathbf x_i)=\delta_{ij},
\qquad
1\leq i,j\leq N_h,
\label{eq:fem-hat-basis}
\end{equation}
where $\{\mathbf x_i\}_{i=1}^{N_h}$ are the mesh nodes. 
A finite element approximation is written as
\begin{equation}
u_h(\mathbf x)
=
\sum_{j=1}^{N_h}a_j\varphi_j(\mathbf x).
\label{eq:fem-approximation}
\end{equation}
Higher-order $P_k$ spaces can be used analogously.

FEM bases are local, conforming, and supported by classical approximation theory. 
They provide a direct interface with existing meshes, degrees of freedom, and numerical solver infrastructure. 
Their main limitations are the need for mesh generation and the rapid growth of the number of degrees of freedom in high-dimensional domains.

\paragraph{Other prescribed bases}

Fourier bases are natural for smooth periodic problems, polynomial bases are useful for globally smooth functions on simple domains, and RBFs provide mesh-free representations for scattered data and geometric fields~\cite{boyd2001chebyshev,trefethen2000spectral,fasshauer2007meshfree,buhmann2003radial}. 
The choice of basis should therefore be regarded as a problem-dependent numerical approximation decision.

\begin{table}[!ht]
\centering
\caption{Prescribed basis systems used in FB-C2CNet.}
\label{tab:basis-comparison-practical}
\renewcommand{\arraystretch}{1.2}
\small
\begin{tabularx}{\textwidth}{@{}l l X X@{}}
\toprule
\textbf{Basis}
& \textbf{Typical form}
& \textbf{Main strength}
& \textbf{Main limitation} \\
\midrule

FEM
& $\varphi_j|_K\in\mathbb P_p(K)$
& Geometry- and boundary-conforming; compatible with existing meshes
& Requires mesh generation and scales poorly with dimension \\

RFM
& $\omega_n(\mathbf x)\sigma(\mathbf w_{n,j}\cdot\mathbf x+b_{n,j})$
& Mesh-free; directly evaluated at scattered locations
& Correlated features may lead to ill-conditioned coefficient recovery \\

RBF
& $\rho(\|\mathbf x-\mathbf c_j\|/\varepsilon_j)$
& Flexible for scattered observations and geometric fields
& Accuracy depends on center placement and width selection \\

\bottomrule
\end{tabularx}
\end{table}

\begin{figure}[!ht]
    \centering
    \begin{minipage}[t]{0.46\textwidth}
        \centering
        \includegraphics[width=\linewidth]{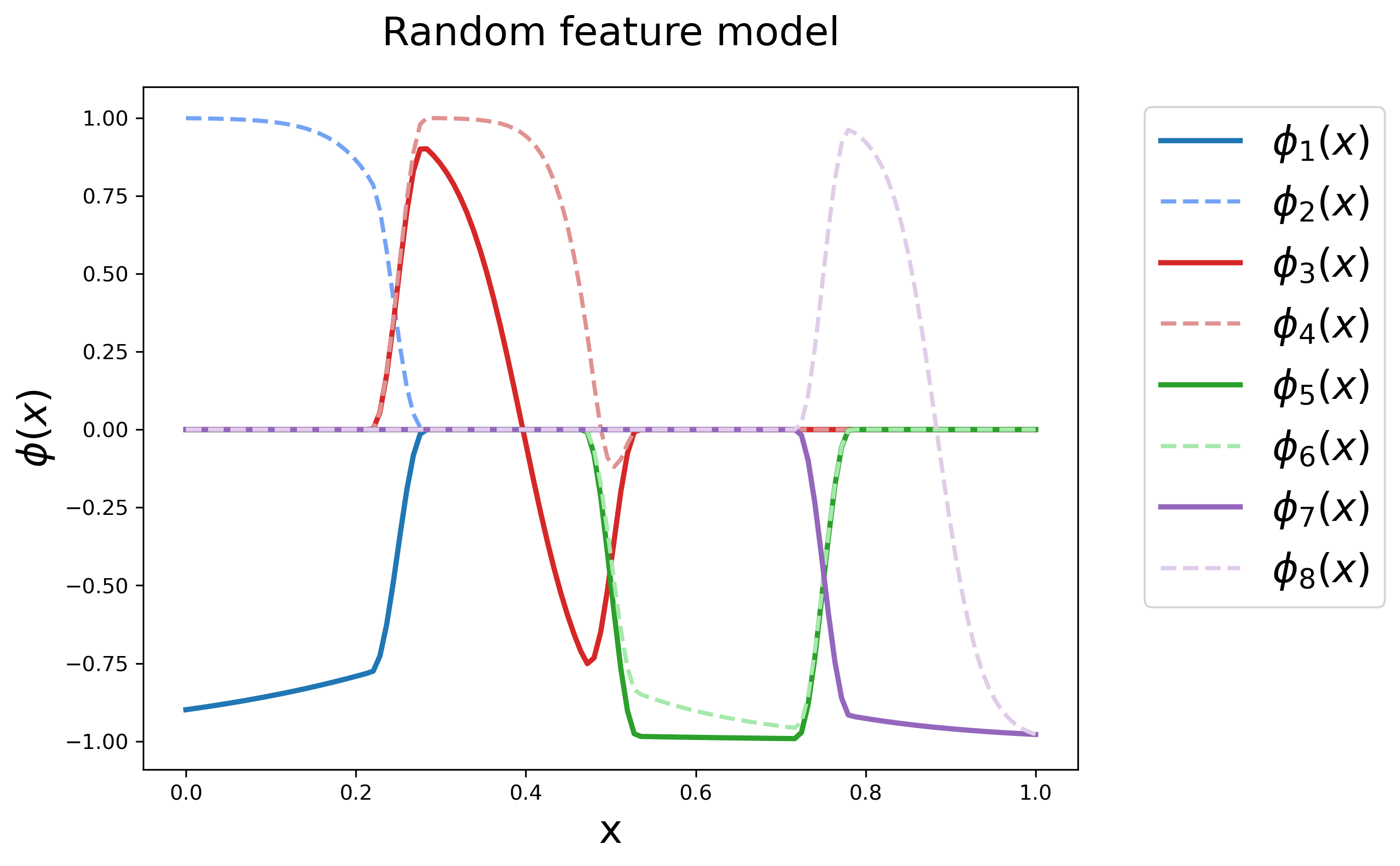}
        \caption*{(a) RFM basis.}
    \end{minipage}
    \hspace{0.035\textwidth}
    \begin{minipage}[t]{0.46\textwidth}
        \centering
        \includegraphics[width=\linewidth]{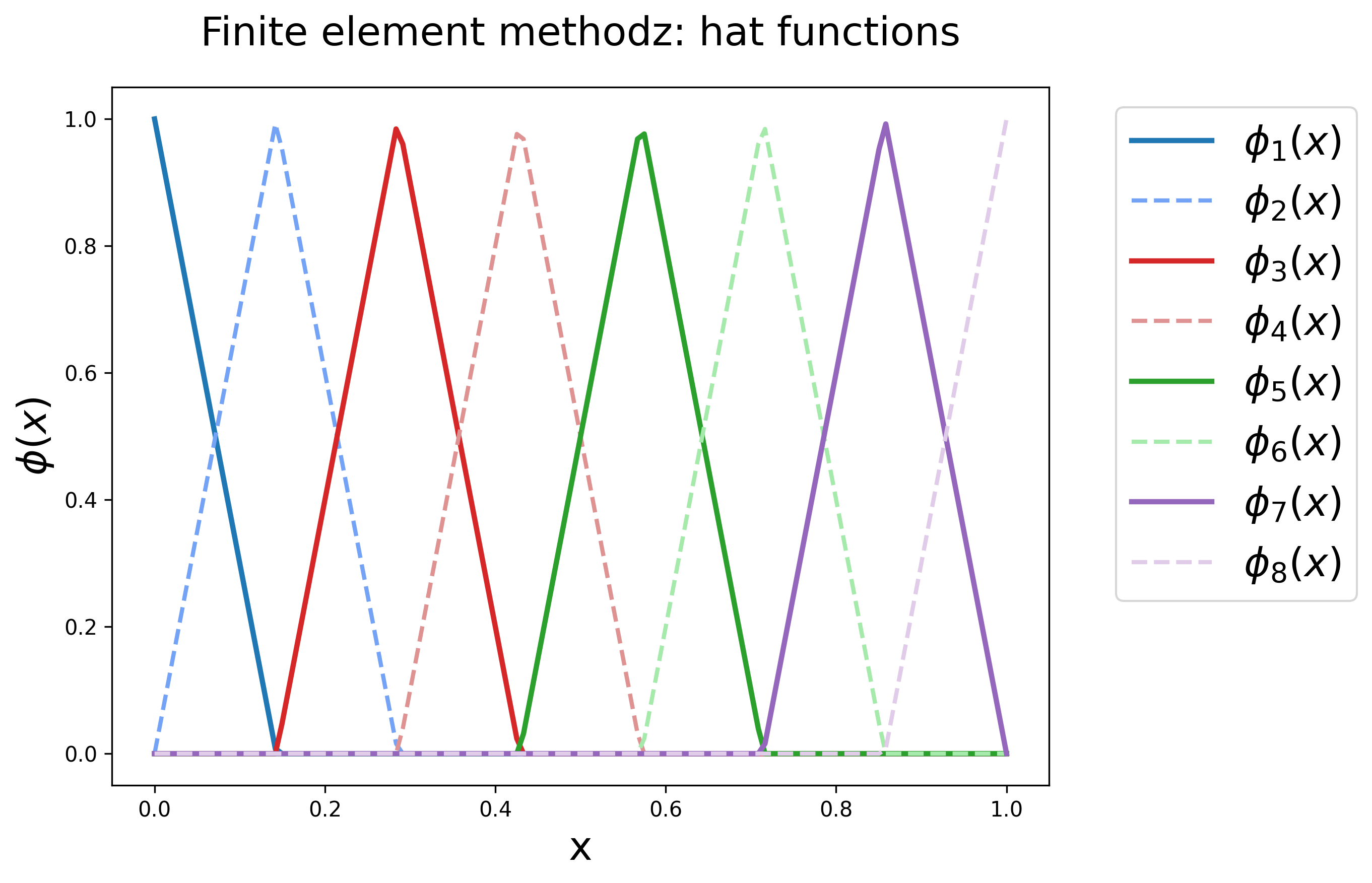}
        \caption*{(b) FEM basis.}
    \end{minipage}
    \caption{
    One-dimensional prescribed basis systems:
    (a) a localized RFM basis with $M_p=4$ partitions and $J_n=2$ features per partition;
    (b) a $P_1$ FEM basis with $N_h=8$ nodal functions.
    }
    \label{fig:1d_basis}
\end{figure}

\subsection{Stability and bias of coefficient recovery}
\label{subsec:coefficient-recovery-analysis}

The input coefficients are obtained by fitting sampled function values in the prescribed input space. 
Even when FEM bases are used, we formulate this encoding step as a least-squares problem rather than assuming that the observations coincide with the FEM nodes. 
The stability of the resulting coefficient representation is therefore determined by the basis evaluation matrix.

For mesh-free or global bases, such as RFM, RBF, and polynomial bases, the evaluation matrix may contain small singular values. 
Regularization or singular-value truncation may then be needed to obtain stable coefficients. 
The same issue can arise for FEM bases when the observations are scattered or when an overdetermined fitting problem is used.

\subsubsection{Regularized coefficient estimation}
\label{subsec:regularized-coefficient-estimation}

Let $\mathbf x_1,\ldots,\mathbf x_{n_1}$ be the input observation points and let
$\mathbf f=(f(\mathbf x_1),\ldots,f(\mathbf x_{n_1}))^\top\in\mathbb R^{n_1}$. 
For $V_{\mathrm{in}}^{m_1}=\operatorname{span}\{\varphi_1,\ldots,\varphi_{m_1}\}$, define
$\Phi\in\mathbb R^{n_1\times m_1}$ by $\Phi_{ij}=\varphi_j(\mathbf x_i)$.

The objective in Section~2 contains the normalization factor $1/n_1$. 
Multiplying that objective by $n_1$ does not change its minimizer. 
For notational convenience, throughout this section we absorb the factor $n_1$ into the regularization parameter and continue to denote the rescaled parameter by $\lambda$. 
The coefficient encoder is therefore written as
\begin{equation}
\mathbf a_\lambda(f)
=
\arg\min_{\mathbf a\in\mathbb R^{m_1}}
\left\{
\|\Phi\mathbf a-\mathbf f\|_2^2
+
\lambda\|\mathbf a\|_2^2
\right\},
\qquad
\lambda>0.
\label{eq:ridge-coefficient-problem}
\end{equation}
Equivalently,
\begin{equation}
\mathbf a_\lambda(f)
=
E_\lambda\mathbf f,
\qquad
E_\lambda
=
(\Phi^\top\Phi+\lambda I)^{-1}\Phi^\top.
\label{eq:ridge-encoder-definition}
\end{equation}

\begin{lemma}[Stability of the regularized coefficient encoder]
\label{lem:regularized-encoder-stability}
For every $\lambda>0$,
\begin{equation}
\|E_\lambda\|_2
=
\max_i
\frac{\sigma_i(\Phi)}
{\sigma_i(\Phi)^2+\lambda}
\leq
\frac{1}{2\sqrt{\lambda}}.
\label{eq:encoder-stability-main}
\end{equation}
Consequently,
\begin{equation}
\|\mathbf a_\lambda(f)-\mathbf a_\lambda(g)\|_2
\leq
\frac{1}{2\sqrt{\lambda}}
\|\mathbf f-\mathbf g\|_2.
\label{eq:encoder-lipschitz}
\end{equation}
\end{lemma}

Lemma~\ref{lem:regularized-encoder-stability} shows that $\lambda$ controls the amplification of perturbations in the sampled input. 
In singular-value coordinates, the ridge encoder applies the filter
$\sigma_i\mapsto\sigma_i/(\sigma_i^2+\lambda)$, which damps poorly conditioned directions.

\begin{remark}[Sample-dependent observations]
\label{rem:sample-dependent-observations}
If the observation locations vary between samples, the evaluation matrix is replaced by a sample-dependent matrix $\Phi^{(k)}$. 
Lemma~\ref{lem:regularized-encoder-stability} applies to each sample separately. 
Moreover, the upper bound $1/(2\sqrt{\lambda})$ is independent of the particular singular values of $\Phi^{(k)}$, and is therefore uniform with respect to the sampling locations.
\end{remark}

Mapping the coefficients back to the sampled input space gives
\begin{equation}
\widetilde{\mathbf f}_\lambda
=
H_\lambda\mathbf f,
\qquad
H_\lambda
=
\Phi(\Phi^\top\Phi+\lambda I)^{-1}\Phi^\top.
\label{eq:ridge-reconstruction-operator}
\end{equation}
The discrepancy $\|\mathbf f-\widetilde{\mathbf f}_\lambda\|_2$ contains both the component that cannot be represented in the prescribed input space and the shrinkage introduced by regularization.

\subsubsection{Encoding bias}
\label{subsec:encoding-bias-lambda}

\begin{proposition}[Spectral form of the encoding error]
\label{prop:encoding-bias-ridge}
Let $\Phi=U\Sigma V^\top$ be a singular value decomposition of $\Phi$. 
Let $\sigma_1,\ldots,\sigma_r>0$ be its nonzero singular values and $u_1,\ldots,u_r$ the corresponding left singular vectors. 
For $\mathbf f\in\mathbb R^{n_1}$, write
\begin{equation}
\mathbf f
=
\sum_{i=1}^r\alpha_i u_i+\mathbf f_\perp,
\qquad
\mathbf f_\perp\perp\operatorname{Range}(\Phi),
\label{eq:f-svd-decomposition}
\end{equation}
where $\alpha_i=\langle\mathbf f,u_i\rangle$. 
Then
\begin{equation}
\widetilde{\mathbf f}_\lambda
=
\sum_{i=1}^r
\frac{\sigma_i^2}{\sigma_i^2+\lambda}
\alpha_i u_i,
\label{eq:regularized-reconstruction-spectral}
\end{equation}
and
\begin{equation}
\|\mathbf f-\widetilde{\mathbf f}_\lambda\|_2^2
=
\sum_{i=1}^r
\left(
\frac{\lambda}{\sigma_i^2+\lambda}
\right)^2
|\alpha_i|^2
+
\|\mathbf f_\perp\|_2^2.
\label{eq:encoding-bias-spectral}
\end{equation}
\end{proposition}

The first term in \eqref{eq:encoding-bias-spectral} is the shrinkage error introduced by regularization within the representable subspace. 
The second term is the unrepresentable component of the sampled input and is independent of $\lambda$. 
Increasing $\lambda$ suppresses poorly conditioned directions and improves stability, but also increases the shrinkage of informative components. 
This is the stability--bias trade-off associated with coefficient recovery.

\subsection{Output-space approximation error}
\label{subsec:output-space-approximation}

We next separate the approximation error imposed by the prescribed output space from the error learned by the coefficient network. 
Assume that $Y$ is a Hilbert space, and let
$P_{\mathrm{out}}:Y\to V_{\mathrm{out}}^{m_2}$ be the orthogonal projection. 
Define
\begin{equation}
\varepsilon_{\mathrm{out}}(u)
=
\|u-P_{\mathrm{out}}u\|_Y
=
\inf_{v\in V_{\mathrm{out}}^{m_2}}
\|u-v\|_Y.
\label{eq:output-projection-error}
\end{equation}

\begin{proposition}[Output-space error decomposition]
\label{prop:output-space-lower-bound}
Let $\widehat G_{\theta,\lambda}(f)\in V_{\mathrm{out}}^{m_2}$ be an FB-C2CNet prediction. 
Then
\begin{equation}
\|G(f)-\widehat G_{\theta,\lambda}(f)\|_Y^2
=
\varepsilon_{\mathrm{out}}(G(f))^2
+
\|P_{\mathrm{out}}G(f)-\widehat G_{\theta,\lambda}(f)\|_Y^2.
\label{eq:output-projection-identity}
\end{equation}
In particular,
\begin{equation}
\varepsilon_{\mathrm{out}}(G(f))
\leq
\|G(f)-\widehat G_{\theta,\lambda}(f)\|_Y.
\label{eq:output-projection-lower-bound}
\end{equation}
\end{proposition}

The first term on the right-hand side of \eqref{eq:output-projection-identity} is the intrinsic approximation error of the prescribed output space. 
The second term is the learnable error associated with the coefficient map. 
Thus, no prediction in $V_{\mathrm{out}}^{m_2}$ can achieve an error below the output projection error.

Let $D_{\mathrm{out}}:\mathbb R^{m_2}\to V_{\mathrm{out}}^{m_2}$ denote the prescribed output decoder. 
Choose a coefficient vector $\mathbf b(f)$ satisfying
$P_{\mathrm{out}}G(f)=D_{\mathrm{out}}\mathbf b(f)$; if the representation is non-unique, $\mathbf b(f)$ denotes the minimum-norm coefficient vector. 
Define
$\widehat{\mathbf b}_{\theta,\lambda}(f)=N_\theta(\mathbf a_\lambda(f))$ and
$\widehat G_{\theta,\lambda}(f)=D_{\mathrm{out}}\widehat{\mathbf b}_{\theta,\lambda}(f)$.

For the theoretical analysis, we introduce the projection-based loss
\begin{equation}
\ell_{\theta,\lambda}(f)
=
\left\|
D_{\mathrm{out}}
\left(
\widehat{\mathbf b}_{\theta,\lambda}(f)-\mathbf b(f)
\right)
\right\|_Y^2
=
\|\widehat G_{\theta,\lambda}(f)-P_{\mathrm{out}}G(f)\|_Y^2.
\label{eq:lambda-loss-definition}
\end{equation}
For i.i.d. samples $f_1,\ldots,f_M$ drawn from a distribution $\mu$, define
\begin{equation}
\mathcal R_\lambda(\theta)
=
\mathbb E_{f\sim\mu}\ell_{\theta,\lambda}(f),
\qquad
\widehat{\mathcal R}_{M,\lambda}(\theta)
=
\frac1M\sum_{i=1}^M\ell_{\theta,\lambda}(f_i).
\label{eq:population-empirical-risk-lambda}
\end{equation}

\begin{remark}[Relation to the training loss]
\label{rem:projection-training-loss}
The loss in \eqref{eq:lambda-loss-definition} is an auxiliary loss introduced to isolate the learnable coefficient error. 
The numerical experiments use the relative physical-space loss defined in Section~2. 
An analogous estimate applies to the relative loss when $\|G(f)\|_Y$ is uniformly bounded away from zero.
\end{remark}

\subsection{Stability-dependent generalization estimate}
\label{subsec:stability-dependent-generalization}

The following estimate makes explicit how the stability of the coefficient encoder enters a standard statistical-learning bound. 
It is not intended as a sharp, architecture-specific generalization result.

\begin{assumption}[Boundedness and complexity]
\label{assump:generalization-main}
The following conditions hold:
\begin{enumerate}
    \item There exists $R_f>0$ such that $\|\mathbf f\|_2\leq R_f$ almost surely.
    \item There exists $B_b>0$ such that $\|\mathbf b(f)\|_2\leq B_b$ almost surely.
    \item The output decoder satisfies $\|D_{\mathrm{out}}\|_{\mathrm{op}}\leq C_D$.
    \item There exist $B_N,L_N>0$ such that
    \begin{equation}
    \sup_{\theta\in\Theta}
    \|N_\theta(\mathbf a)\|_2
    \leq
    B_N+L_N\|\mathbf a\|_2
    \qquad
    \text{for all }\mathbf a\in\mathbb R^{m_1}.
    \label{eq:network-affine-growth}
    \end{equation}
    \item The loss class
    $\mathcal F_\lambda=\{\ell_{\theta,\lambda}:\theta\in\Theta\}$
    satisfies
    $\sup_{\lambda\in\Lambda}\operatorname{Pdim}(\mathcal F_\lambda)\leq\Pi$.
\end{enumerate}
\end{assumption}

\begin{lemma}[Uniform loss bound]
\label{lem:lambda-loss-bound}
Under Assumption~\ref{assump:generalization-main}, define
\begin{equation}
B_\lambda
=
C_D
\left(
B_b+B_N+\frac{L_NR_f}{2\sqrt{\lambda}}
\right).
\label{eq:B-lambda-definition}
\end{equation}
Then $0\leq\ell_{\theta,\lambda}(f)\leq B_\lambda^2$ for every $\theta\in\Theta$ and almost every $f$.
\end{lemma}

\begin{theorem}[Stability-dependent generalization estimate]
\label{thm:lambda-dependent-generalization}
Suppose Assumption~\ref{assump:generalization-main} holds. 
Then, for $M\geq2$,
\begin{equation}
\mathbb E
\sup_{\theta\in\Theta}
\left|
\mathcal R_\lambda(\theta)
-
\widehat{\mathcal R}_{M,\lambda}(\theta)
\right|
\leq
C
B_\lambda^2
\sqrt{\frac{\Pi}{M}}
\log(eM),
\label{eq:lambda-dependent-generalization-bound}
\end{equation}
where $C>0$ is a universal constant.
\end{theorem}

The theorem preserves the standard $M^{-1/2}$ dependence on sample size, while encoder regularization enters the constant through $B_\lambda$. 
Larger values of $\lambda$ reduce the possible amplification of sampled inputs, but this effect must be balanced against the encoding error in Proposition~\ref{prop:encoding-bias-ridge}.

\begin{corollary}[Excess risk with optimization error]
\label{cor:erm-excess-risk-lambda}
Let $\widehat\theta_\lambda\in\Theta$ satisfy
\begin{equation}
\widehat{\mathcal R}_{M,\lambda}(\widehat\theta_\lambda)
\leq
\inf_{\theta\in\Theta}
\widehat{\mathcal R}_{M,\lambda}(\theta)
+
\eta_{\mathrm{opt}},
\label{eq:approximate-erm-condition}
\end{equation}
where $\eta_{\mathrm{opt}}\geq0$. 
Then
\begin{equation}
\mathbb E\mathcal R_\lambda(\widehat\theta_\lambda)
\leq
\inf_{\theta\in\Theta}
\mathcal R_\lambda(\theta)
+
C
B_\lambda^2
\sqrt{\frac{\Pi}{M}}
\log(eM)
+
\eta_{\mathrm{opt}},
\label{eq:erm-excess-risk-lambda}
\end{equation}
where $C$ denotes a possibly different universal constant.
\end{corollary}

Corollary~\ref{cor:erm-excess-risk-lambda} separates the coefficient-learning error into network approximation, statistical generalization, and optimization contributions.

\begin{remark}[Schematic total-error decomposition]
\label{rem:schematic-total-error}
Suppose that the sampled reconstruction $\widetilde{\mathbf f}_\lambda$ is identified with a function in the prescribed input space and that the solution operator is Lipschitz continuous with constant $L_G$ in the corresponding norm. 
The preceding results then motivate the informal decomposition
\begin{equation}
\begin{aligned}
\mathcal E_{\mathrm{test}}(\lambda)
\lesssim{}&
\underbrace{
\mathbb E\varepsilon_{\mathrm{out}}(G(f))^2
}_{\text{output approximation}}
+
\underbrace{
L_G^2
\mathbb E
\|\mathbf f-\widetilde{\mathbf f}_\lambda\|_2^2
}_{\text{input encoding}}
\\
&+
\underbrace{
\inf_{\theta\in\Theta}
\mathcal R_\lambda(\theta)
}_{\text{network approximation}}
+
\underbrace{
B_\lambda^2
\sqrt{\frac{\Pi}{M}}
\log(eM)
}_{\text{generalization}}
+
\underbrace{
\eta_{\mathrm{opt}}
}_{\text{optimization}}.
\end{aligned}
\label{eq:lambda-bias-generalization-decomposition}
\end{equation}
This expression is interpretive and is not used as an additional theorem.
\end{remark}

The analysis identifies two competing effects of the regularization parameter. 
Increasing $\lambda$ improves the stability of coefficient recovery and decreases the stability-dependent generalization constant, but it also increases the encoding bias by shrinking informative singular directions. 
A moderate regularization level is therefore preferable in practice.

The proofs of Lemma~\ref{lem:regularized-encoder-stability},
Proposition~\ref{prop:encoding-bias-ridge},
Proposition~\ref{prop:output-space-lower-bound},
Lemma~\ref{lem:lambda-loss-bound},
Theorem~\ref{thm:lambda-dependent-generalization}, and
Corollary~\ref{cor:erm-excess-risk-lambda}
are provided in~\ref{app:analysis-proofs}.

\subsection{Numerical illustration: spectral cut-off and coefficient stability}
\label{subsec:svd-cutoff-illustration}

We illustrate the preceding stability mechanism using the one-dimensional Darcy flow problem introduced in Section~\ref{subsubsec:1d-darcy-c2c-p2c}. 
This diagnostic experiment uses RFM bases for both input and output functions, with \(m_1=m_2=128\), and the functions are evaluated at 1000 uniformly spaced points. 
Instead of ridge regularization, the input coefficients are computed by truncated-SVD least squares with different thresholds \(\sigma_{\mathrm{cut}}\).

Let \(\Phi=U\Sigma V^\top\) be the singular value decomposition of the input basis evaluation matrix, with singular values \(\{\sigma_i\}\). 
The unregularized least-squares coefficient vector can be written formally as
\begin{equation}
\mathbf a_{\rm LS}
=
\sum_{\sigma_i>0}
\frac{\langle \mathbf f,u_i\rangle}{\sigma_i}v_i .
\label{eq:svd-ls-coeff}
\end{equation}
Small singular values may strongly amplify noise, sampling errors, or components of \(\mathbf f\) that are poorly represented by the basis. 
The truncated-SVD encoder keeps only singular directions above a prescribed threshold:
\begin{equation}
\mathbf a_{\rm cut}
=
E_{\rm cut}\mathbf f
=
\sum_{\sigma_i\geq \sigma_{\rm cut}}
\frac{\langle \mathbf f,u_i\rangle}{\sigma_i}v_i .
\label{eq:svd-cutoff-coeff}
\end{equation}
Equivalently, singular directions with \(\sigma_i<\sigma_{\rm cut}\) are discarded. 
Therefore,
\begin{equation}
\|E_{\rm cut}\|_2
\leq
\sigma_{\rm cut}^{-1}.
\label{eq:svd-cutoff-stability}
\end{equation}

The connection with ridge regularization is through the singular-value filters. 
For the ridge encoder,
$
E_\lambda
=
(\Phi^\top\Phi+\lambda I)^{-1}\Phi^\top,
$
the coefficient filter is
$
g_\lambda(\sigma)
=
\frac{\sigma}{\sigma^2+\lambda},
$
whereas the truncated-SVD coefficient filter is
$
g_{\rm cut}(\sigma)
=
\begin{cases}
\sigma^{-1}, & \sigma\geq\sigma_{\rm cut},\\
0, & \sigma<\sigma_{\rm cut}.
\end{cases}
$
Thus, ridge regularization smoothly damps small singular directions, while truncated SVD removes them sharply.

There is no exact one-to-one correspondence between \(\lambda\) and \(\sigma_{\rm cut}\), but they play analogous roles. 
From the reconstruction viewpoint, the ridge reconstruction filter is
$
r_\lambda(\sigma)
=
\frac{\sigma^2}{\sigma^2+\lambda},
$
which equals \(1/2\) at \(\sigma=\sqrt{\lambda}\). 
Therefore, \(\sqrt{\lambda}\) can be viewed as the characteristic singular-value scale below which ridge regularization strongly suppresses information. 
In this sense, the hard cut-off threshold \(\sigma_{\rm cut}\) is roughly analogous to \(\sqrt{\lambda}\). 
From the stability viewpoint, ridge satisfies \(\|E_\lambda\|_2\leq 1/(2\sqrt{\lambda})\), whereas truncated SVD satisfies \(\|E_{\rm cut}\|_2\leq 1/\sigma_{\rm cut}\). 
Matching these worst-case amplification factors gives the alternative correspondence \(\sigma_{\rm cut}\approx 2\sqrt{\lambda}\). 
In either interpretation, decreasing \(\sigma_{\rm cut}\) corresponds to weaker stabilization, similar to decreasing \(\lambda\).

\begin{figure}[!ht]
    \centering
    \begin{minipage}[t]{0.42\textwidth}
        \centering
        \includegraphics[width=\linewidth]{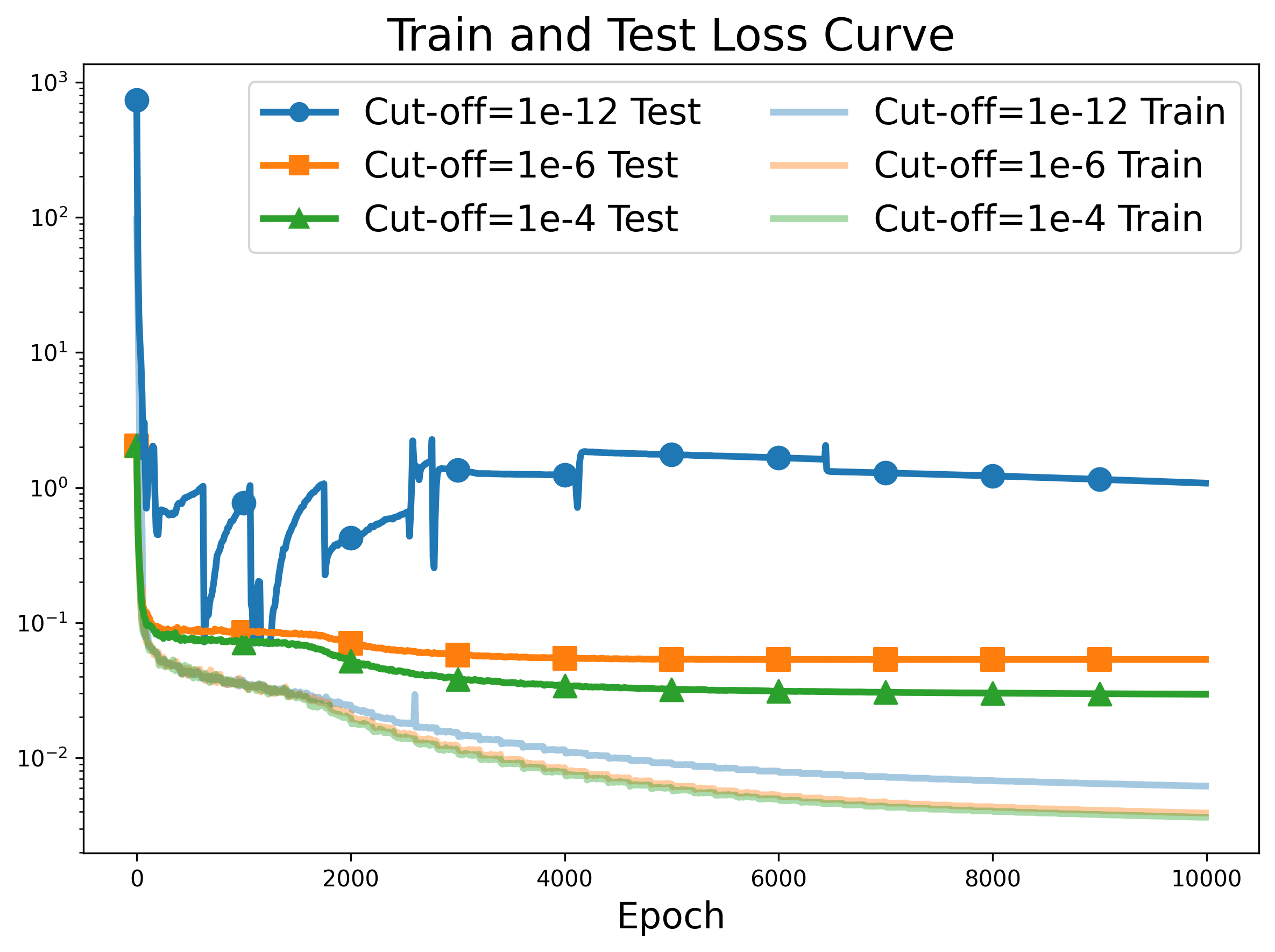}
        \caption*{(a) Training and testing curves.}
    \end{minipage}
    \hspace{0.05\textwidth}
    \begin{minipage}[t]{0.42\textwidth}
        \centering
        \includegraphics[width=\linewidth]{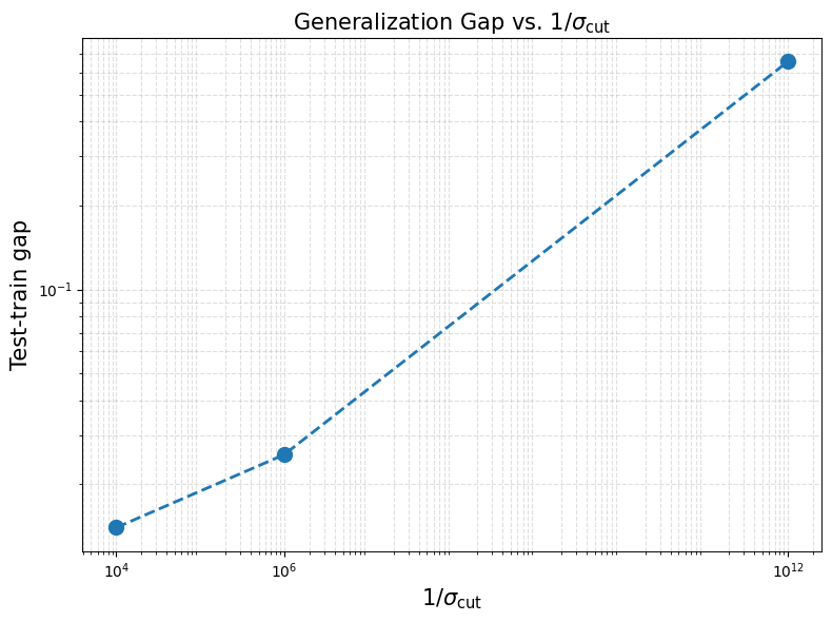}
        \caption*{(b) Gap versus \(1/\sigma_{\mathrm{cut}}\).}
    \end{minipage}
    \caption{
    Numerical illustration of spectral cut-off in coefficient recovery. 
    Smaller thresholds \(\sigma_{\mathrm{cut}}\) retain more singular directions but increase the encoder amplification factor \(\|E_{\mathrm{cut}}\|_2\leq\sigma_{\mathrm{cut}}^{-1}\). 
    Panel (b) plots the empirical train--test RL2E gap against \(1/\sigma_{\mathrm{cut}}\) on a log--log scale.
    }
    \label{fig:gen_curve}
\end{figure}

Figure~\ref{fig:gen_curve} shows this trade-off empirically. 
Smaller cut-off thresholds retain more singular directions and often reduce the training error, but they also increase the encoder amplification factor and lead to a larger empirical train--test gap. 
This behavior is qualitatively consistent with the stability estimate: retaining less stable singular directions can make the coefficient representation more sensitive and degrade test performance.

\section{Numerical Experiments}
\label{sec:numerical-experiments}

We present numerical experiments to evaluate FB-C2CNet from four perspectives: accuracy, training efficiency, sampling flexibility, and applicability to different PDE solution operators. The experiments are organized around the main design questions of the proposed framework. First, we compare coefficient-space learning with pointwise-input variants to quantify the reduction in trainable dimension, number of parameters, and training time. Second, we test non-aligned and resolution-varying observations to examine the flexibility of prescribed-basis encoding and decoding. Third, we study how basis choice, coefficient recovery, and output approximation error affect the final prediction accuracy. Finally, we apply FB-C2CNet to more complex PDE settings, including multi-component operators, nonlinear time-dependent problems, weak-solution learning, high-dimensional problems, and an inverse Stokes boundary-recovery problem from sparse velocity observations.

\subsection{Experimental setup and evaluation metrics}
\label{subsec:exp-setup}

Unless otherwise stated, the coefficient map \(N_\theta\) in FB-C2CNet is implemented as a fully connected neural network with GELU as activation function. 
Its input and output dimensions are determined by the numbers of prescribed input and output basis functions, denoted by \(m_1\) and \(m_2\), respectively. 
For comparisons between C2C and pointwise-input variants, we use the same training data, testing data, and output evaluation points whenever possible. 
Pointwise-input variants take sampled function values as network inputs, whereas FB-C2CNet first encodes input observations into coefficient vectors and then trains the coefficient map.

All neural networks are trained with the Adam optimizer. 
The learning-rate schedules, network architectures, basis parameters, and regularization or cut-off parameters are specified in the corresponding experiments or in the appendix. 
Unless otherwise stated, all timing results are measured on a single NVIDIA RTX 4090 GPU. 
For FB-C2CNet, we report the coefficient-encoding time and the neural-network training time separately. 
The encoding time includes basis evaluation and coefficient recovery for the training samples, while the training time refers only to the optimization of \(N_\theta\). 
When total offline time is reported for FB-C2CNet, it is the sum of encoding time and training time.

For each benchmark, we report the mean relative \(L^2\) error on the testing set,
\begin{equation}
\mathrm{RL2E}
=
\frac{1}{N_{\mathrm{test}}}
\sum_{k=1}^{N_{\mathrm{test}}}
\frac{\|\widehat{\mathbf u}^{(k)}-\mathbf u^{(k)}\|_2}{\|\mathbf u^{(k)}\|_2},
\label{eq:rl2e-exp}
\end{equation}
where \(\mathbf u^{(k)}\) and \(\widehat{\mathbf u}^{(k)}\) denote the reference and predicted output vectors evaluated at the target points. 
When mean squared error is reported, we use
\begin{equation}
\mathrm{MSE}
=
\frac{1}{N_{\mathrm{test}}n_2}
\sum_{k=1}^{N_{\mathrm{test}}}
\|\widehat{\mathbf u}^{(k)}-\mathbf u^{(k)}\|_2^2 .
\label{eq:mse-exp}
\end{equation}
For multi-component outputs, all components are concatenated before computing the error, and \(n_2\) denotes the total number of scalar values in the concatenated output vector.

\subsection{Effect of coefficient-space input representation}
\label{subsec:efficiency}

This subsection evaluates whether replacing pointwise input values by prescribed-basis coefficients reduces the trainable dimension and training cost. 
We compare FB-C2CNet with point-to-coefficient (P2C) variants on one- and two-dimensional Darcy flow benchmarks. 
In all C2C/P2C comparisons below, the output representation, training objective, and optimizer are kept the same; the only difference is whether the neural network input is the coefficient vector or the sampled pointwise values.

\subsubsection{1D Darcy flow: C2C versus P2C}
\label{subsubsec:1d-darcy-c2c-p2c}

We first consider the nonlinear one-dimensional Darcy flow problem
\begin{equation}
\begin{aligned}
\frac{d}{dx}\left(a(u)\frac{du}{dx}\right) &= f(x), \qquad x\in(0,1),\\
u(0)=u(1)&=0,
\end{aligned}
\label{eqn:Exp_1d_darcy}
\end{equation}
where \(a(u(x))=0.2+u^2(x)\). 
The input \(f\) is sampled from a Gaussian random field with covariance kernel
$
k(x,x')=\sigma^2\exp\left(-\frac{|x-x'|^2}{\ell^2}\right),
\quad
\ell=0.04,\quad \sigma=1.0.
$
The dataset is taken from~\cite{ingebrand2025basis,bahmani2025resolution}. 
We use 800 training samples and 200 testing samples on a uniform grid of 2000 points, and learn the operator
$
G:f(x)\mapsto u(x).
$

\begin{figure}[!ht]
    \centering
    \begin{minipage}[t]{0.4\textwidth}
        \centering
        \includegraphics[width=\linewidth]{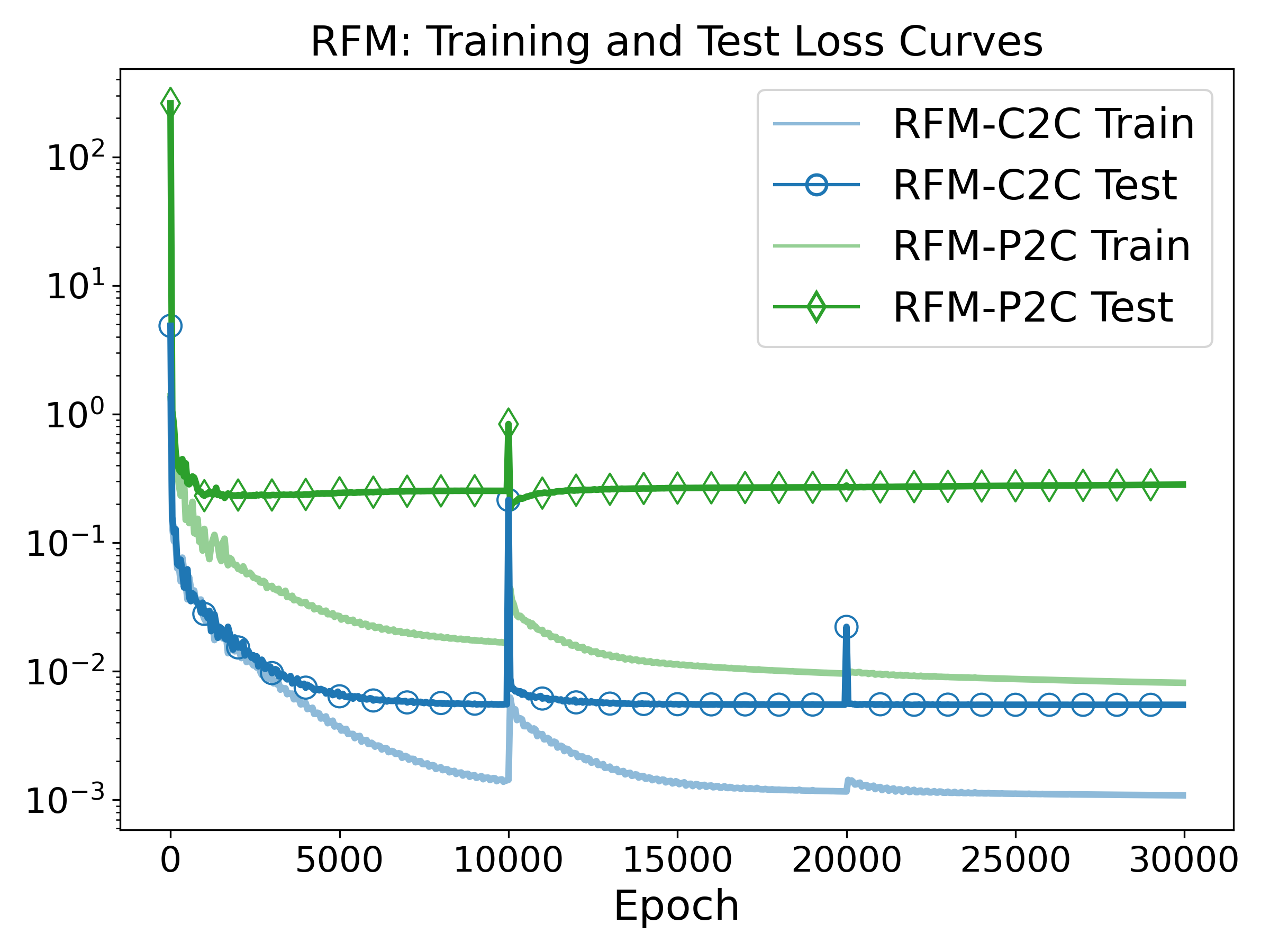}
        \caption*{(a) RFM basis.}
    \end{minipage}
\hspace{0.1\textwidth}
    \begin{minipage}[t]{0.4\textwidth}
        \centering
        \includegraphics[width=\linewidth]{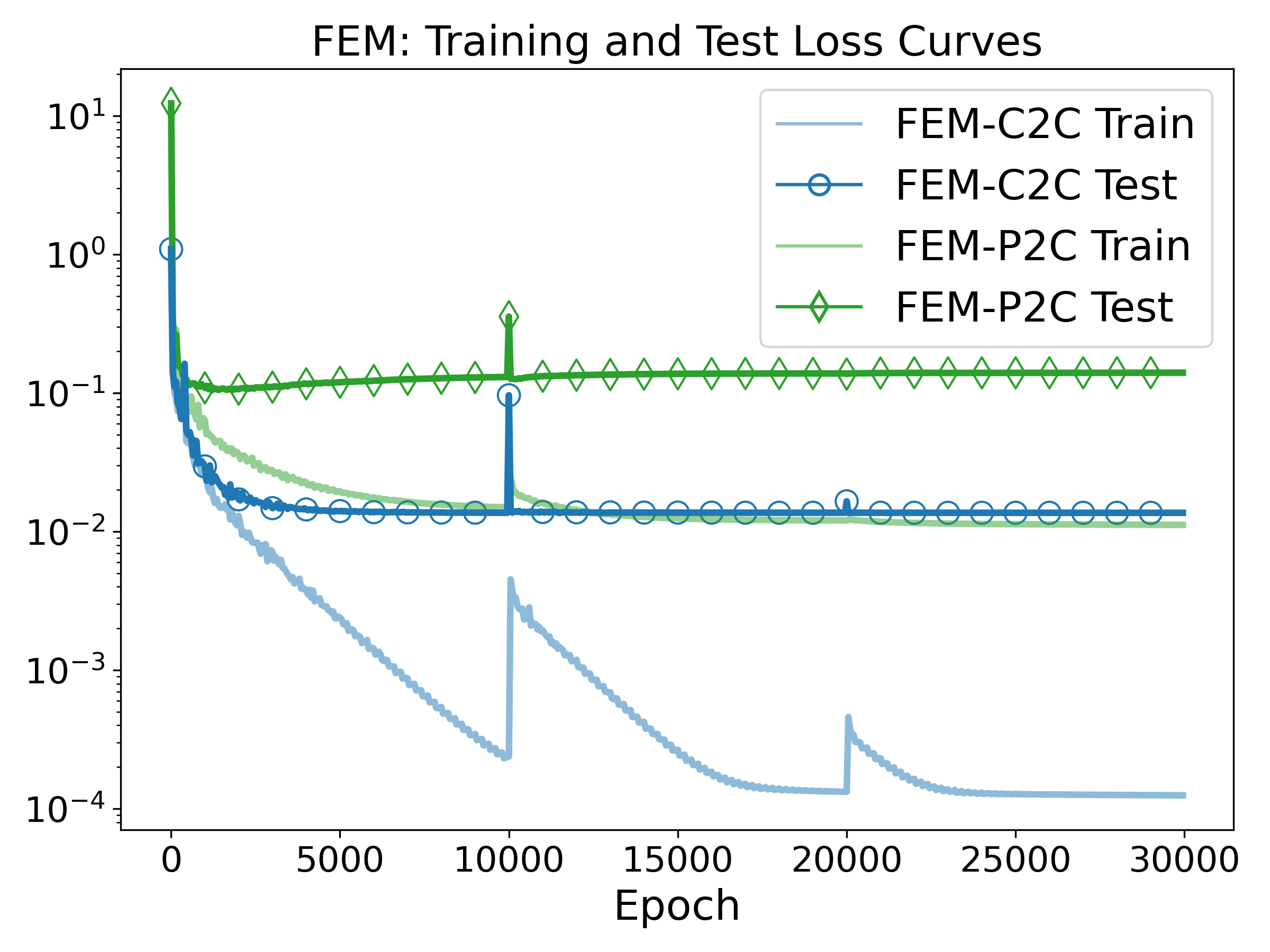}
        \caption*{(b) FEM basis.}
    \end{minipage}
    \caption{
    \textbf{1D Darcy flow.} 
    Training curve comparison between C2C and P2C using RFM and FEM bases. 
    The pointwise-input P2C network takes \(n_1=2000\) sampled values as input, whereas C2C uses \(m_1=128\) input coefficients. 
    The output basis dimension is \(m_2=128\) for both C2C and P2C.
    }
    \label{fig:1d_darcy_curve_2000}
\end{figure}

We compare RFM-C2C and FEM-C2C with their P2C counterparts. 
For C2C, the neural network input dimension is \(m_1=128\), determined by the number of prescribed input basis functions. 
For P2C, the input dimension is \(n_1=2000\), given by the number of sampled grid values. 
As shown in Figure~\ref{fig:1d_darcy_curve_2000} and Table~\ref{tab:1d_darcy_time}, the C2C formulation achieves substantially lower test errors while using a much lower-dimensional input representation. 
One possible reason is that the coefficient input leads to a simpler and better-conditioned learning problem: the C2C network has a much smaller input layer than the P2C network, which makes the optimization easier and allows the training loss to decrease more effectively.

\begin{table}[!htpb]
    \centering
    \caption{
    1D Darcy flow.
    Comparison of C2C and P2C using RFM and FEM bases. 
    The encoder cost includes basis generation and projection of input functions into the prescribed basis space. 
    The training time refers to the optimization time of the neural operator.
    Within each basis, the lower C2C/P2C training time and test error are shown in bold.
    }
    \label{tab:1d_darcy_time}
    \footnotesize
    \renewcommand{\arraystretch}{1.12}
    \setlength{\tabcolsep}{3.2pt}
    \begin{tabular}{@{}lccccccc@{}}
        \toprule
        & \multicolumn{1}{c}{Encoding time}
        & \multicolumn{2}{c}{Training time}
        & \multicolumn{2}{c}{Test RL2E}
        & \multicolumn{2}{c}{Network architecture}\\
        \cmidrule(lr){2-2}\cmidrule(lr){3-4}\cmidrule(lr){5-6}\cmidrule(l){7-8}
        Basis & Encoder & C2C & P2C & C2C & P2C & C2C & P2C \\ 
        \midrule
        RFM 
        & 0.93 s
        & \textbf{49.32 s}
        & 53.54 s
        & \textbf{3.870e-3} 
        & 2.829e-1 
        & [128,400\(\times\)3,128] 
        & [2000,400\(\times\)3,128] \\ 
        FEM 
        & 1.75 s
        & \textbf{52.26 s}
        & 53.17 s
        & \textbf{1.367e-2}
        & 1.405e-1 
        & [128,400\(\times\)3,128]
        & [2000,400\(\times\)3,128] \\ 
        \bottomrule
    \end{tabular}
\end{table}

\subsubsection{2D Darcy flow: C2C versus P2C}
\label{subsubsec:2d-darcy-c2c-p2c}

We next consider the two-dimensional Darcy flow equation on \(\Omega=(0,1)^2\):
\begin{equation}
\begin{aligned}
-\nabla\cdot\left(\kappa(x,y)\nabla u(x,y)\right) &= f(x,y), \qquad (x,y)\in\Omega,\\
u(x,y)&=0, \qquad (x,y)\in\partial\Omega.
\end{aligned}
\label{eqn:2d_darcy}
\end{equation}
Here \(f(x,y)\equiv 1\), and the diffusion coefficient is sampled as \(\kappa=\exp(\psi)\), where
$
\psi\sim\mathcal N\left(0,(-\Delta+9I)^{-2}\right).
$
The dataset is taken from~\cite{li2021fourier,hua2023basis}. 
It contains 2400 training samples and 600 testing samples on a uniform \(141\times141\) grid. 
The target operator is
$
G:\kappa(x,y)\mapsto u(x,y).
$

\begin{figure}[!ht]
    \centering
    \begin{minipage}[t]{0.4\textwidth}
        \centering
        \includegraphics[width=\linewidth]{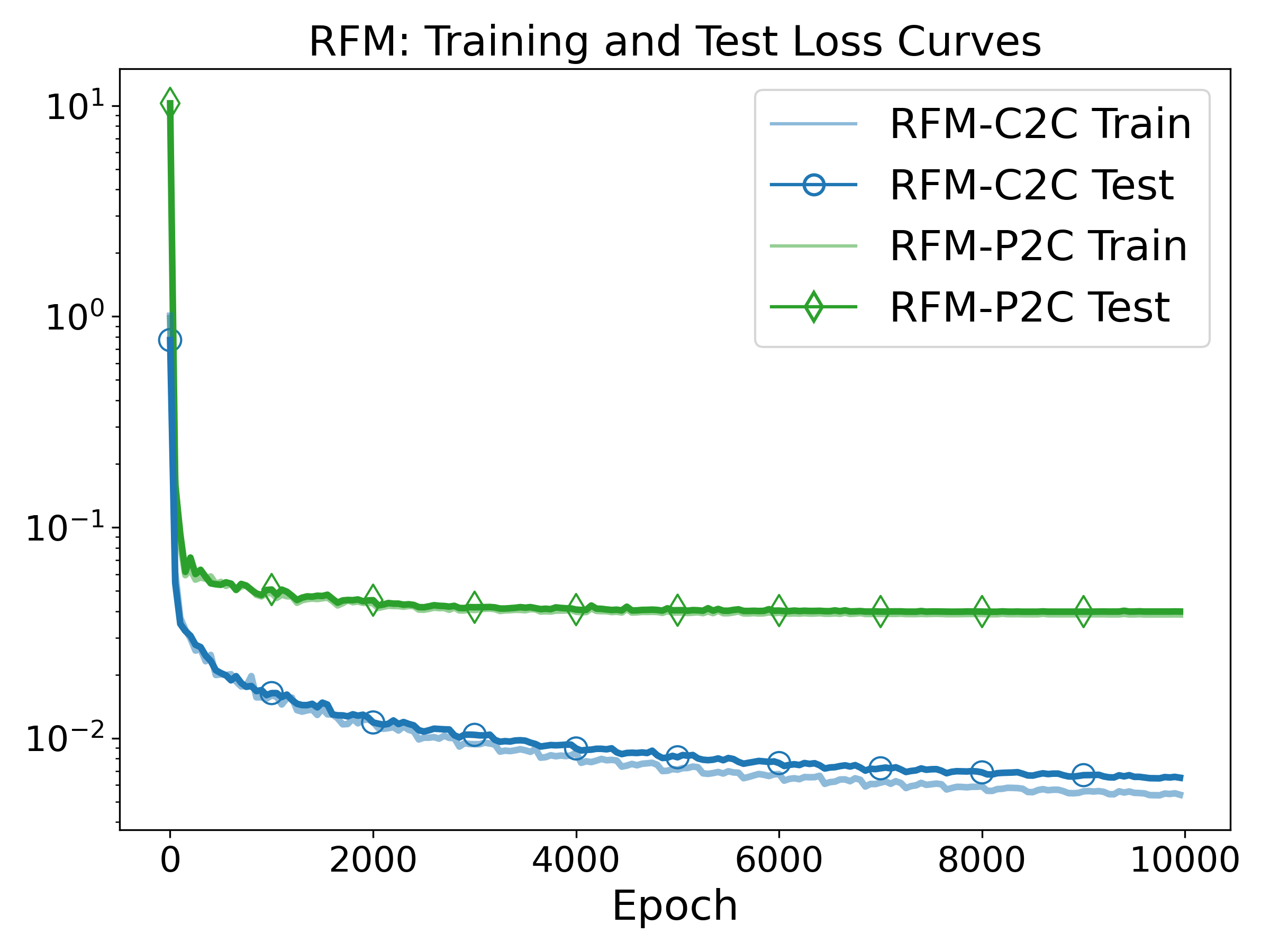}
        \caption*{(a) RFM basis.}
    \end{minipage}
\hspace{0.1\textwidth}
    \begin{minipage}[t]{0.4\textwidth}
        \centering
        \includegraphics[width=\linewidth]{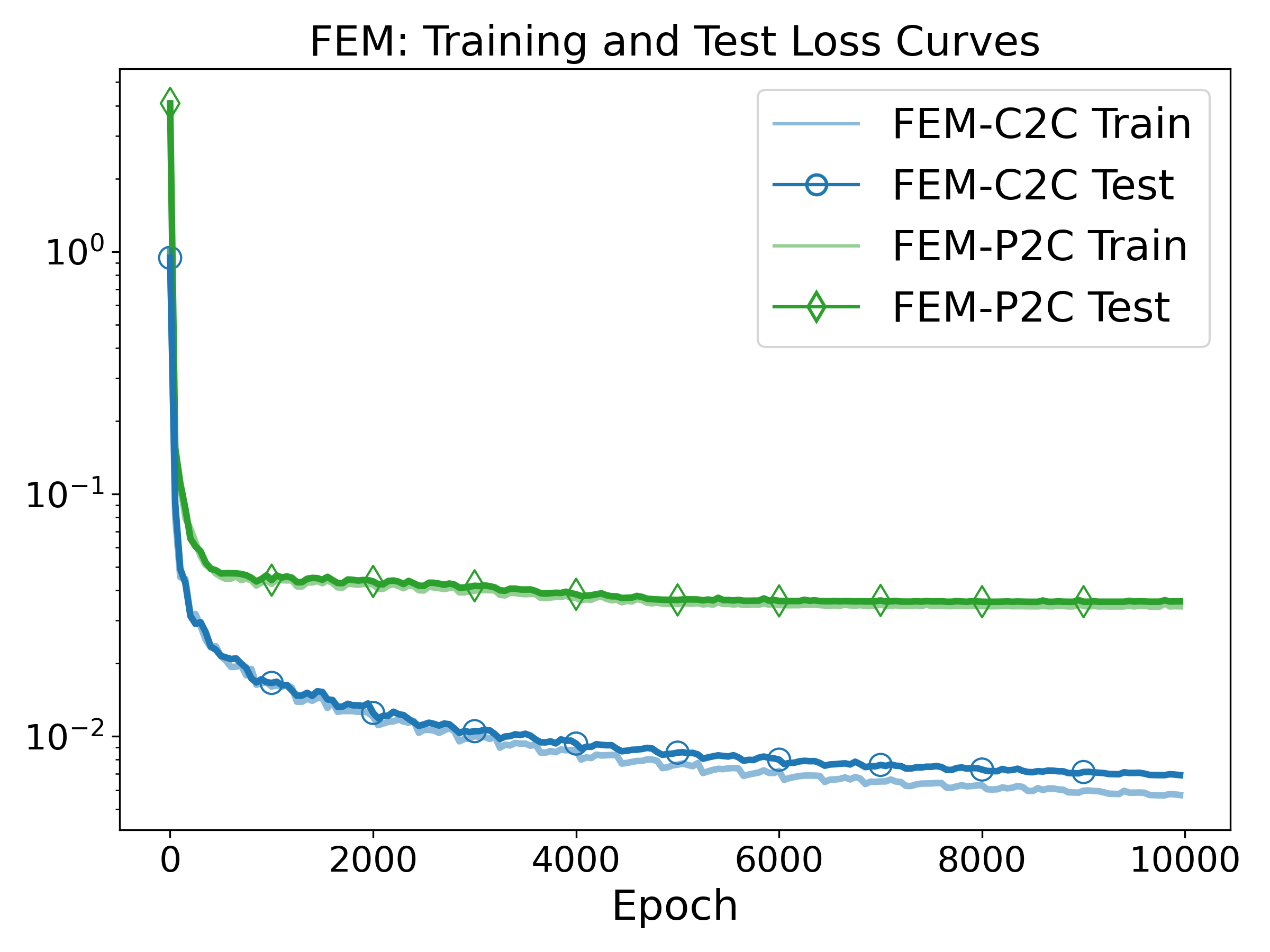}
        \caption*{(b) FEM basis.}
    \end{minipage}
    \caption{
    \textbf{2D Darcy flow.} 
    Training curve comparison between C2C and P2C using RFM and FEM bases. 
    The number of grid points is \(n_1=n_2=141\times141\). 
    The prescribed basis dimensions are \((m_1,m_2)=(512,2048)\) for RFM and \((m_1,m_2)=(553,2059)\) for FEM.
    }
    \label{fig:2d_darcy_exp_curve}
\end{figure}

We first compare C2C and P2C using the same RFM and FEM output bases. 
For RFM, we use \((m_1,m_2)=(512,2048)\); for FEM, we use \((m_1,m_2)=(553,2059)\). 
The P2C variants instead take all \(141^2\) pointwise input values as the network input. 
Figure~\ref{fig:2d_darcy_exp_curve} and Table~\ref{tab:2d_darcy_time} show that the coefficient-space representation improves both accuracy and training efficiency. 
The FEM encoder is more expensive than the RFM encoder in this two-dimensional setting because it requires mesh construction and finite element basis evaluation.

\begin{table}[!htpb]
    \centering
    \caption{
    2D Darcy flow.
    Comparison of C2C and P2C using RFM and FEM bases. 
    The encoder cost includes basis generation and projection of input functions into the prescribed basis space. 
    The training time refers to the optimization time of the neural operator.
    Within each basis, the lower C2C/P2C training time and test error are shown in bold.
    }
    \label{tab:2d_darcy_time}
    \footnotesize
    \renewcommand{\arraystretch}{1.12}
    \setlength{\tabcolsep}{3pt}
    \begin{tabular}{@{}lccccccc@{}}
        \toprule
        & \multicolumn{1}{c}{Encoding time}
        & \multicolumn{2}{c}{Training time}
        & \multicolumn{2}{c}{Test RL2E}
        & \multicolumn{2}{c}{Network architecture}\\
        \cmidrule(lr){2-2}\cmidrule(lr){3-4}\cmidrule(lr){5-6}\cmidrule(l){7-8}
        Basis & Encoder & C2C & P2C & C2C & P2C & C2C & P2C \\ 
        \midrule
        RFM 
        & 3.67 s
        & \textbf{274.40 s}
        & 302.11 s
        & \textbf{6.740e-3} 
        & 3.990e-2 
        & [512,180\(\times\)2,2048]
        & [\(141^2\),180\(\times\)2,2048] \\ 
        FEM 
        & 447.05 s
        & \textbf{291.17 s}
        & 317.65 s
        & \textbf{6.920e-3}
        & 3.600e-2 
        & [553,180\(\times\)2,2059]
        & [\(141^2\),180\(\times\)2,2059] \\ 
        \bottomrule
    \end{tabular}
\end{table}

% \subsubsection{Comparison with other neural operators}
% We further compare FB-C2CNet with representative neural operator methods on the same dataset. 
% Table~\ref{tab:2d_darcy_comparison} reports the training time, number of trainable parameters, RL2E, and MSE. 
% RFM-C2C achieves accuracy comparable to BasisONet and FNO2D while requiring substantially less training time. 
% FNO2D attains the smallest error in this test, but uses more parameters and a longer training time. 
% DeepONet trains quickly but has significantly larger error.

% \begin{table}[!htbp]
% \centering
% \caption{
% \textbf{2D Darcy flow.} 
% Performance comparison on the rectangular-domain Darcy benchmark. 
% For C2C methods, the reported time includes the coefficient-encoding time and the neural-network training time.
% }
% \label{tab:2d_darcy_comparison}
% \begin{tabular}{|c|c|c|c|c|}
% \hline 
% Method & Time & Params & RL2E & MSE \\
% \hline 
% BasisONet~\cite{hua2023basis} 
% & 88m1s 
% & 2,712,174 
% & 7.02e-3 
% & 9.12e-8 \\
% \hline 
% B2B~\cite{ingebrand2025basis} 
% & 94m2s 
% & 499,758 
% & 1.174e-2 
% & 2.731e-7 \\
% \hline 
% DeepONet~\cite{ingebrand2025basis} 
% & 4m10s 
% & 504,926 
% & 1.151e-1 
% & 2.588e-5 \\
% \hline
% FNO2D~\cite{lu2022comprehensive} 
% & 31m22s 
% & 1,188,353 
% & \textbf{5.43e-3} 
% & \textbf{5.19e-8} \\
% \hline 
% RFM-C2C 
% & \textbf{5s+3m17s} 
% & \textbf{492,864} 
% & 6.47e-3 
% & 8.68e-8 \\
% \hline 
% FEM-C2C 
% & 7m27s+4m51s 
% & 504,979 
% & 6.92e-3 
% & 9.03e-8 \\
% \hline 
% \end{tabular}
% \end{table}

\subsection{Comparison with representative neural operators}
\label{subsec:nonaligned-resolution}

This subsection assesses the flexibility of FB-C2CNet when input and output functions are not sampled on a common fixed grid and when the inference resolution differs from the training resolution. These evaluations use the 1D Darcy flow benchmark introduced in Section~\ref{subsubsec:1d-darcy-c2c-p2c}. We also use the 2D Darcy flow problem from Section~\ref{subsubsec:2d-darcy-c2c-p2c} to evaluate the accuracy and efficiency of the proposed method.

% This subsection evaluates the flexibility of FB-C2CNet when the input and output functions are not observed on a common fixed grid, and when the evaluation resolution differs from the training resolution. 
% Both tests are performed on the 1D Darcy flow benchmark introduced in Section~\ref{subsubsec:1d-darcy-c2c-p2c}.

\subsubsection{Accuracy and efficiency on the 2D Darcy benchmark}
\label{subsubsec:standard-2d-darcy-comparison}

We compare FB-C2CNet with representative neural operator architectures on the aligned-grid 2D Darcy benchmark introduced in Section~\ref{subsubsec:2d-darcy-c2c-p2c}. 
This experiment compares point-sensor, grid-based, data-driven reduced-basis, learned-basis, and prescribed-basis C2C methods on the same dataset and output grid. 
For C2C methods, the reported time includes both coefficient encoding and neural-network training.
For PCA-Net, we retain 16 input and 16 output PCA coefficients, i.e., \(m_1=m_2=16\), and use the coefficient network \([16,512\times3,16]\), which has 542,224 trainable parameters.

Table~\ref{tab:2d_darcy_comparison} reports training time, parameter count, RL2E, and MSE. 
FB-C2CNet achieves competitive accuracy with much lower total offline cost. 
In particular, RFM-C2C reaches an RL2E of \(6.47\mathrm{e}{-3}\) using only \(5\mathrm{s}+3\mathrm{m}17\mathrm{s}\) for encoding and training, highlighting the efficiency gained by prescribing the basis and training only the coefficient map.

\begin{table}[!htbp]
\centering
\caption{
2D Darcy flow.
Performance comparison on the rectangular-domain Darcy benchmark. 
For C2C methods, the reported time includes the coefficient-encoding time and the neural-network training time.
}
\label{tab:2d_darcy_comparison}
\small
\renewcommand{\arraystretch}{1.12}
\begin{tabular}{@{}lcccc@{}}
\toprule
Method & Time & Parameters & Test RL2E & MSE \\
\midrule
PCA-Net~\cite{bhattacharya2021model} & 54 s & 542,224 & 1.07e-2& 2.22e-7 \\
BasisONet~\cite{hua2023basis} 
& 88 min 1 s
& 2,712,174 
& 7.02e-3 
& 9.12e-8 \\
B2B~\cite{ingebrand2025basis} 
& 94 min 2 s
& 499,758 
& 1.17e-2 
& 2.73e-7 \\
DeepONet~\cite{ingebrand2025basis} 
& 4 min 10 s
& 504,926 
& 1.15e-1 
& 2.59e-5 \\
FNO2D~\cite{lu2022comprehensive} 
& 31 min 22 s
& 1,188,353 
& \textbf{5.43e-3} 
& \textbf{5.19e-8} \\
RFM-C2C 
& \textbf{5 s + 3 min 17 s}
& \textbf{492,864} 
& 6.47e-3 
& 8.68e-8 \\
FEM-C2C 
& 7 min 27 s + 4 min 51 s
& 504,979 
& 6.92e-3 
& 9.03e-8 \\
\bottomrule
\end{tabular}
\end{table}

\newpage
\subsubsection{Non-aligned input-output observations}
\label{subsubsec:nonaligned-random}

\begin{figure}[!ht]
    \centering
    \begin{minipage}[t]{0.4\textwidth}
        \centering
        \includegraphics[width=\linewidth]{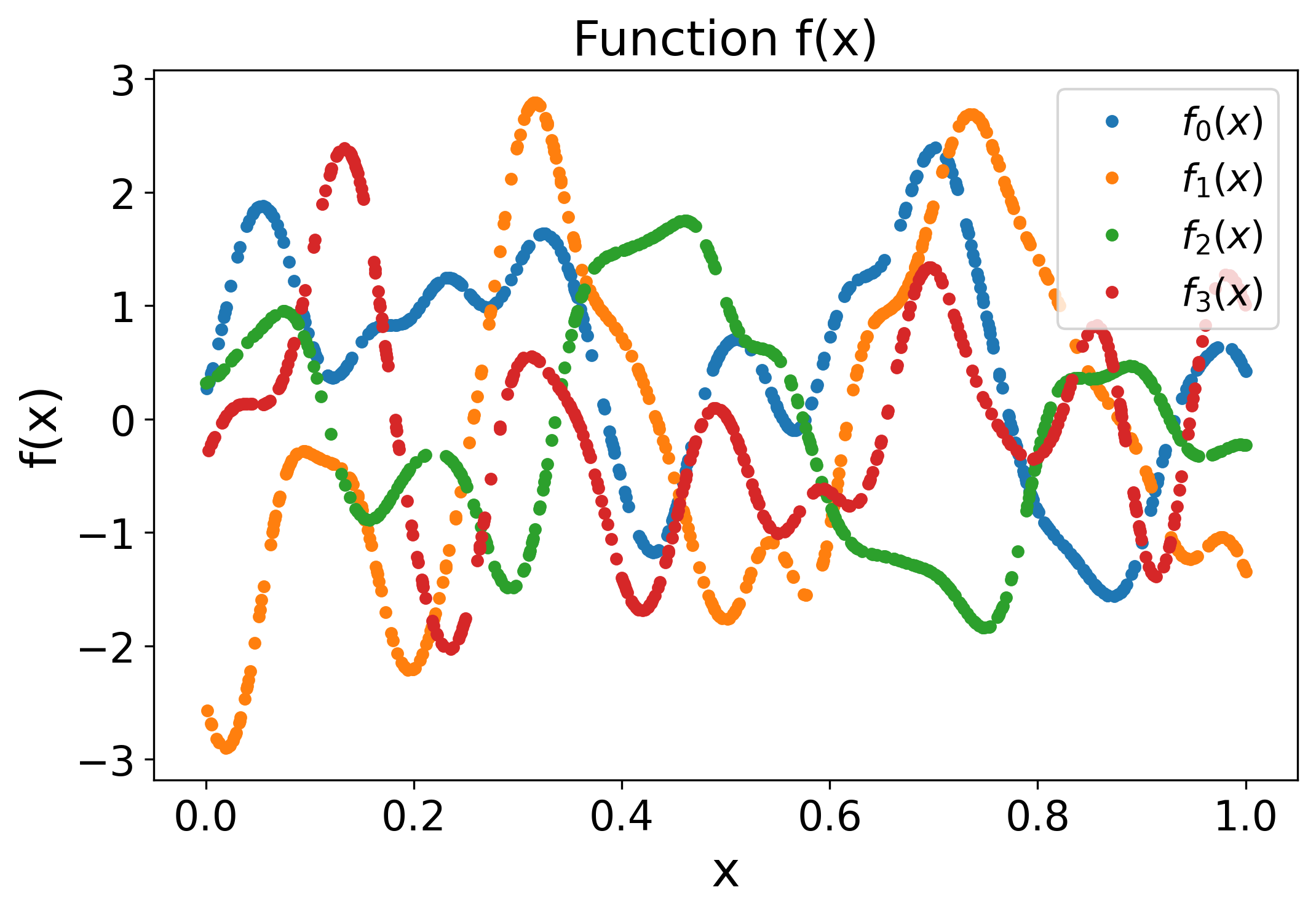}
        \caption*{(a) Input function \(f(x)\).}
    \end{minipage}
    \hspace{0.1\textwidth}
    \begin{minipage}[t]{0.4\textwidth}
        \centering
        \includegraphics[width=\linewidth]{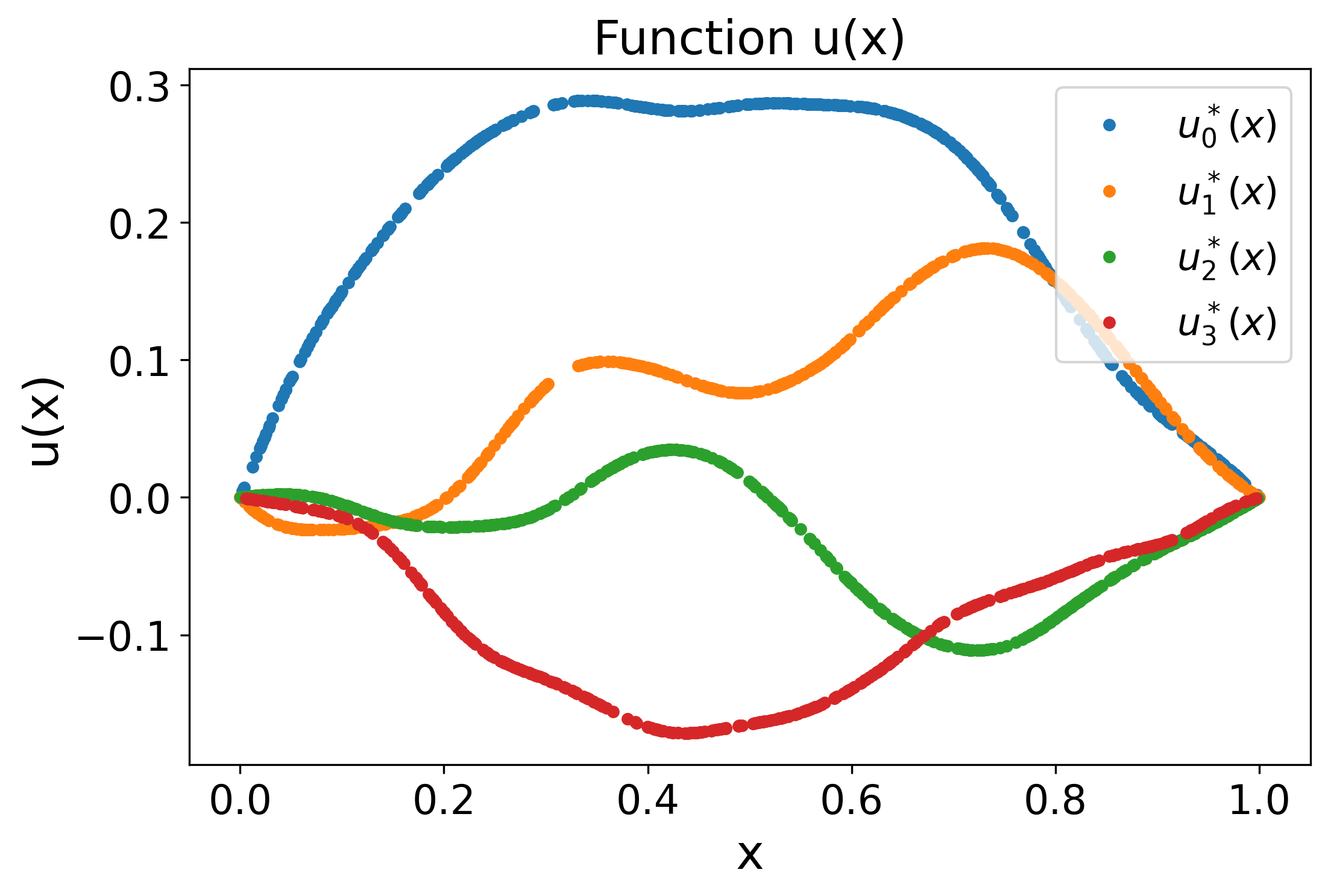}
        \caption*{(b) Output function \(u(x)\).}
    \end{minipage}
    \caption{
    \textbf{1D Darcy flow with non-aligned observations.} 
    For each sample, 400 random locations are independently selected for the input function \(f(x)\), and another 400 random locations are independently selected for the output function \(u(x)\). 
    The two observation sets are non-aligned and have no pointwise correspondence.
    }
    \label{fig:1d_darcy_random_2000}
\end{figure}

We first consider a non-uniform and non-aligned observation setting. 
For each input-output pair, we independently draw 400 random points from the original 2000-point uniform grid for the input function \(f\), and another 400 random points for the output function \(u\). 
Thus, the input and output observations have different locations and no pointwise correspondence as shown in Figure~\ref{fig:1d_darcy_random_2000}.
This setting directly tests whether the prescribed-basis encoding and decoding procedures can handle observations beyond a common fixed grid.

\begin{figure}[!ht]
    \centering
    \begin{minipage}[t]{0.4\textwidth}
        \centering
        \includegraphics[width=\linewidth]{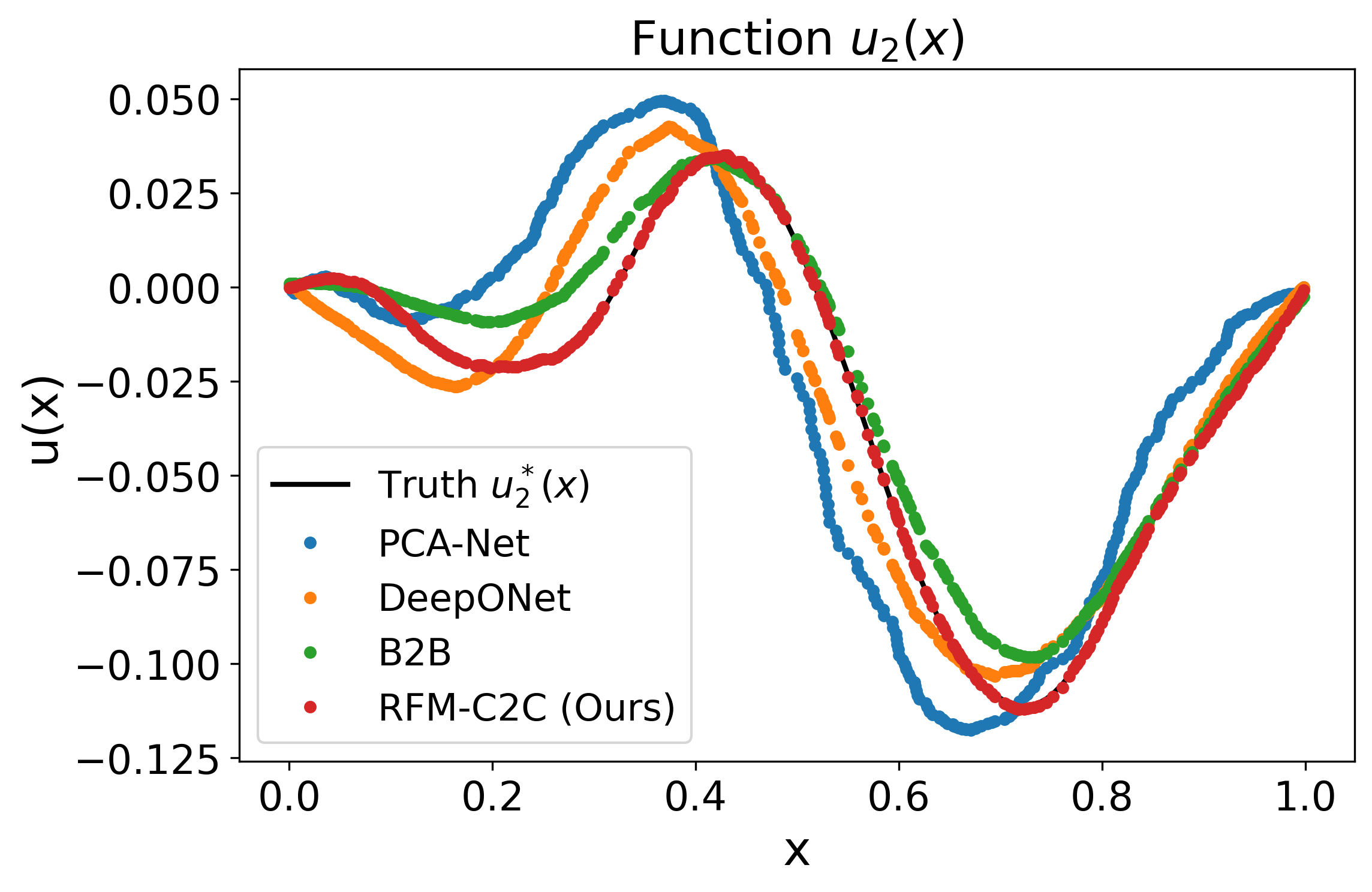}
        \caption*{(a) Prediction of \(u(x)\).}
    \end{minipage}
    \hspace{0.1\textwidth}
    \begin{minipage}[t]{0.4\textwidth}
        \centering
        \includegraphics[width=\linewidth]{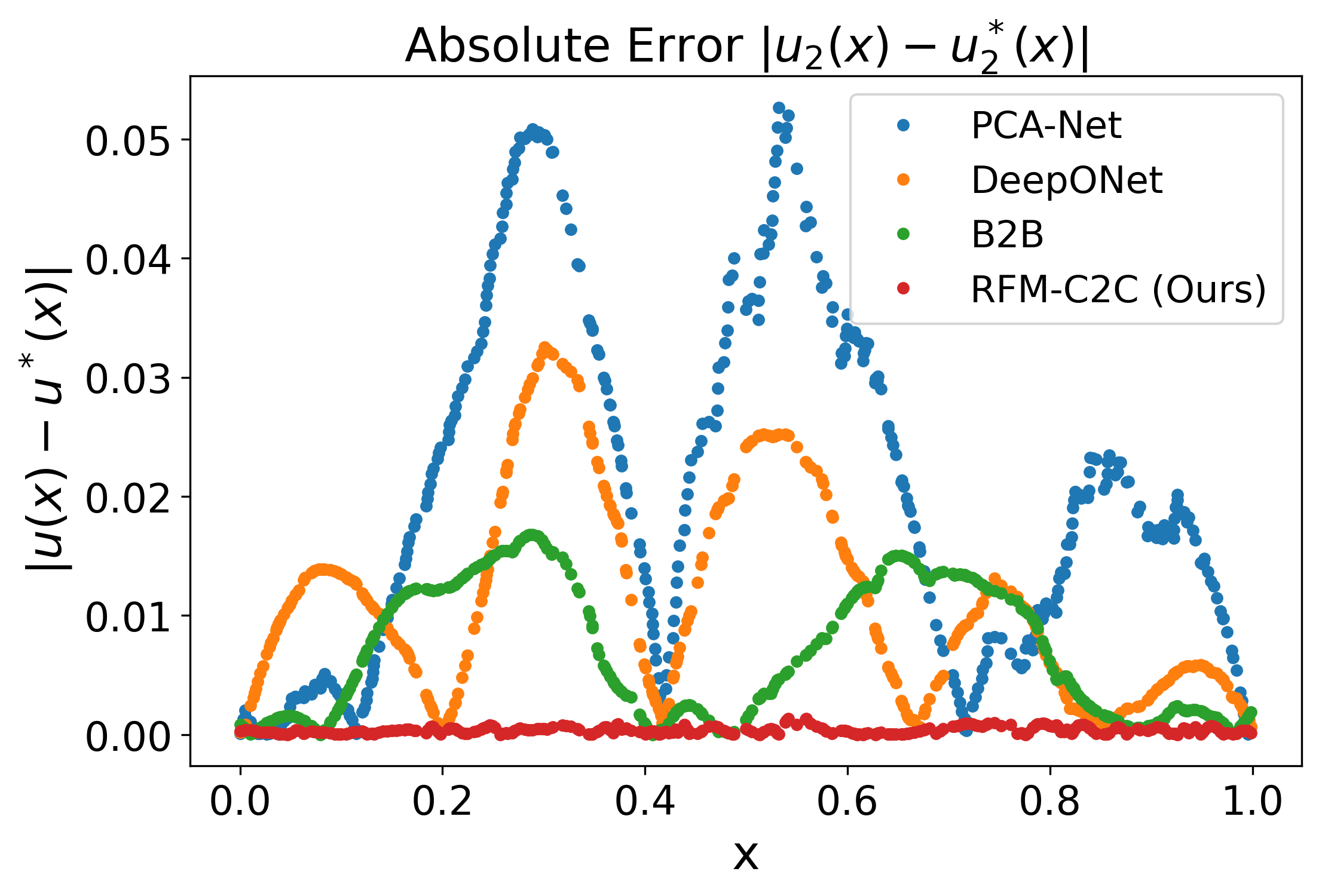}
        \caption*{(b) Pointwise absolute error.}
    \end{minipage}
    \caption{
    \textbf{1D Darcy flow with non-aligned observations.} 
    Pointwise predictions and absolute errors for a representative test sample using PCA-Net, DeepONet, B2B, and RFM-C2C. 
    The RFM-C2C prediction is obtained from prescribed-basis coefficient recovery rather than from aligned pointwise input-output pairs.
    }
    \label{fig:1d_darcy_random_2000_u2}
\end{figure}

We apply RFM-C2C with \(m_1=m_2=128\) and compare it with RFM-P2C, RFM-C2P, PCA-Net, DeepONet, and B2B. 
For the compared reduced models, the tuned dimensions are \(m_1=m_2=16\) for PCA-Net, \(m_1=m_2=60\) for DeepONet, and \(m_1=m_2=80\) for B2B. 
Standard FNO is not included in this comparison because it assumes a fixed structured grid. 
As shown in Table~\ref{tab:1d_darcy_random_comparison}, RFM-C2C achieves the lowest test error under non-aligned observations. 
A representative prediction is shown in Figure~\ref{fig:1d_darcy_random_2000_u2}.

\begin{table}[!htbp]
\centering
\caption{
1D Darcy flow with non-aligned observations.
Performance comparison on the dataset with independently sampled input and output observation locations. 
``Parameters'' denotes the number of trainable parameters.
}
\label{tab:1d_darcy_random_comparison}
\footnotesize
\renewcommand{\arraystretch}{1.12}
\begin{tabularx}{\textwidth}{@{}lcc>{\raggedright\arraybackslash}X@{}}
\toprule
Method & Parameters & Test RL2E & Network architecture\\
\midrule
RFM-P2C  
& 1,172,528 
& 2.47e-1 
& [2000,400\(\times\)3,128] \\
RFM-C2P 
& 1,174,400 
& 8.87e-2 
& [128,400\(\times\)3,2000]  \\
PCA-Net~\cite{bhattacharya2021model} 
& 542,224 
& 2.88e-1 
& [16,512\(\times\)3,16]\\
DeepONet~\cite{lu2021deeponet} 
& 500,803 
& 3.02e-1 
& [400,294\(\times\)3,60] + [1,294\(\times\)3,60] + 1  \\
B2B~\cite{ingebrand2025basis} 
& 499,352 
& 9.15e-2 
& [1,262\(\times\)3,80] + [80,262\(\times\)3,81] + [1,262\(\times\)3,81] \\
\textbf{RFM-C2C (Ours)}  
& \textbf{423,728} 
& \textbf{1.18e-2} 
& [128,400\(\times\)3,128] \\
\bottomrule
\end{tabularx}
\end{table}

\subsubsection{Resolution transfer}
\label{subsubsec:resolution-transfer}

\begin{figure}[!ht]
    \centering
    \begin{minipage}[t]{0.45\textwidth}
        \centering
        \includegraphics[width=\linewidth]{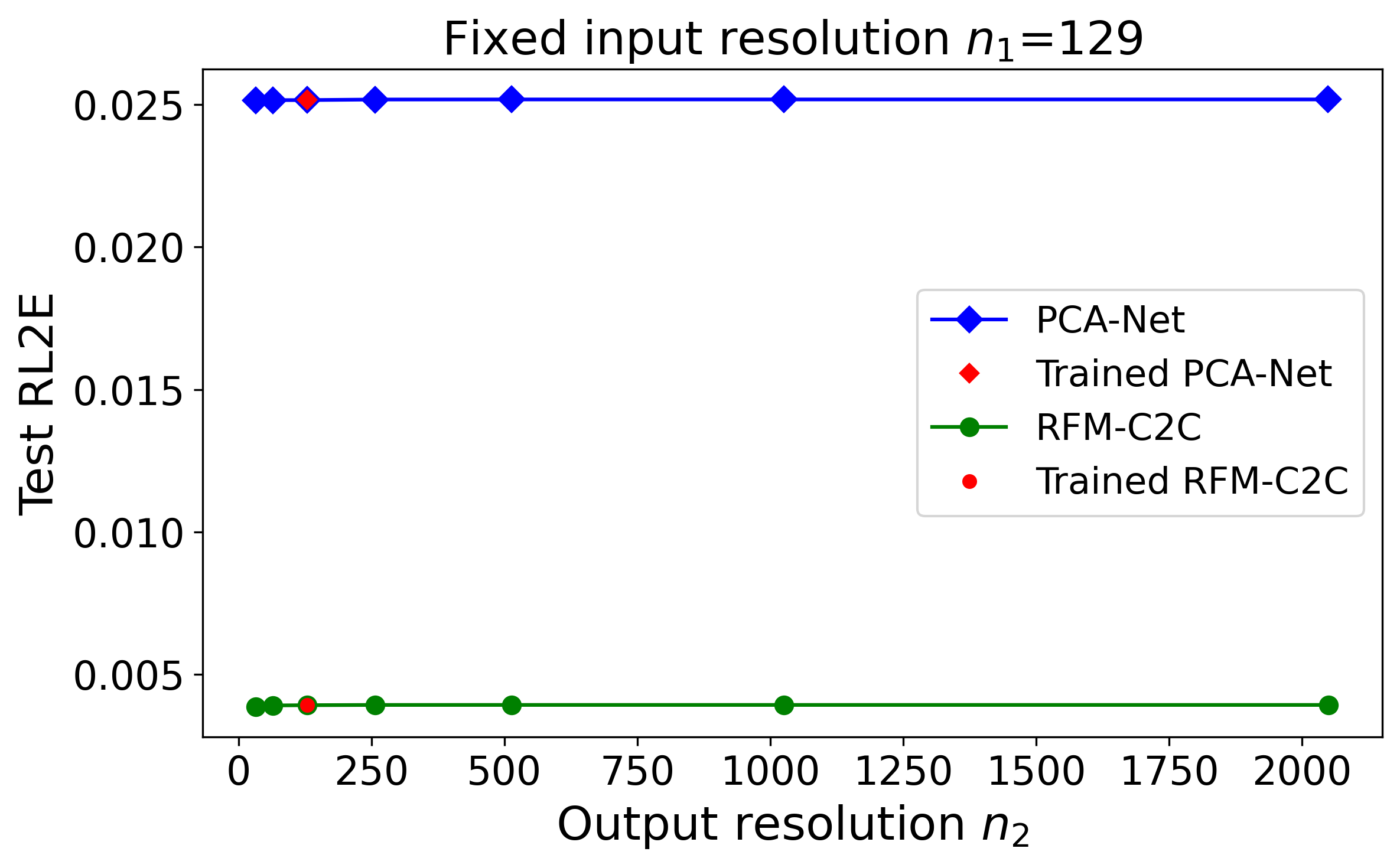}
        \caption*{(a) Input resolution fixed.}
    \end{minipage}
    \hfill
    \begin{minipage}[t]{0.45\textwidth}
        \centering
        \includegraphics[width=\linewidth]{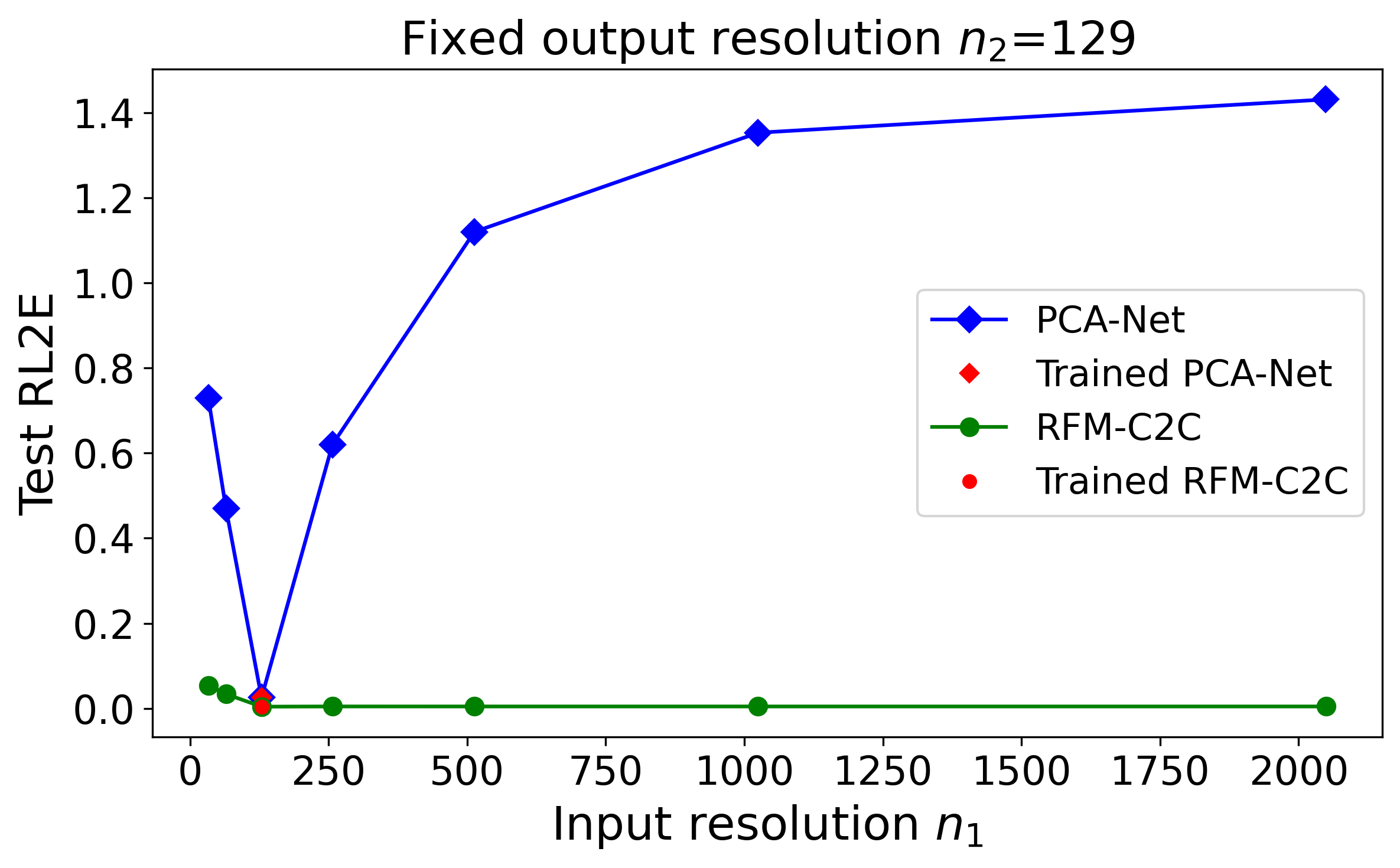}
        \caption*{(b) Output resolution fixed.}
    \end{minipage}
    \caption{
    \textbf{Resolution transfer on 1D Darcy flow.} 
    PCA-Net and RFM-C2C are trained at \(n_1=n_2=129\). 
    In panel (a), the input resolution is fixed at \(n_1=129\), while the output resolution is varied. 
    In panel (b), the output resolution is fixed at \(n_2=129\), while the input resolution is varied. 
    The tested resolutions are $\{33,65,129,257,513,1025,2049\}$.
    }
    \label{fig:1D_Darcy_resolution_PCA}
\end{figure}

We next test whether a trained model can be evaluated at resolutions different from the training resolution. 
First, PCA-Net and RFM-C2C are trained on \(n_1=n_2=129\) uniformly sampled points on \([0,1]\). 
We then vary either the output resolution or the input resolution over
$
\{33,\ 65,\ 129,\ 257,\ 513,\ 1025,\ 2049
\}$. For PCA-Net, we use the resolution-transfer strategy in~\cite{bhattacharya2021model}; for RFM-C2C, the same continuous RFM basis functions are directly evaluated at the new input or output locations. 
Figure~\ref{fig:1D_Darcy_resolution_PCA} shows that RFM-C2C is more robust under both input and output resolution changes.

\begin{figure}[!ht]
    \centering
    \includegraphics[width=0.4\textwidth]{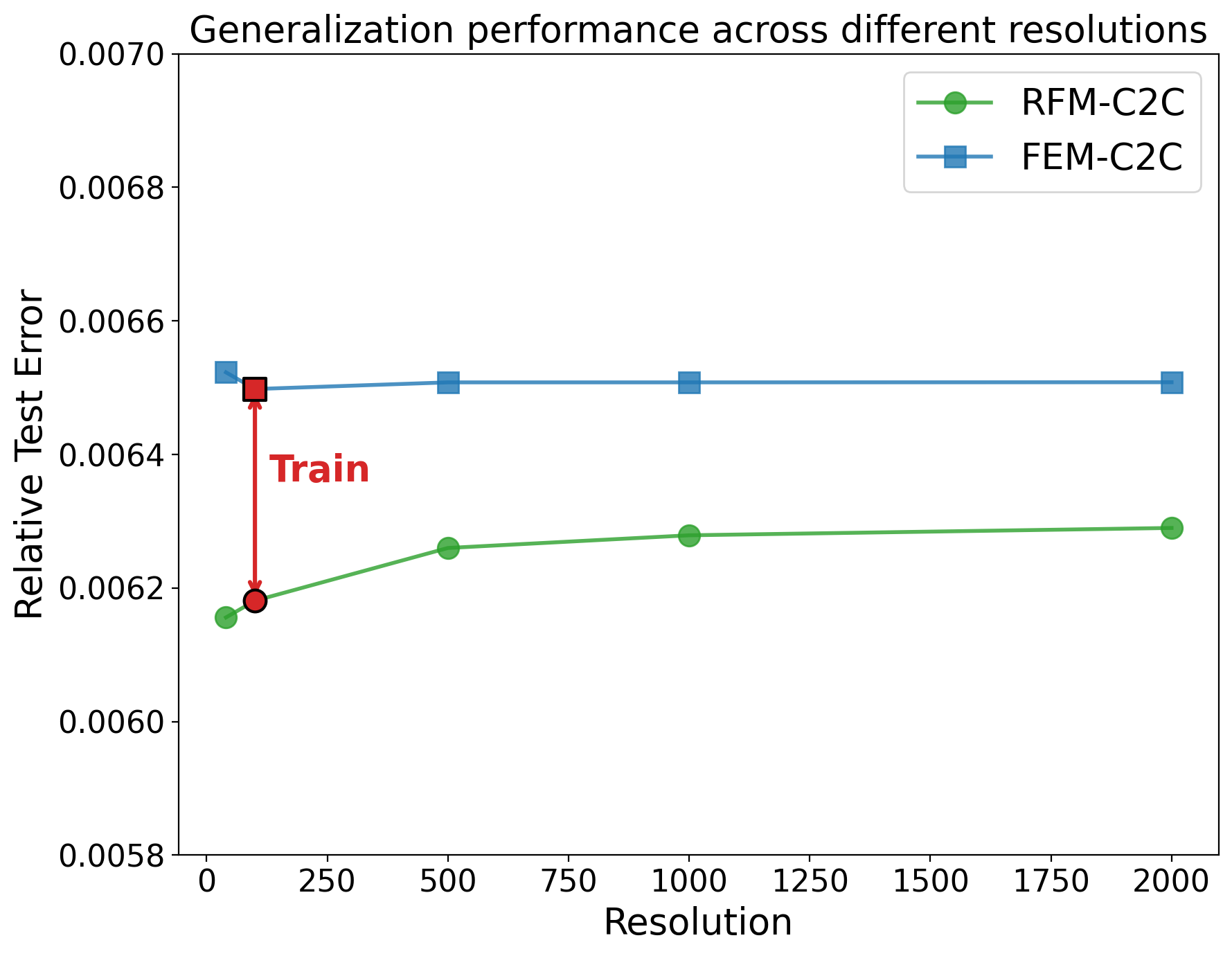}
    \caption{
    \textbf{Resolution transfer with RFM and FEM bases.} 
    RFM-C2C and FEM-C2C are trained on a uniform \(n_1=n_2=100\) grid and tested at output resolutions \(n_2=40,500,1000,2000\). 
    Both prescribed-basis methods remain stable under unseen output resolutions.
    }
    \label{fig:1D_Darcy_gen_res}
\end{figure}
We also compare RFM-C2C and FEM-C2C in a separate output-resolution transfer test. 
Both models are trained on a 100-point uniform grid obtained by taking every 20th point from the original 2000-point grid. 
They are then evaluated at output resolutions \(n_2=40,500,1000,2000\). 
As shown in Figure~\ref{fig:1D_Darcy_gen_res}, both prescribed-basis methods remain stable across unseen output resolutions, with RFM-C2C achieving lower error in this experiment.

\begin{comment}
\subsubsection{Coefficient recovery and regularization}
\label{subsubsec:coefficient-recovery-exp}

The coefficient-recovery step is important whenever the basis evaluation matrix is ill-conditioned or when observations are scattered. 
In the experiments, RFM and RBF coefficients are computed by regularized least squares or truncated-SVD fitting, while FEM coefficients are obtained by the same prescribed-basis fitting principle on the chosen mesh. 
The regularization or cut-off parameter controls the trade-off between stable coefficient recovery and approximation bias, in agreement with the analysis in Section~3. 
The specific basis parameters and cut-off values used in each benchmark are reported in the appendix.
\end{comment}

\newpage
\subsection{Complex geometries and general PDE operators}
\label{subsec:complex-operators}

This subsection examines whether FB-C2CNet can be extended beyond scalar elliptic operators on regular domains.
We consider irregular-domain problems, multi-component operators, nonlinear time-dependent dynamics, weak-solution learning, and a high-dimensional proof-of-concept.
The goal is not to exhaustively optimize each benchmark, but to demonstrate that the prescribed-basis coefficient formulation can be adapted to different PDE settings by choosing suitable input and output approximation spaces.

\subsubsection{Irregular-domain Poisson Equation}

In this example, we illustrate this point on a complex-domain Poisson problem and then examine the output projection error induced by the prescribed decoder.

The Darcy experiments above show that RFM bases can provide accurate and efficient coefficient representations on regular domains and under non-aligned sampling. 
In contrast, FEM bases can be advantageous when the computational geometry is nontrivial or when the mesh already encodes boundary information. 
This complementary behavior is consistent with the design philosophy of FB-C2CNet: the basis is not learned from data, but chosen according to the geometry, sampling pattern, and regularity of the PDE solution space.

\begin{figure}[!ht]
    \centering
    \begin{minipage}[t]{0.3\textwidth}
        \centering
        \includegraphics[width=\linewidth]{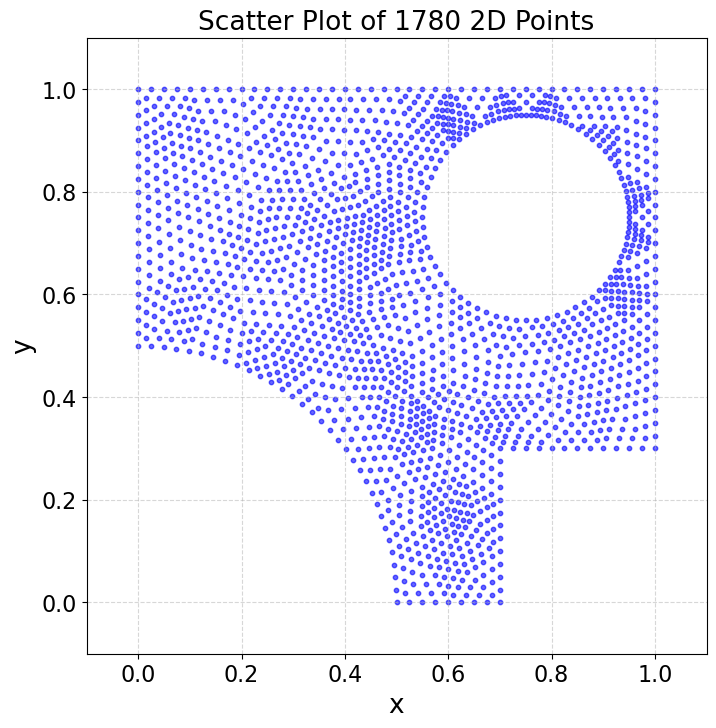}
        \caption*{(a) Sample points.}
    \end{minipage}
    \hfill
    \begin{minipage}[t]{0.3\textwidth}
        \centering
        \includegraphics[width=\linewidth]{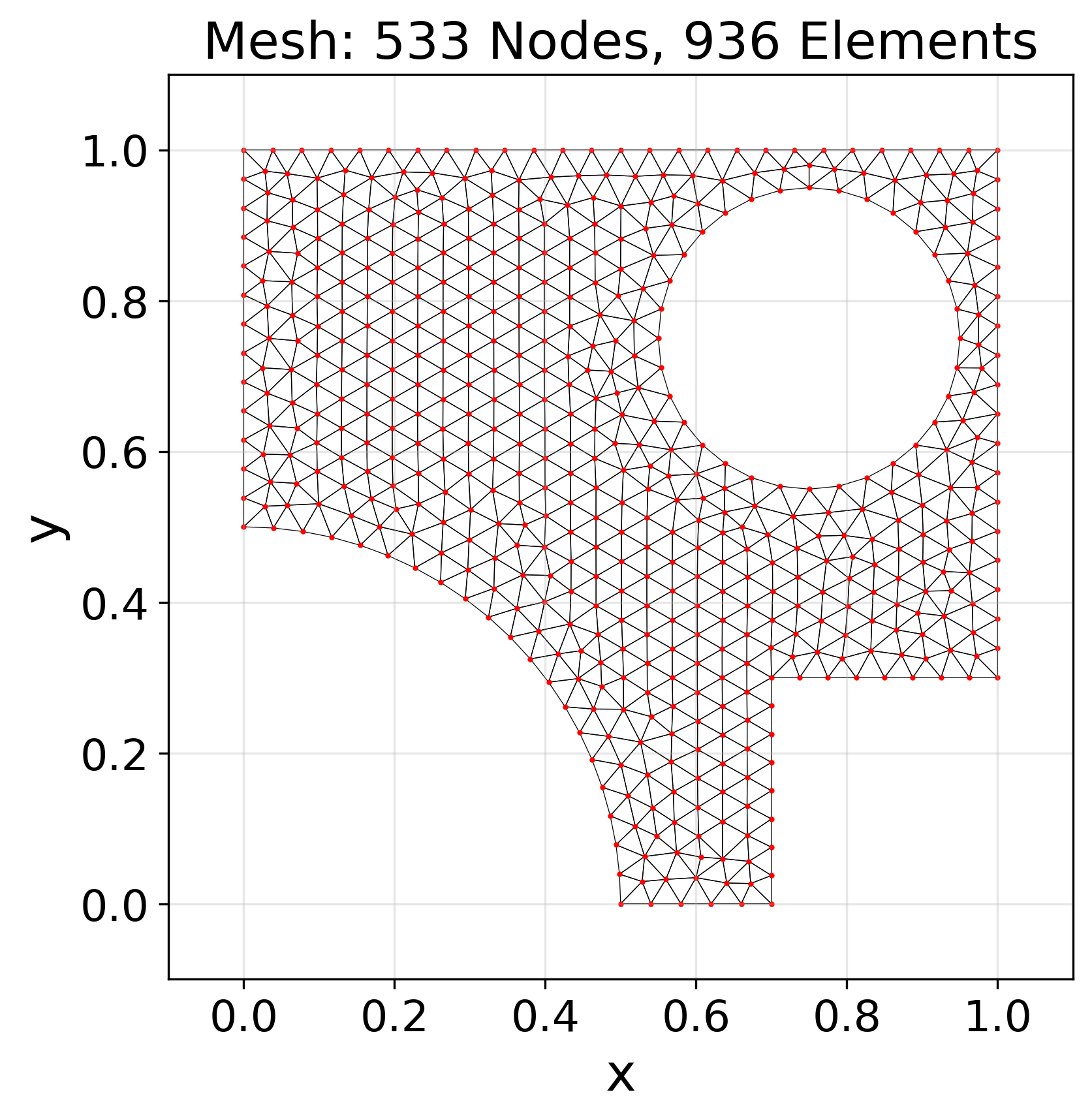}
        \caption*{(b) FEM mesh for encoder.}
    \end{minipage}
    \hfill
    \begin{minipage}[t]{0.3\textwidth}
        \centering
        \includegraphics[width=\linewidth]{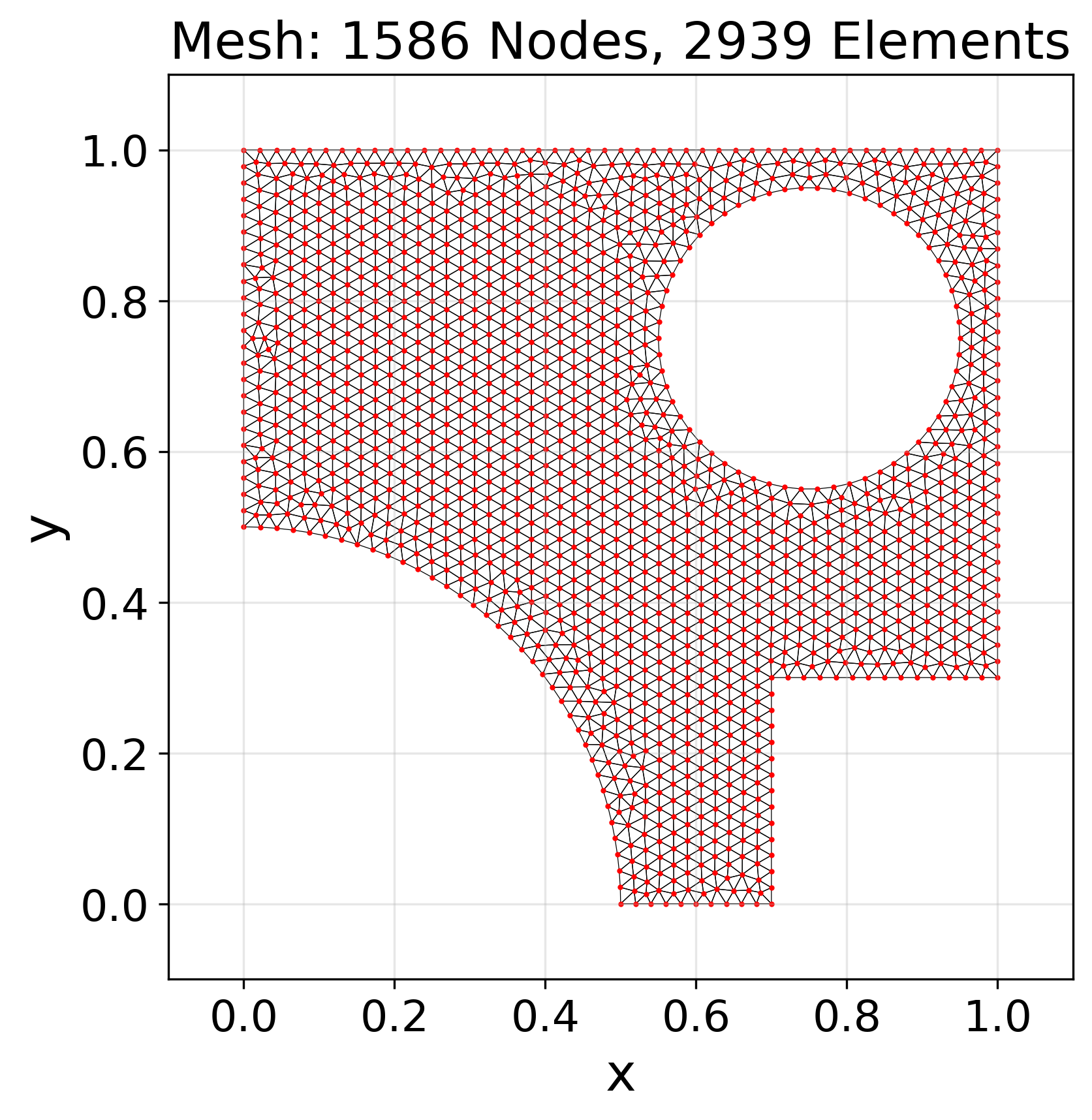}
        \caption*{(c) FEM mesh for decoder.}
    \end{minipage}
    \caption{
    \textbf{Complex-domain Poisson equation.}
    (a) Sample points used to evaluate the input and output fields.
    (b) FEM mesh used for the input encoder.
    (c) FEM mesh used for the output decoder.
    The domain contains three interior holes, so a boundary-conforming mesh provides useful geometric information for the prescribed FEM basis.
    }
    \label{fig:2d_Poi_point}
\end{figure}
We consider the two-dimensional Poisson equation on a complex domain \(\Omega\) with three interior holes:
\begin{equation}
\begin{aligned}
\Delta u(\mathbf{x}) &= f(\mathbf{x}), \qquad \mathbf{x}\in\Omega,\\
u(\mathbf{x}) &= 0, \qquad \mathbf{x}\in\partial\Omega .
\end{aligned}
\label{eq:poisson_complex_domain}
\end{equation}

The operator is
$
G:f(\mathbf{x})\mapsto u(\mathbf{x}).
$
The forcing term is sampled from
$
f(\mathbf{x})
=
\sum_{i=0}^{4}\sum_{j=0}^{4}
\alpha_{ij}
\sin\!\left([i\pi,j\pi]\cdot\mathbf{x}^{\top}\right)$, $
\alpha_{ij}\sim \mathcal{U}[-1,1].
$
Reference solutions are generated by FEM on a mesh with 1780 nodes. 
We use 1600 training samples and 400 testing samples. 
For FEM-C2C, the input and output basis dimensions are \((m_1,m_2)=(533,1586)\); for RFM-C2C, they are \((m_1,m_2)=(512,1600)\).

\begin{figure}[!ht]
    \centering
    \begin{minipage}[t]{0.33\textwidth}
        \centering
        \includegraphics[width=\linewidth]{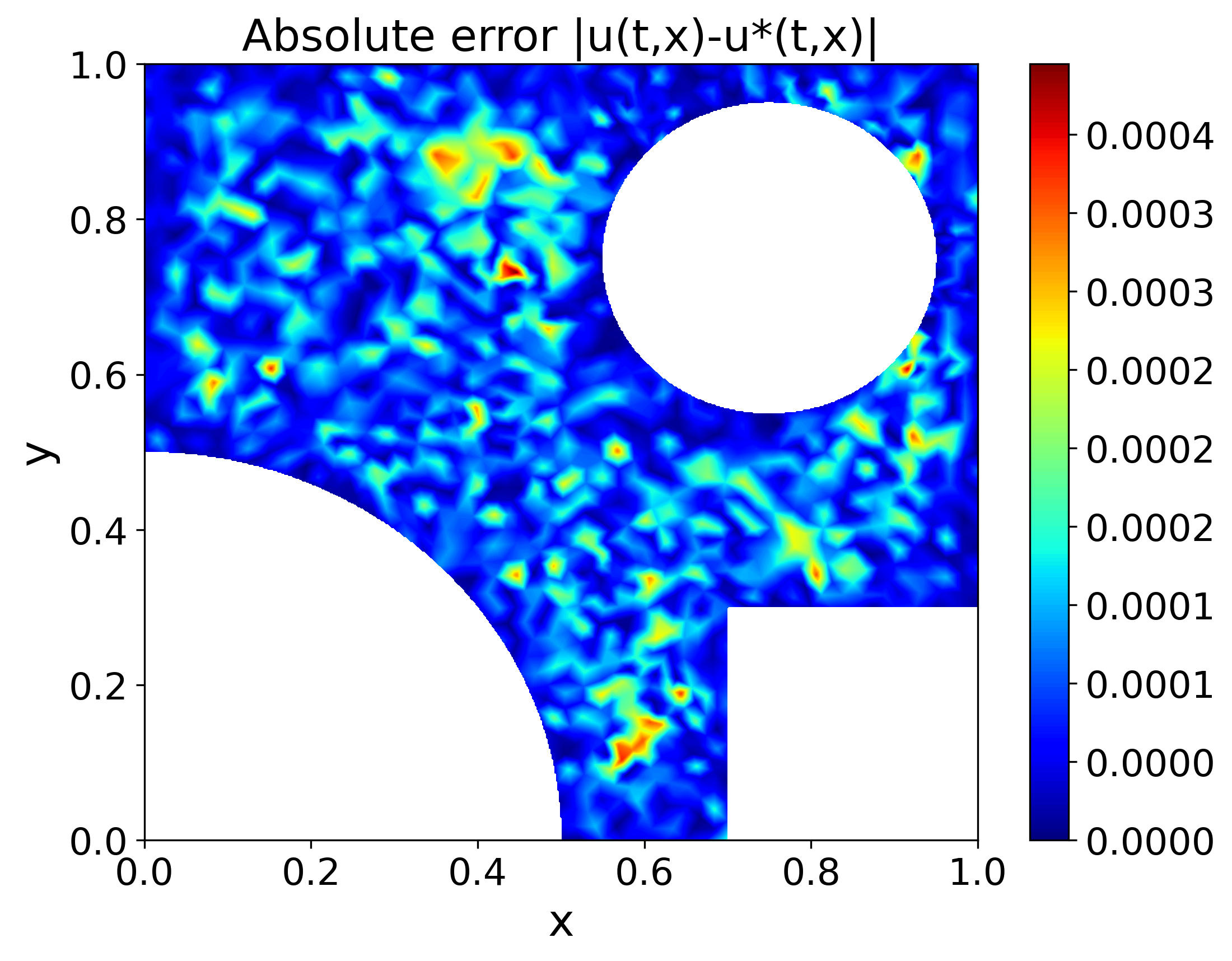}
        \caption*{(a) Pointwise error.}
    \end{minipage}
    \hfill
    \begin{minipage}[t]{0.33\textwidth}
        \centering
        \includegraphics[width=\linewidth]{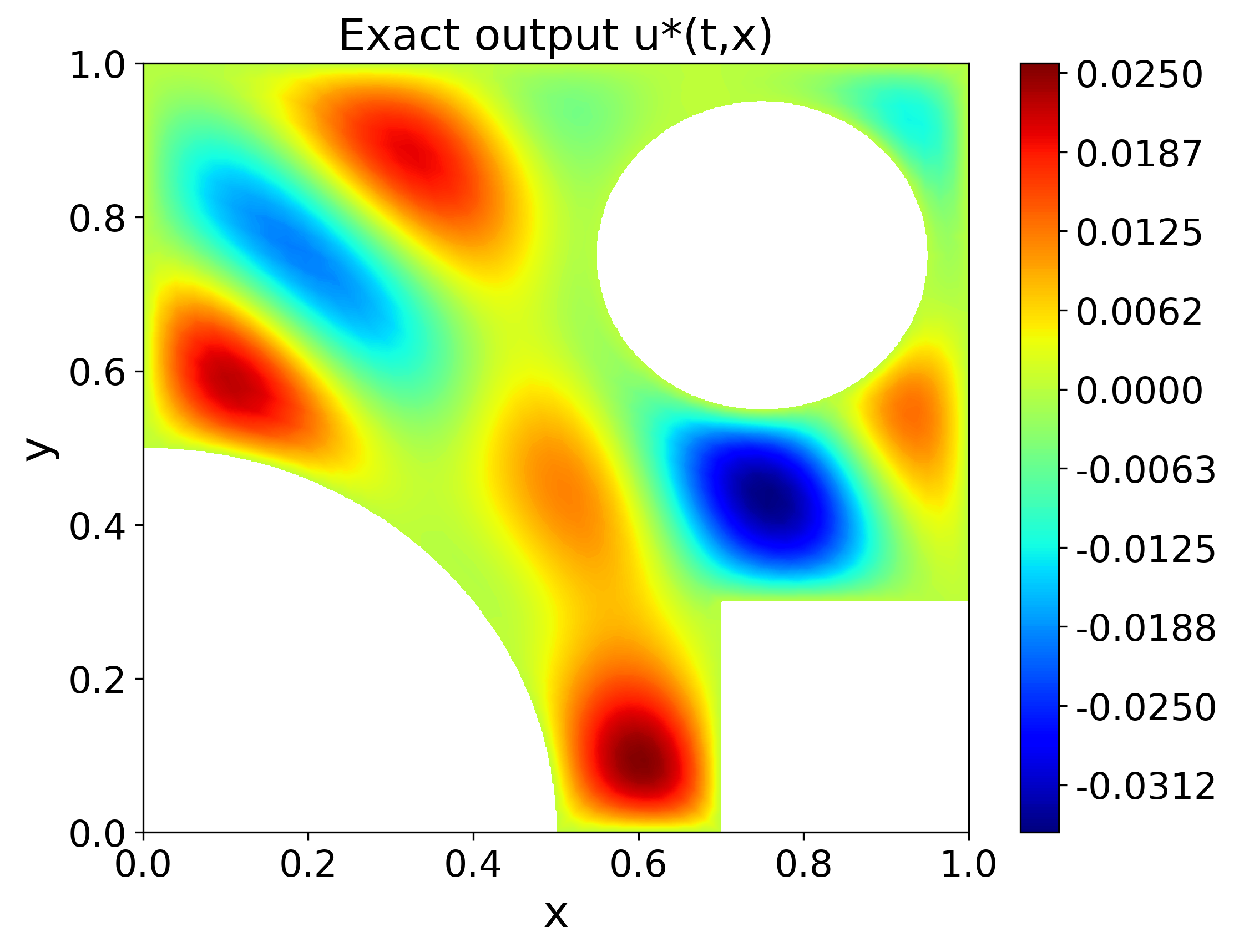}
        \caption*{(b) Exact solution.}
    \end{minipage}
    \hfill
    \begin{minipage}[t]{0.33\textwidth}
        \centering
        \includegraphics[width=\linewidth]{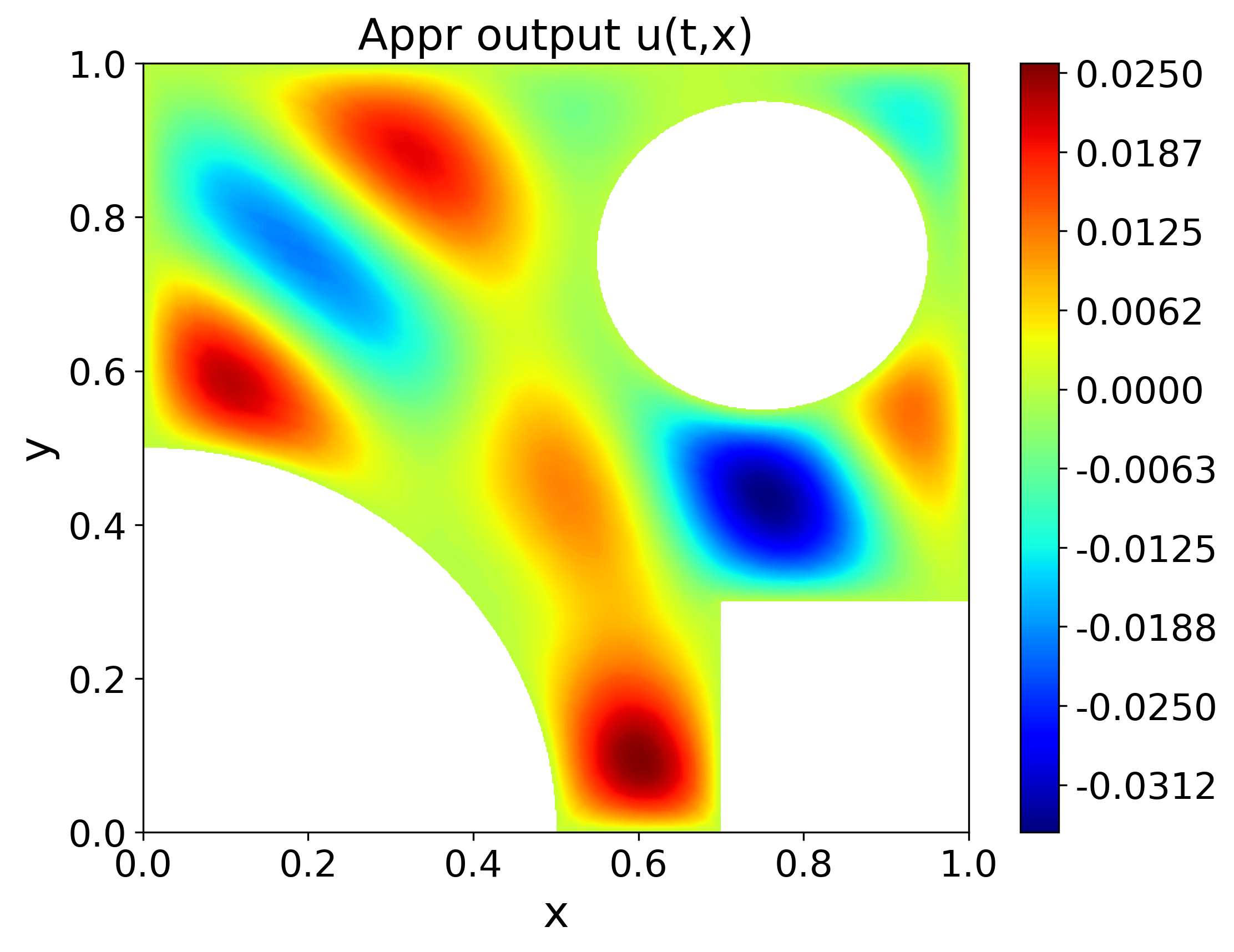}
        \caption*{(c) Predicted solution.}
    \end{minipage}
    \caption{
    \textbf{Complex-domain Poisson equation solved by FEM-C2C.}
    Worst-case test sample: (a) pointwise error, (b) exact solution, and (c) predicted solution.
    }
    \label{fig:2d_Poi_FEM_worst_case}
\end{figure}

\begin{figure}[!ht]
    \centering
    \begin{minipage}[t]{0.33\textwidth}
        \centering
        \includegraphics[width=\linewidth]{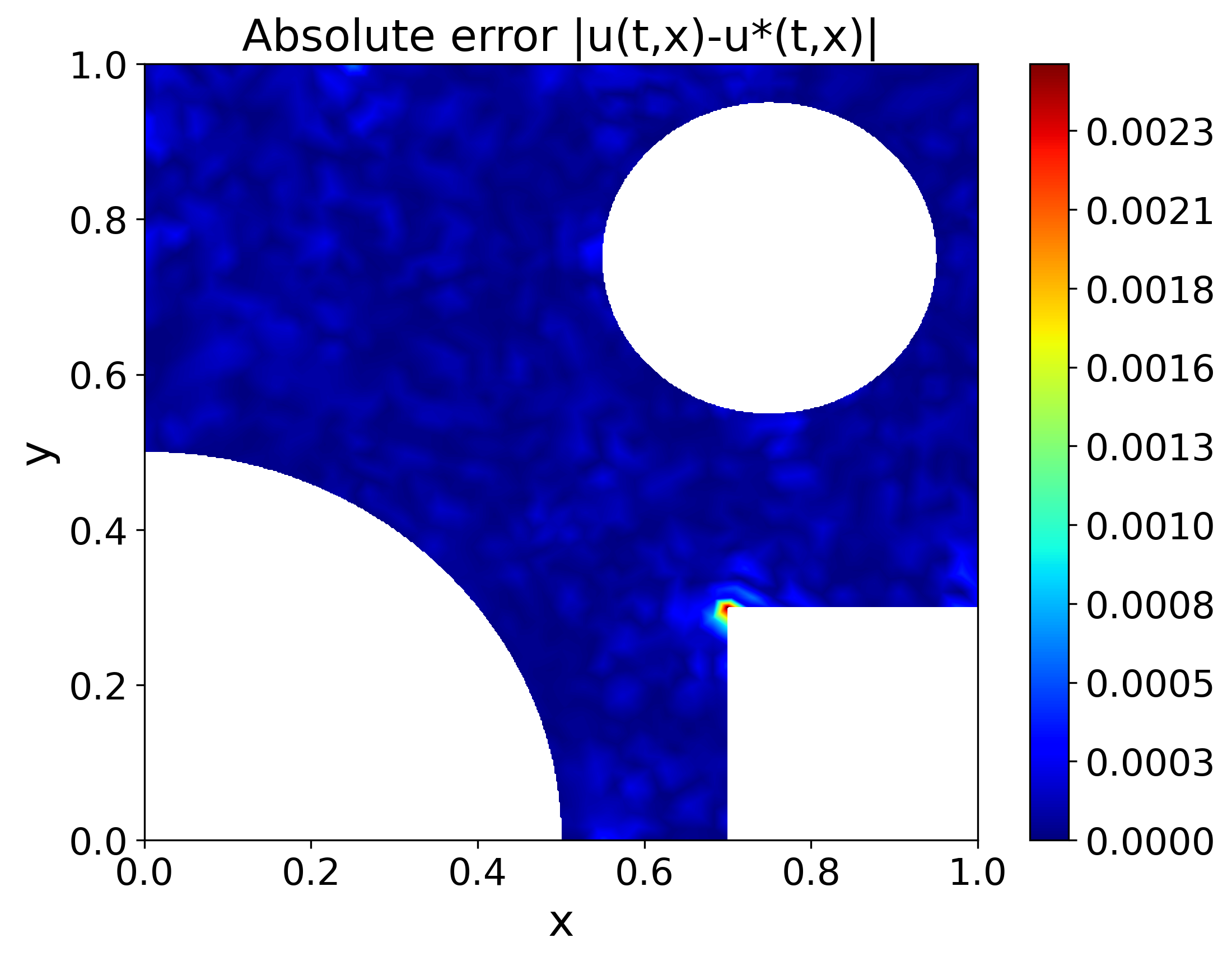}
        \caption*{(a) Pointwise error.}
    \end{minipage}
    \hfill
    \begin{minipage}[t]{0.33\textwidth}
        \centering
        \includegraphics[width=\linewidth]{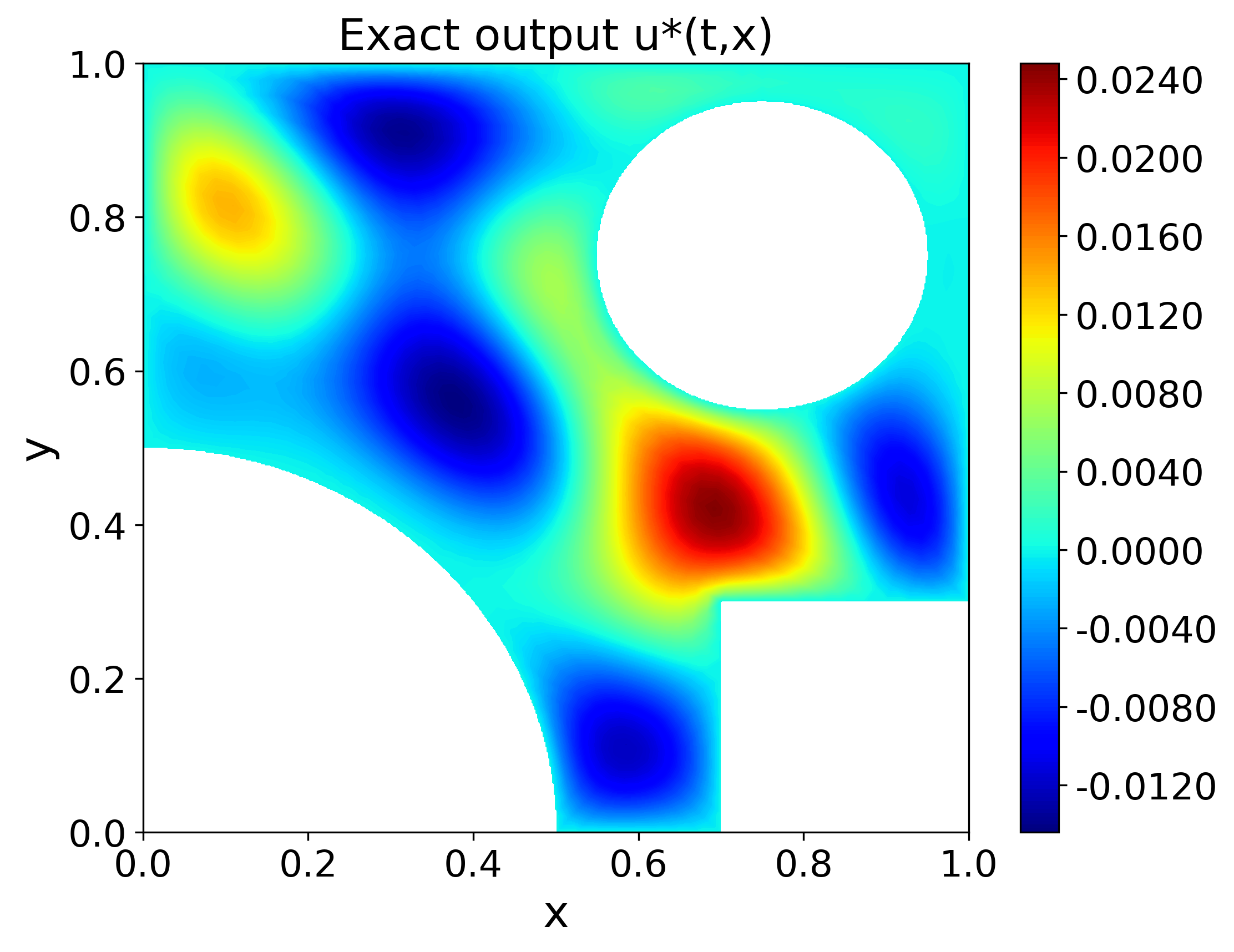}
        \caption*{(b) Exact solution.}
    \end{minipage}
    \hfill
    \begin{minipage}[t]{0.33\textwidth}
        \centering
        \includegraphics[width=\linewidth]{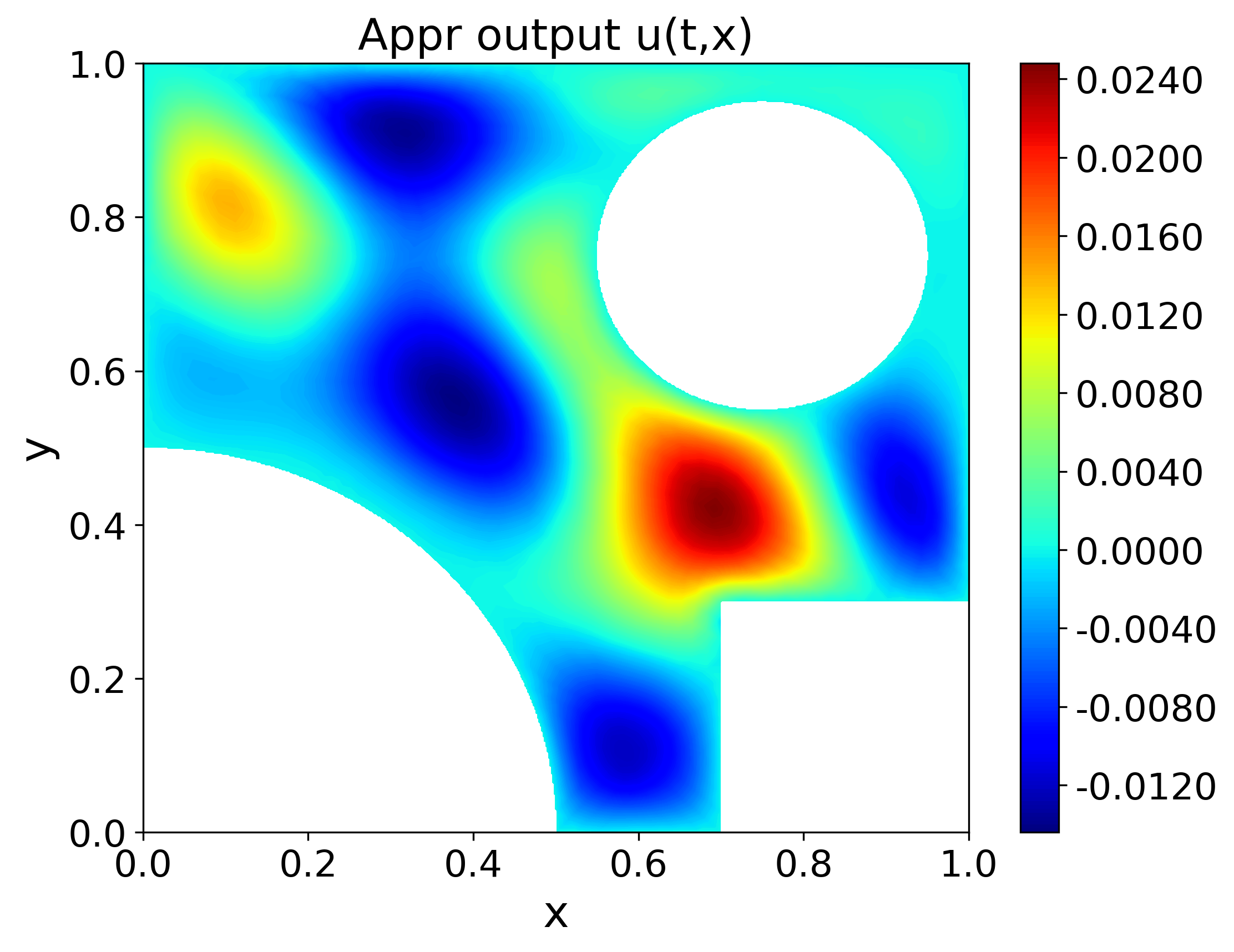}
        \caption*{(c) Predicted solution.}
    \end{minipage}
    \caption{
    \textbf{Complex-domain Poisson equation solved by RFM-C2C.}
    Worst-case test sample: (a) pointwise error, (b) exact solution, and (c) predicted solution.
    }
    \label{fig:2d_Poi_RFM_worst_case}
\end{figure}

The test RL2E values are 4.13e-3 for FEM-C2C and 7.69e-3 for RFM-C2C. 
The better performance of FEM-C2C reflects the advantage of a boundary-conforming basis on this complex geometry. 
Figure~\ref{fig:2d_Poi_point} shows the sampling points and FEM meshes, while Figures~\ref{fig:2d_Poi_FEM_worst_case} and~\ref{fig:2d_Poi_RFM_worst_case} show worst-case predictions for the two bases.

\begin{figure}[!ht]
    \centering
    \begin{minipage}[t]{0.4\textwidth}
        \centering
        \includegraphics[width=\linewidth]{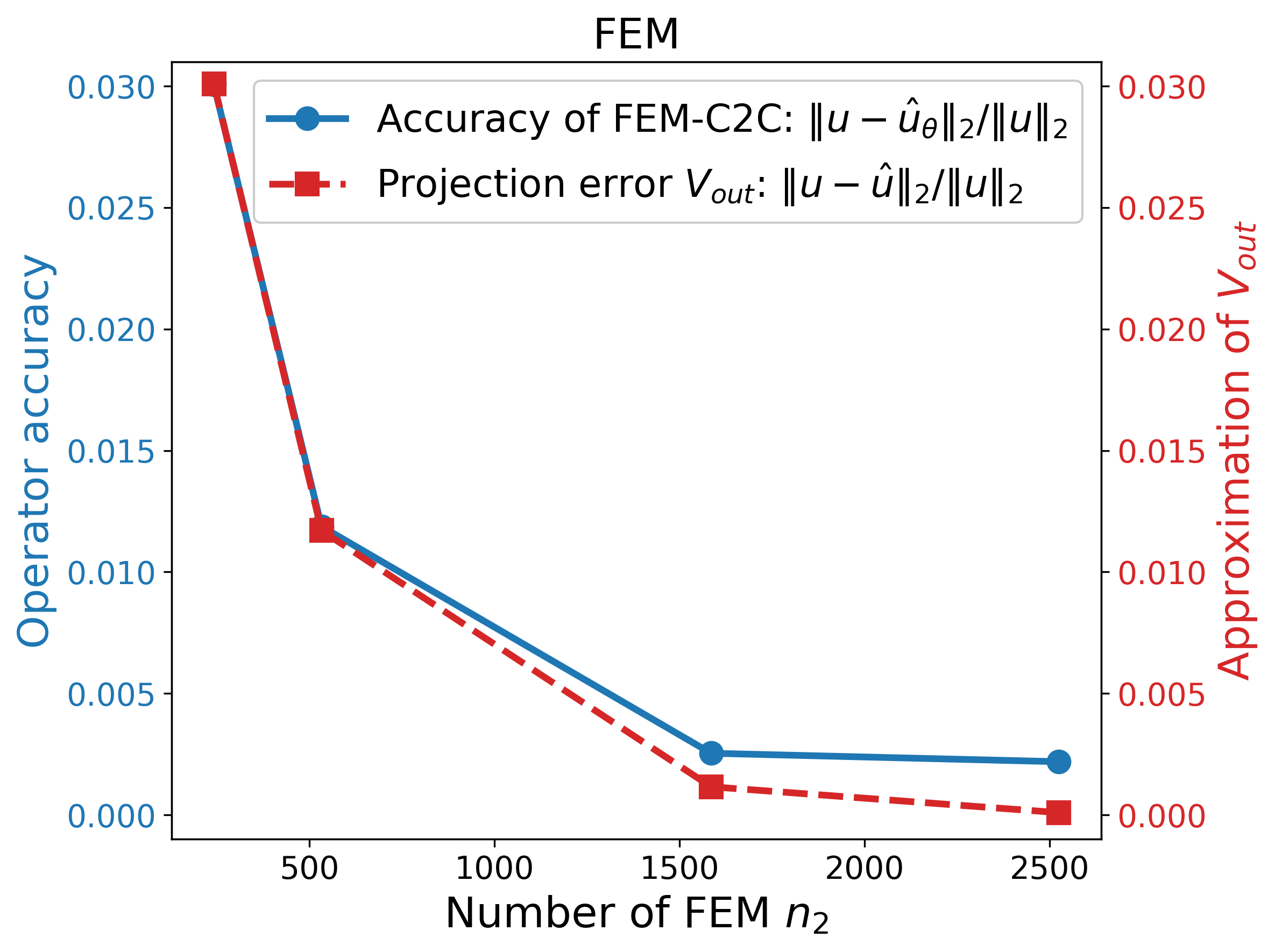}
        \caption*{(a) FEM basis.}
    \end{minipage}
\hspace{0.1\textwidth}
    \begin{minipage}[t]{0.4\textwidth}
        \centering
        \includegraphics[width=\linewidth]{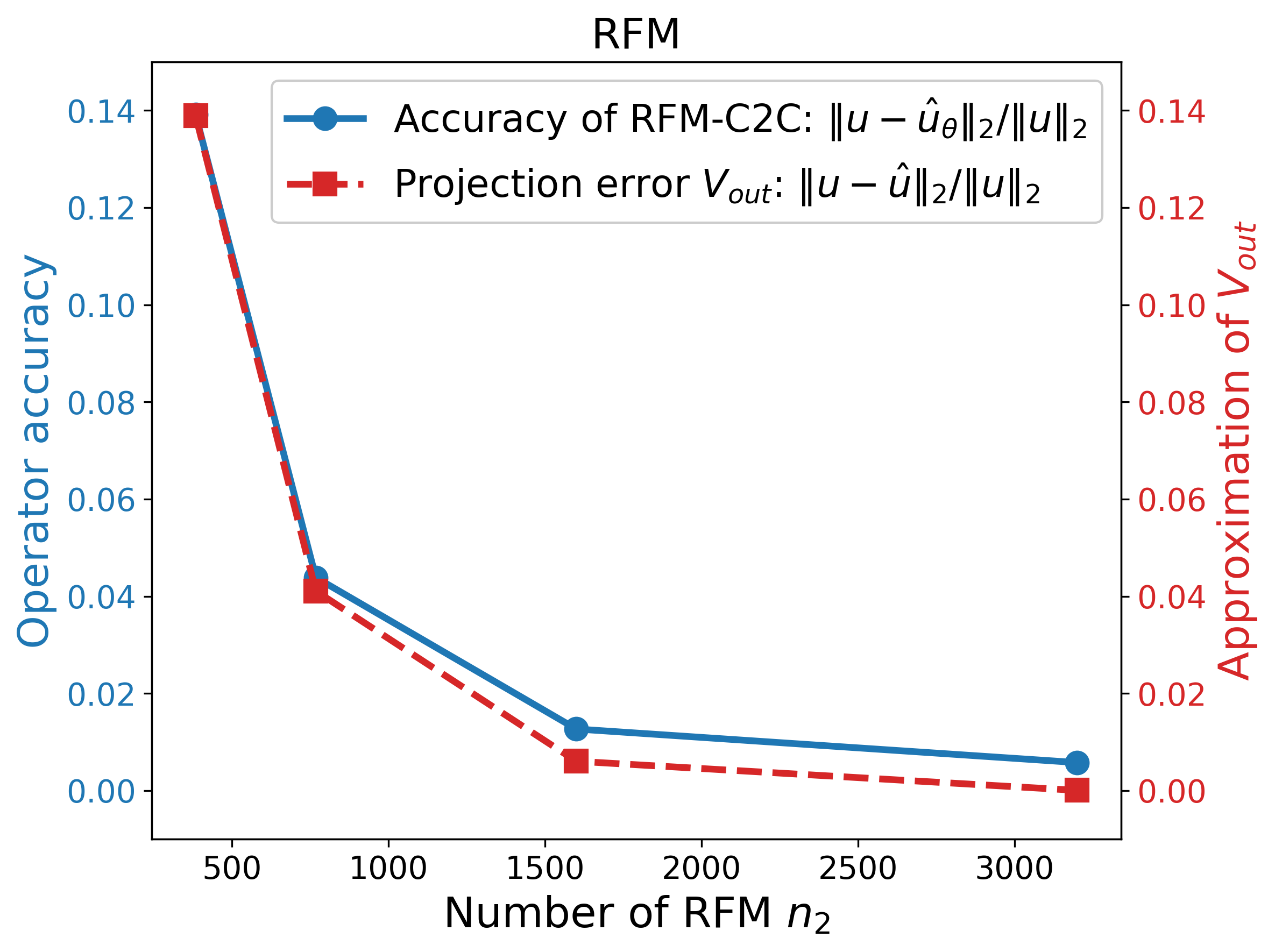}
        \caption*{(b) RFM basis.}
    \end{minipage}
    \caption{
    \textbf{Output projection error on the complex-domain Poisson problem.}
    The blue curve shows the FB-C2C prediction error, while the red dashed curve shows the best approximation error in the prescribed output basis space.
    The prediction error is bounded below by the output projection error, and both decrease as the number of output basis functions \(m_2\) increases.
    }
    \label{fig:2d_Poi_error}
\end{figure}

We next isolate the effect of the prescribed output space. 
For each output basis size \(m_2\), we compare the FB-C2C prediction error with the best approximation error attainable in the same output space. 
As shown in Figure~\ref{fig:2d_Poi_error}, the learned-operator error remains above the output projection error, consistent with the error decomposition in Section~3. 
Both errors decrease as \(m_2\) increases, and the two curves are close in this experiment. 
This suggests that, for this experiment, the coefficient-map learning error is small relative to the approximation error imposed by the prescribed output space.

\newpage
\subsubsection{Multi-component PDE operators}
\label{subsubsec:multi-component-operators}

We first consider PDE operators with multiple input or output components. 
The elastic plate example maps one input loading function to two displacement fields, while the L-shaped Darcy example maps two input fields to one output solution. 
These examples show that FB-C2CNet can handle multi-component operators through either concatenated scalar-valued coefficient representations or vector-valued basis functions.

\paragraph{Elastic plate: one input and two output fields}

\begin{figure}[!ht]
    \centering
    \begin{minipage}[t]{0.3\textwidth}
        \centering
        \includegraphics[width=\linewidth]{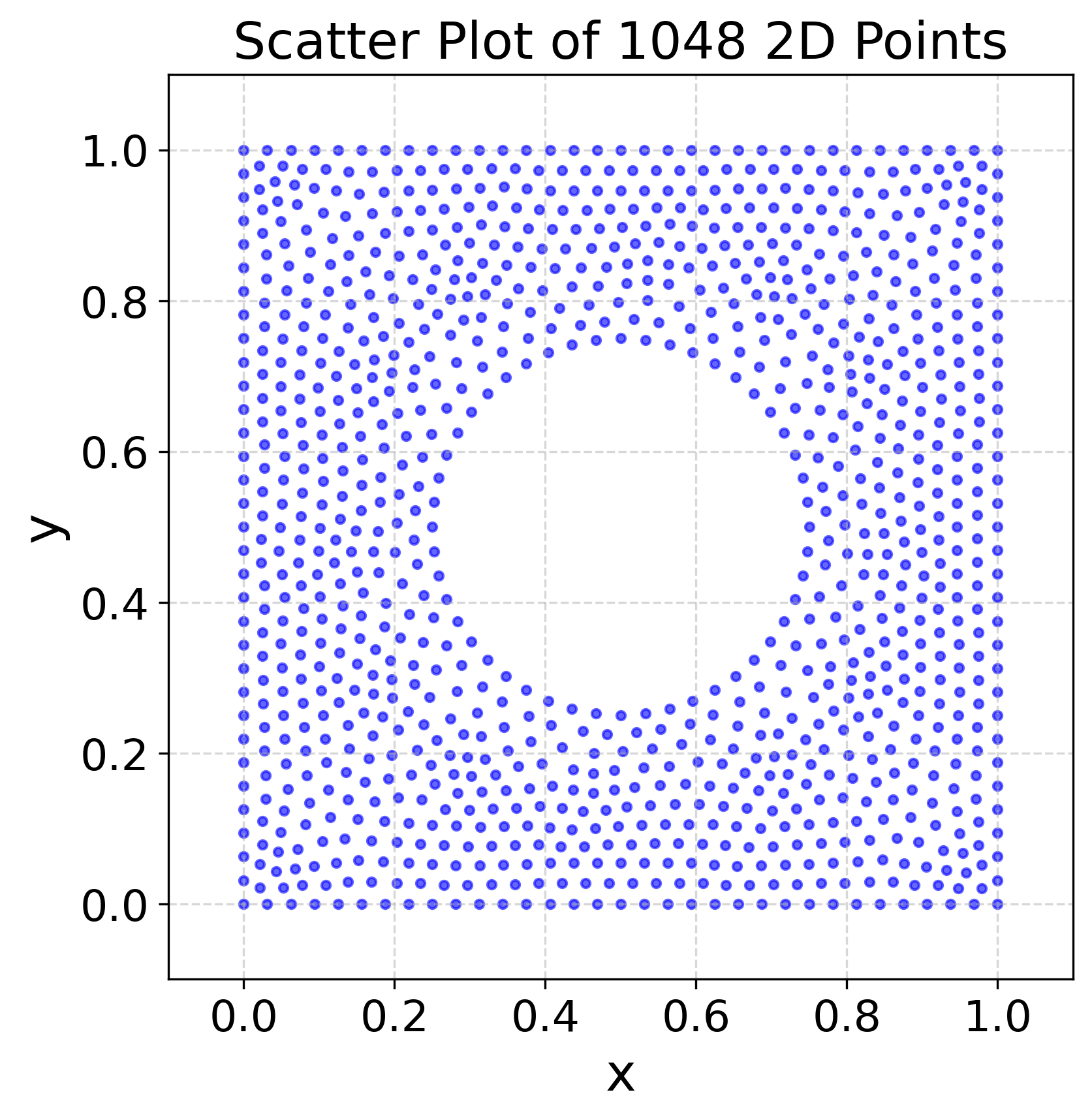}
        \caption*{(a) Output sample points.}
    \end{minipage}
    \hfill
    \begin{minipage}[t]{0.3\textwidth}
        \centering
        \includegraphics[width=\linewidth]{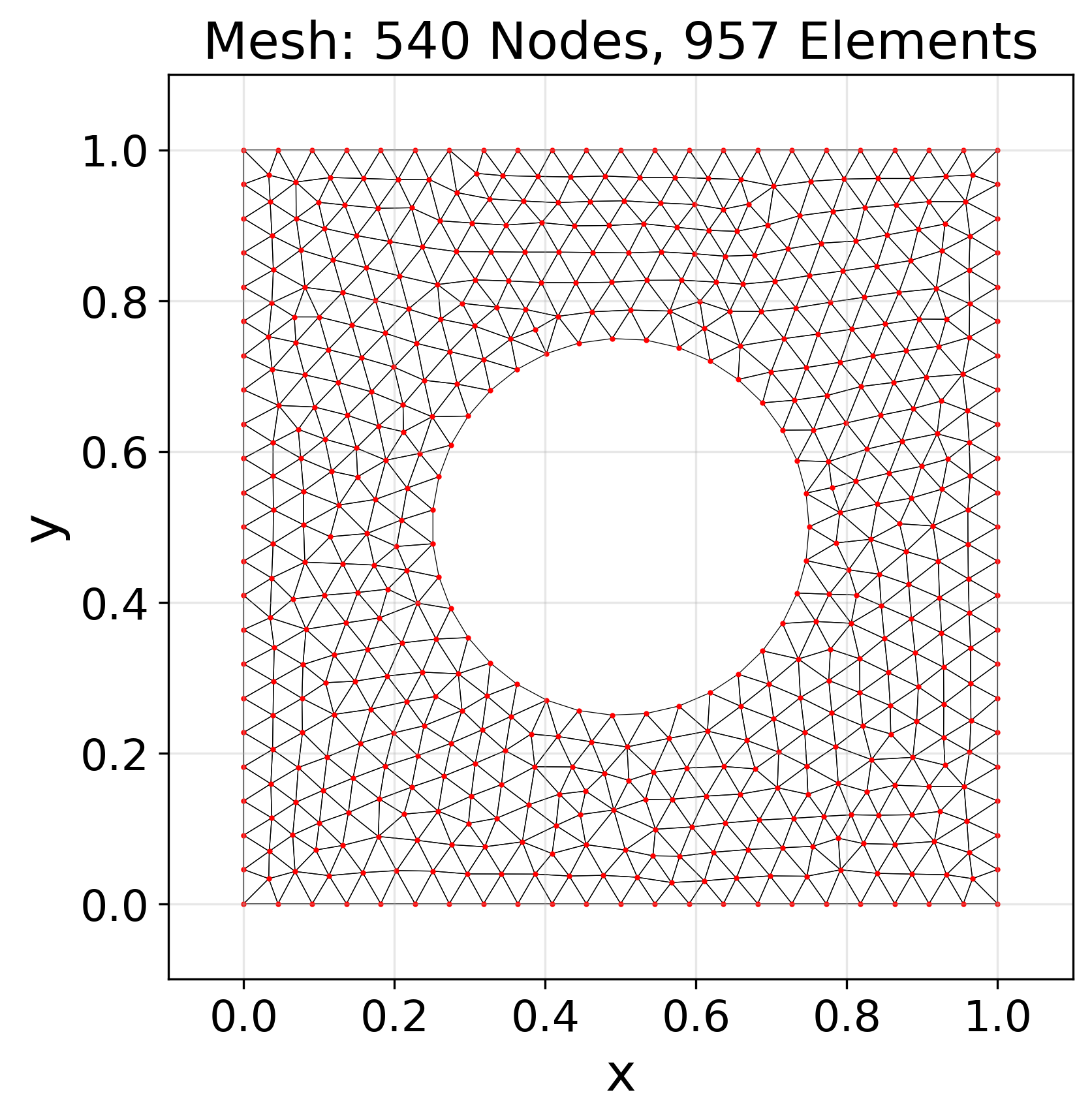}
        \caption*{(b) \(\psi^v\) FEM mesh.}
    \end{minipage}
    \hfill
    \begin{minipage}[t]{0.3\textwidth}
        \centering
        \includegraphics[width=\linewidth]{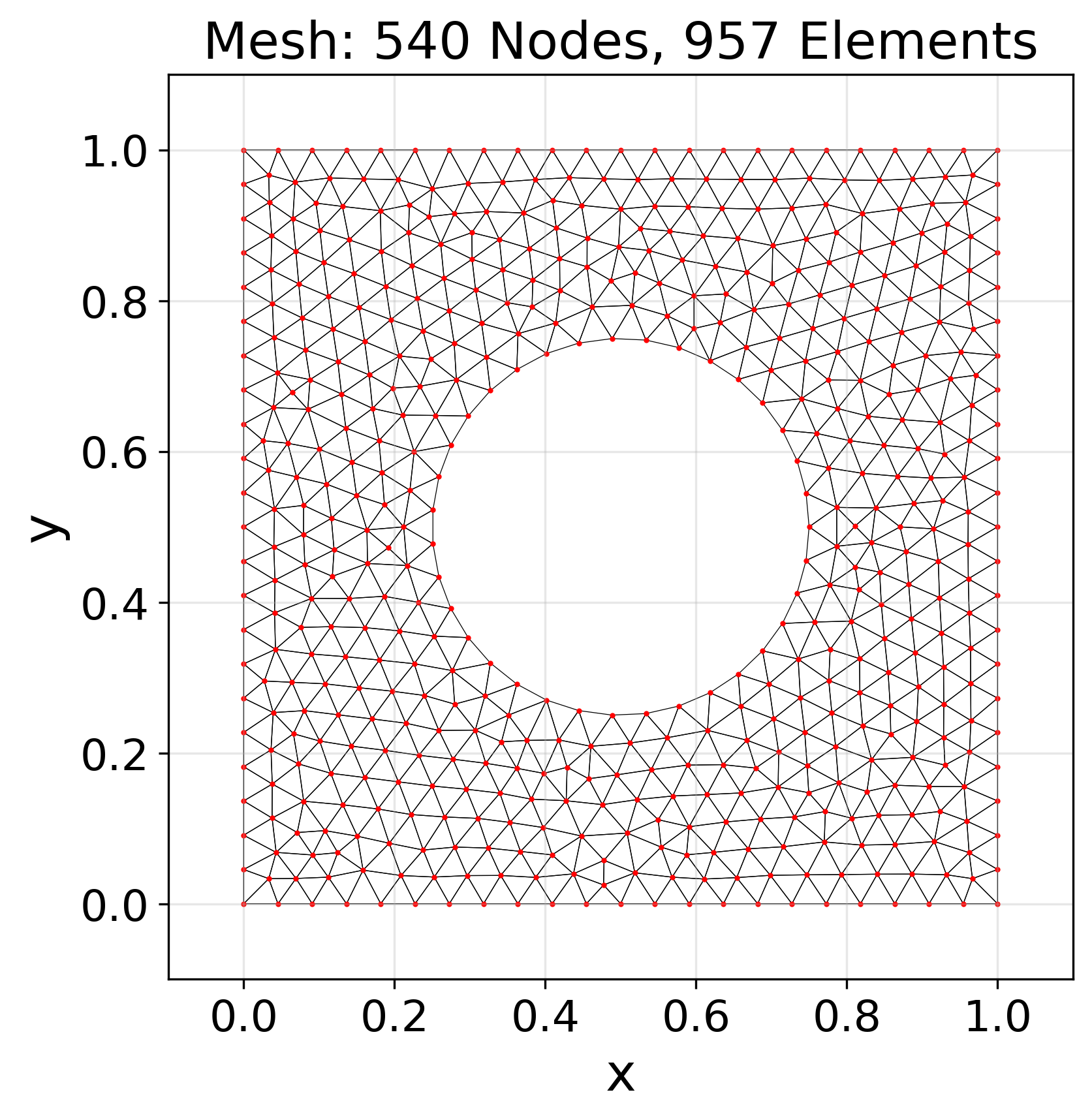}
        \caption*{(c) \(\psi^u\) FEM mesh, also used for \(\psi\).}
    \end{minipage}
    \caption{
    \textbf{Elastic plate.}
    Output sample points and the FEM meshes used to construct scalar- and vector-valued output bases. 
    }
    \label{fig:elasticplateoutputgridandelement}
\end{figure}

We consider a thin rectangular plate under in-plane loading, modeled by the two-dimensional plane-stress elasticity equations on a domain \(\Omega=(0,1)\times(0,1)\) with an interior cut-out. 
The equilibrium equation is
\begin{equation}
\nabla\cdot\boldsymbol\sigma+\mathbf b(\mathbf x)=\mathbf 0,
\qquad \mathbf x\in\Omega,
\label{eq:elastic-equilibrium}
\end{equation}
where \(\mathbf b=\mathbf 0\). 
The constitutive relation is
\begin{equation}
\begin{bmatrix}
\sigma_{xx}\\
\sigma_{yy}\\
\tau_{xy}
\end{bmatrix}
=
\frac{E}{1-\nu^2}
\begin{bmatrix}
1&\nu&0\\
\nu&1&0\\
0&0&(1-\nu)/2
\end{bmatrix}
\begin{bmatrix}
\epsilon_{xx}\\
\epsilon_{yy}\\
\gamma_{xy}
\end{bmatrix},
\qquad
E=300\times 10^5,\quad \nu=0.3,
\label{eq:elastic-constitutive}
\end{equation}
with
$
\epsilon_{xx}=\frac{\partial u}{\partial x},
\epsilon_{yy}=\frac{\partial v}{\partial y},
\gamma_{xy}=\frac{\partial u}{\partial y}+\frac{\partial v}{\partial x}.
$ The plate is fixed on the left edge, and a traction \(\mathbf f(\mathbf x)=[0,f(1,y)]^\top\) is applied on the right edge. 
The input loading \(f(1,y)\) is modeled as a Gaussian random field, and the target operator is $
G:f(1,y)\mapsto [u(\mathbf x),v(\mathbf x)].
$ The dataset is from~\cite{ingebrand2025basis}, with 1850 training samples and 100 testing samples. 
The input uses 101 grid points, and the output grid contains 1048 points, as shown in Figure~\ref{fig:elasticplateoutputgridandelement}.

For multi-output approximation, we compare scalar-valued and vector-valued output bases. 
In the scalar-valued case, the two displacement components are represented separately:
\begin{equation}
u(\mathbf x)\approx\sum_{j=1}^{m_2}b_j^u\psi_j(\mathbf x),
\qquad
v(\mathbf x)\approx\sum_{j=1}^{m_2}b_j^v\psi_j(\mathbf x).
\label{eq:elastic-scalar-basis}
\end{equation}
In the vector-valued case, we use
\[
V_{\mathrm{out}}^{m_2}
=
\operatorname{span}\{\boldsymbol\psi_j(\mathbf x)\}_{j=1}^{m_2},
\qquad
\boldsymbol\psi_j(\mathbf x)=[\psi_j^u(\mathbf x),\psi_j^v(\mathbf x)],
\]
and approximate
\begin{equation}
[u(\mathbf x),v(\mathbf x)]
\approx
\sum_{j=1}^{m_2}b_j\boldsymbol\psi_j(\mathbf x).
\label{eq:elastic-vector-basis}
\end{equation}

\begin{figure}[!ht]
    \centering
    \begin{minipage}[t]{0.24\textwidth}
        \centering
        \includegraphics[width=\linewidth]{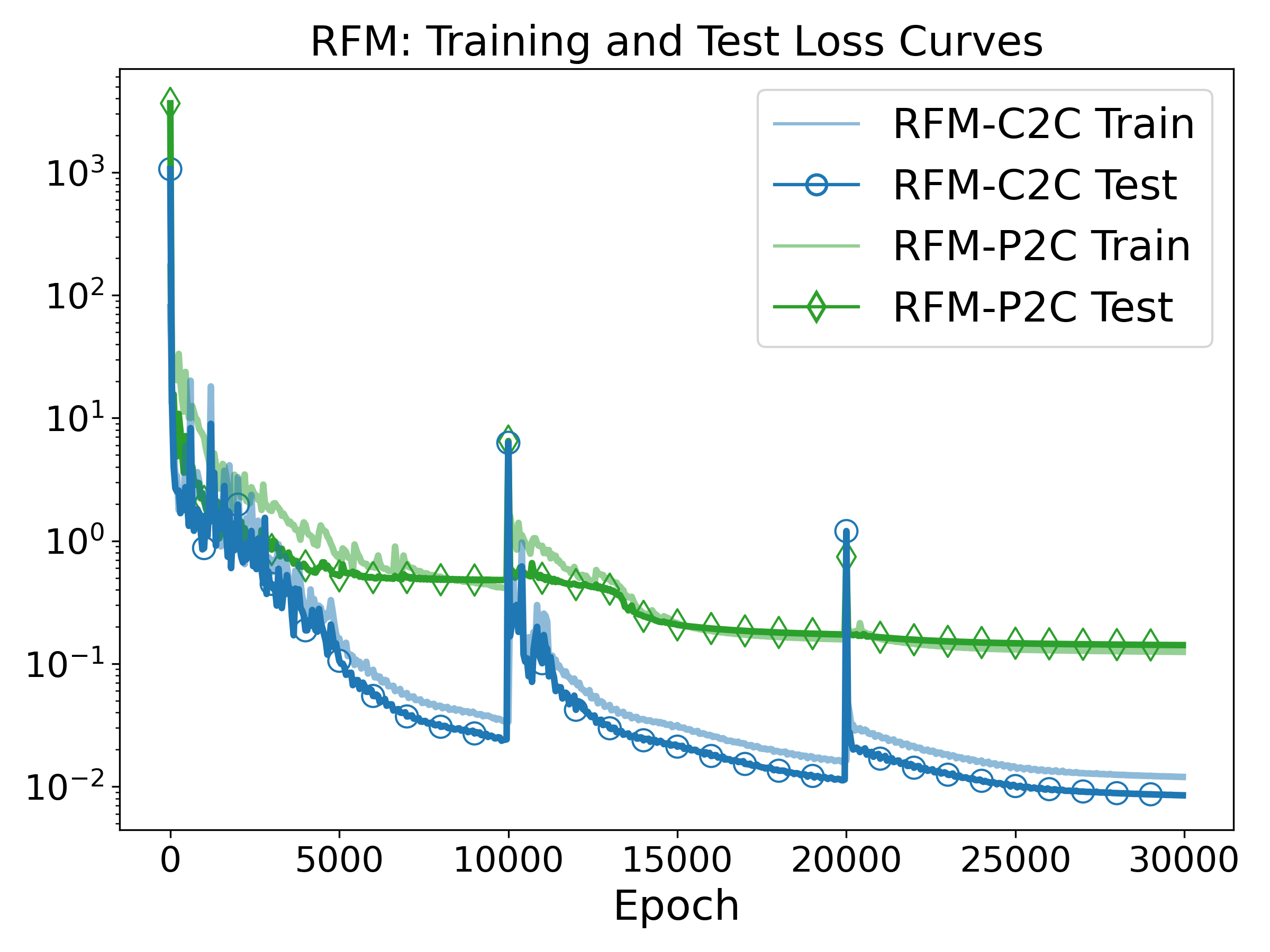}
        \caption*{(a) Scalar-valued RFM basis.}
    \end{minipage}
    \hfill
    \begin{minipage}[t]{0.24\textwidth}
        \centering
        \includegraphics[width=\linewidth]{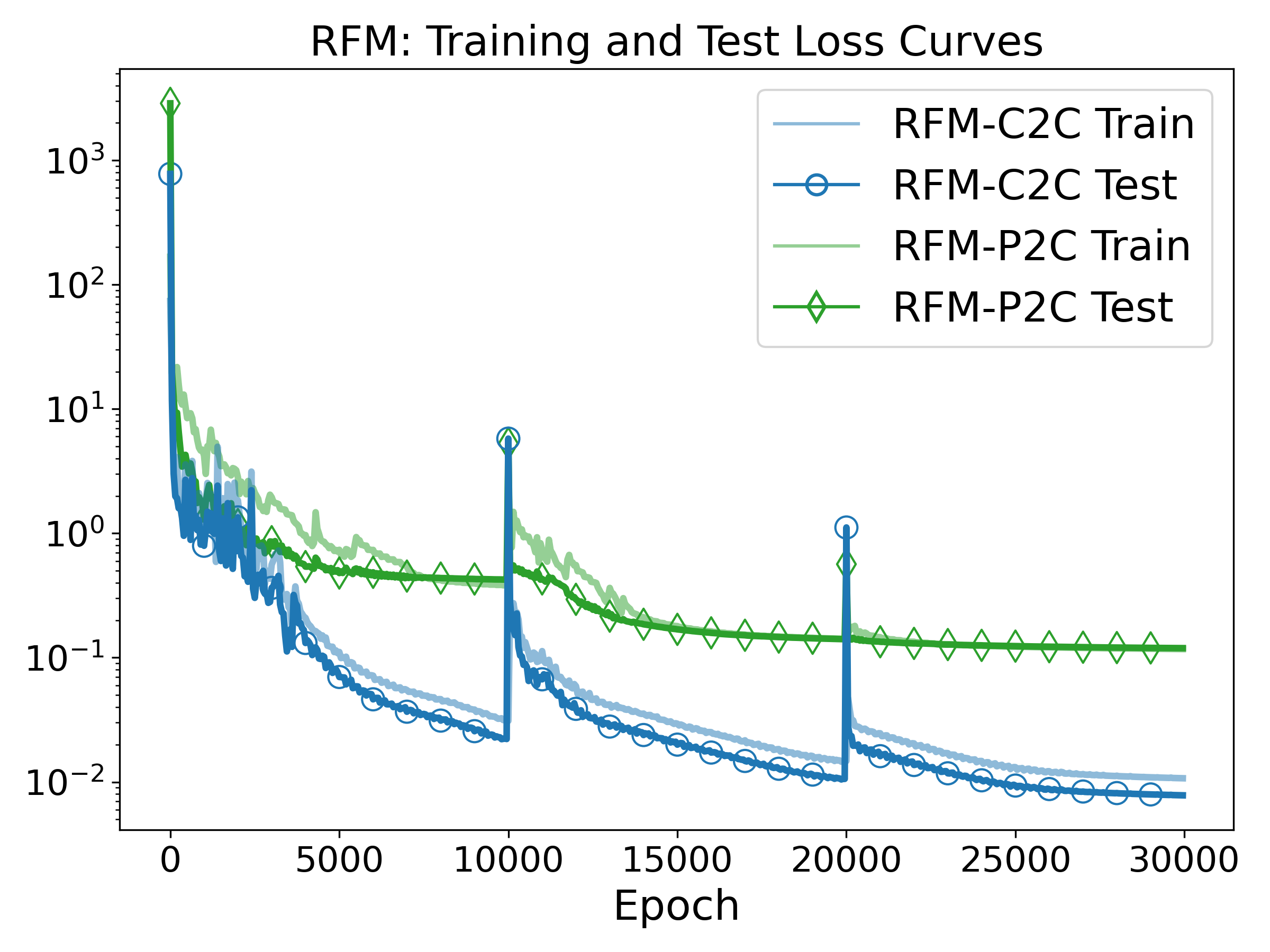}
        \caption*{(b) Vector-valued RFM basis.}
    \end{minipage}
    \hfill
    \begin{minipage}[t]{0.24\textwidth}
        \centering
        \includegraphics[width=\linewidth]{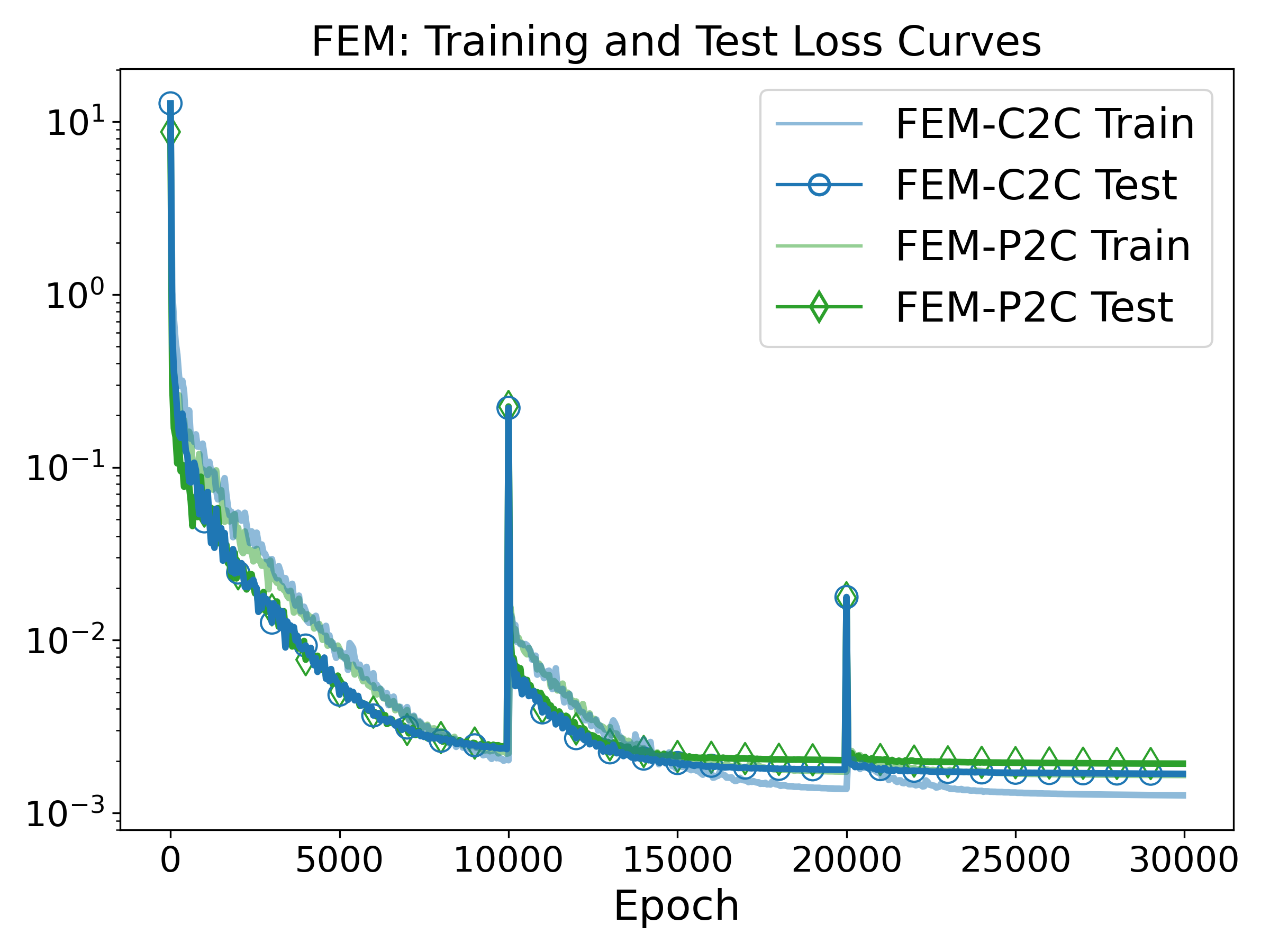}
        \caption*{(c) Scalar-valued FEM basis.}
    \end{minipage}
    \hfill
    \begin{minipage}[t]{0.24\textwidth}
        \centering
        \includegraphics[width=\linewidth]{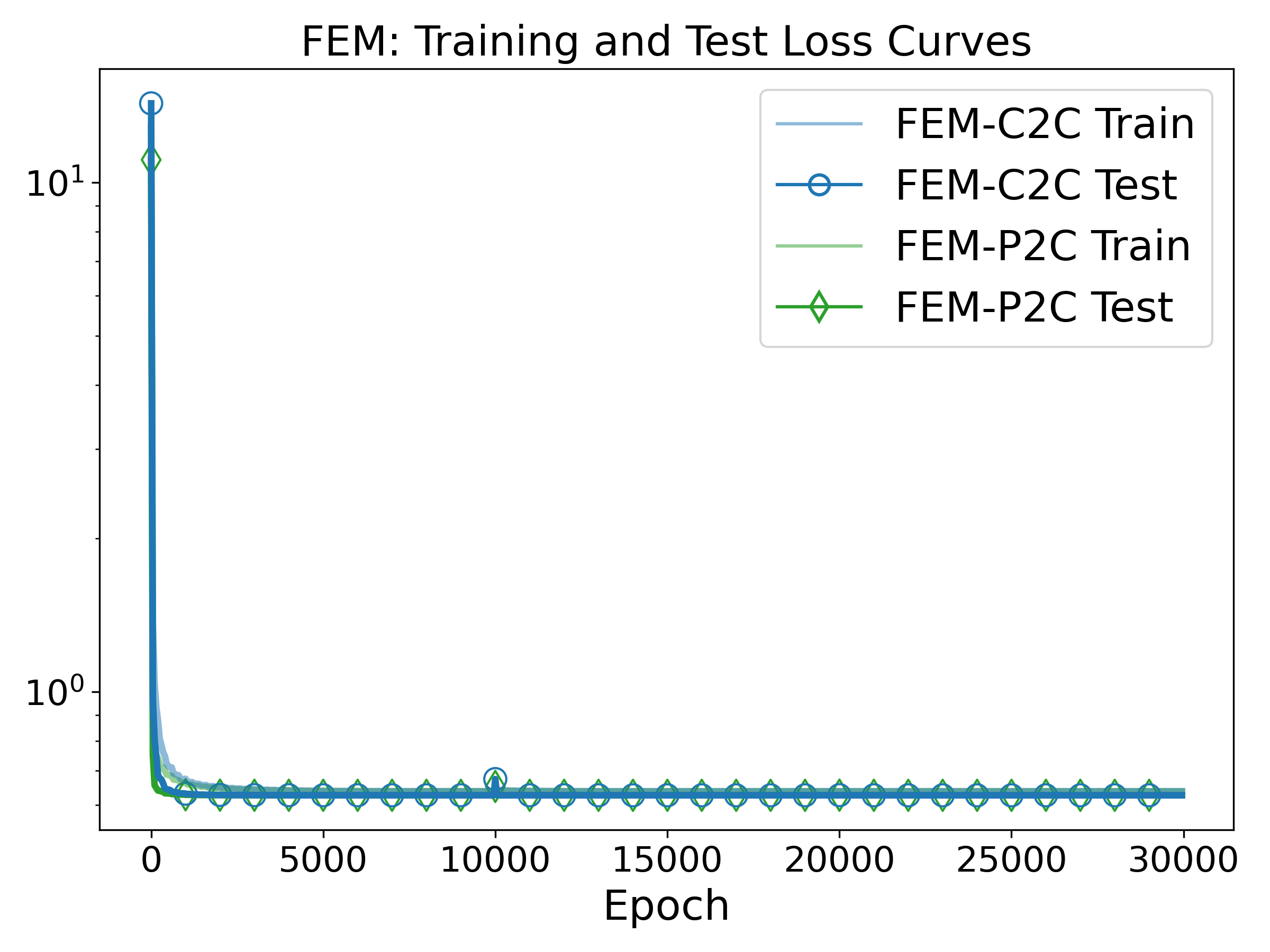}
        \caption*{(d) Vector-valued FEM basis.}
    \end{minipage}
    \caption{
    \textbf{Elastic plate with RFM bases.}
    Training curves for scalar- and vector-valued RFM and FEM output bases.
    %Both variants achieve similar accuracy, while the vector-valued representation uses fewer output coefficients.
    }
    \label{fig:elastic_curve_RFM}
\end{figure}

The vector-valued representation uses fewer output coefficients, but its accuracy depends on whether the chosen vector-valued basis provides a good joint approximation of the two displacement components. In this example, RFM attains $m_2=512$ while FEM yields $m_2=540$, as the latter cannot generate a mesh with exactly 512 basis functions. In the vector-valued FEM case, we construct different FEM bases $\{\psi_i^u\}$ and $\{\psi_i^v\}$ for the displacement components $u$ and $v$, respectively. Both bases have $m_2=540$ and are shown in Figure~\ref{fig:elasticplateoutputgridandelement}(b)--(c). In the scalar-valued FEM case, we take $\psi_i=\psi_i^u$. Experimental results and training curves are shown in Table~\ref{tab:2d_elasticity_time} and Figure~\ref{fig:elastic_curve_RFM}. 

\begin{table}[!htpb]
    \centering
    \caption{
    Elastic plate.
    Comparison of C2C and P2C using scalar- and vector-valued RFM and FEM output bases. 
    Within each basis, the lower C2C/P2C training time and test error are shown in bold.
    }
    \label{tab:2d_elasticity_time}
    \footnotesize
    \renewcommand{\arraystretch}{1.12}
    \setlength{\tabcolsep}{3pt}
    \begin{tabular}{@{}lccccccc@{}}
        \toprule
        & \multicolumn{1}{c}{Encoding time}
        & \multicolumn{2}{c}{Training time}
        & \multicolumn{2}{c}{Test RL2E}
        & \multicolumn{2}{c}{Network architecture}\\
        \cmidrule(lr){2-2}\cmidrule(lr){3-4}\cmidrule(lr){5-6}\cmidrule(l){7-8}
        Basis & Encoder & C2C & P2C & C2C & P2C & C2C & P2C \\ 
        \midrule
        Scalar RFM 
        & 0.82 s
        & \textbf{43.40 s}
        & 48.11 s
        & \textbf{8.47e-3}
        & 1.42e-1 
        & [64,512,1024] 
        & [101,512,1024]\\ 
        Scalar FEM 
        & 11.75 s
        & \textbf{49.01 s}
        & 52.80 s
        & \textbf{1.68e-3} 
        & 1.93e-3 
        & [64,512,1080] 
        & [101,512,1080]\\ 
        Vector RFM 
        & 0.82 s
        & \textbf{42.11 s}
        & 64.86 s
        & \textbf{7.81e-3}
        & 1.41e-1 
        & [64,512,512]
        & [101,512,512]\\ 
        Vector FEM 
        & 11.75 s
        & \textbf{40.22 s}
        & 42.84 s
        & \textbf{6.26e-1}
        & \textbf{6.26e-1}
        & [64,512,540] 
        & [101,512,540]\\ 
        \bottomrule
    \end{tabular}
\end{table}

% \begin{figure}[!ht]
%     \centering
%     \begin{minipage}[t]{0.3\textwidth}
%         \centering
%         \includegraphics[width=\linewidth]{2D_Elasticity_input_1048.png}
%         \caption*{(a) Output sample points.}
%     \end{minipage}
%     \hfill
%     \begin{minipage}[t]{0.3\textwidth}
%         \centering
%         \includegraphics[width=\linewidth]{2D_Elasticity_FEM_mesh1.png}
%         \caption*{(b) \(\psi^v\) FEM mesh.}
%     \end{minipage}
%     \hfill
%     \begin{minipage}[t]{0.3\textwidth}
%         \centering
%         \includegraphics[width=\linewidth]{2D_Elasticity_FEM_mesh2.png}
%         \caption*{(c) \(\psi^u\) FEM mesh.}
%     \end{minipage}
%     \caption{
%     \textbf{Elastic plate.}
%     Output sample points and the two FEM meshes used to construct scalar-valued output bases for the displacement components.
%     }
%     \label{fig:elasticplateoutputgridandelement}
% \end{figure}

% \begin{figure}[!ht]
%     \centering
%     \begin{minipage}[t]{0.4\textwidth}
%         \centering
%         \includegraphics[width=\linewidth]{2D_Elasticity_scalar_FEM_loss.png}
%         \caption*{(a) Scalar-valued FEM basis.}
%     \end{minipage}
%     \hfill
%     \begin{minipage}[t]{0.4\textwidth}
%         \centering
%         \includegraphics[width=\linewidth]{2D_Elasticity_FEM_loss.png}
%         \caption*{(b) Vector-valued FEM basis.}
%     \end{minipage}
%     \caption{
%     \textbf{Elastic plate with FEM bases.}
%     Training curves for scalar-valued and vector-valued FEM output bases.
%     In this example, the scalar-valued FEM basis is substantially more accurate than the vector-valued FEM construction.
%     }
%     \label{fig:elastic_curve_FEM}
% \end{figure}
\paragraph{L-shaped Darcy flow: two input fields and one output field}

We also consider Darcy flow in an L-shaped domain,
\begin{equation}
\begin{aligned}
\nabla\cdot\big(k(\mathbf x)\nabla u(\mathbf x)\big)+f(\mathbf x)&=0,
\qquad
\mathbf x\in\Omega=(0,1)^2\setminus[0.5,1)^2,\\
u(\mathbf x)&=g(\mathbf x),
\qquad
\mathbf x\in\partial\Omega.
\end{aligned}
\label{eq:lshape-darcy}
\end{equation}
The goal is to learn the multi-input operator $G:[k(\mathbf x),f(\mathbf x)]\mapsto u(\mathbf x).$
The dataset is from~\cite{ingebrand2025basis}. 
It contains 51000 training samples and 7000 testing samples. 
The input and output functions are sampled on grids with 736 and 450 points, respectively.

We use vector-valued RFM bases to jointly encode the two input components:
\[
V_{\mathrm{in}}^{m_1}
=
\operatorname{span}\{\boldsymbol\varphi_i(\mathbf x)\}_{i=1}^{m_1},
\qquad
\boldsymbol\varphi_i(\mathbf x)
=
[\varphi_i^f(\mathbf x),\varphi_i^k(\mathbf x)].
\]
Thus,
\begin{equation}
[f(\mathbf x),k(\mathbf x)]
\approx
\sum_{j=1}^{m_1}a_j\boldsymbol\varphi_j(\mathbf x).
\label{eq:lshape-vector-input-basis}
\end{equation}
This joint coefficient representation reduces the input dimension relative to concatenating pointwise observations from both components.

\begin{figure}[!ht]
    \centering
    \begin{minipage}[t]{0.3\textwidth}
        \centering
        \includegraphics[width=\linewidth]{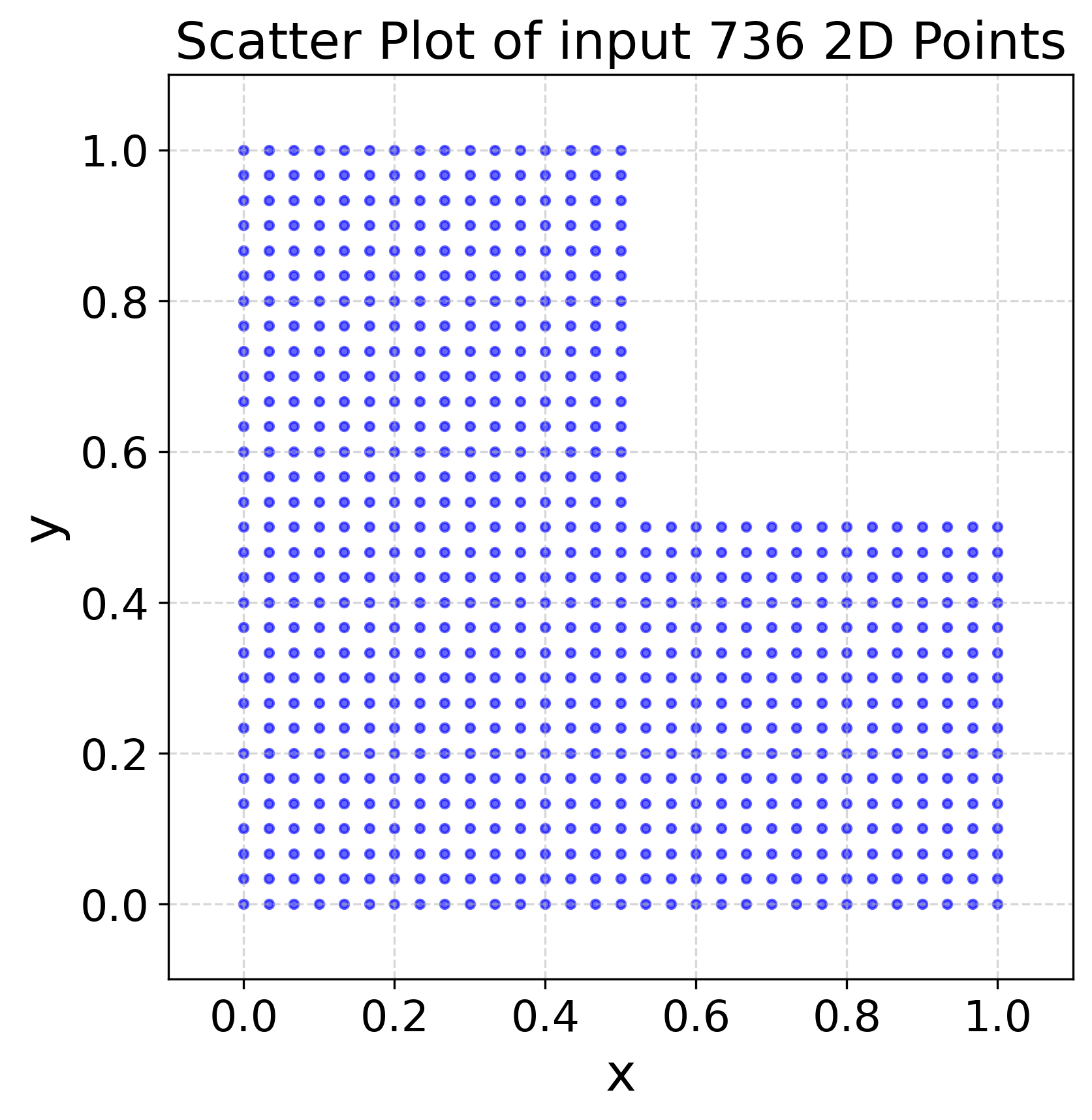}
        \caption*{(a) Input sample points.}
    \end{minipage}
    \hspace{0.05\textwidth}
    \begin{minipage}[t]{0.3\textwidth}
        \centering
        \includegraphics[width=\linewidth]{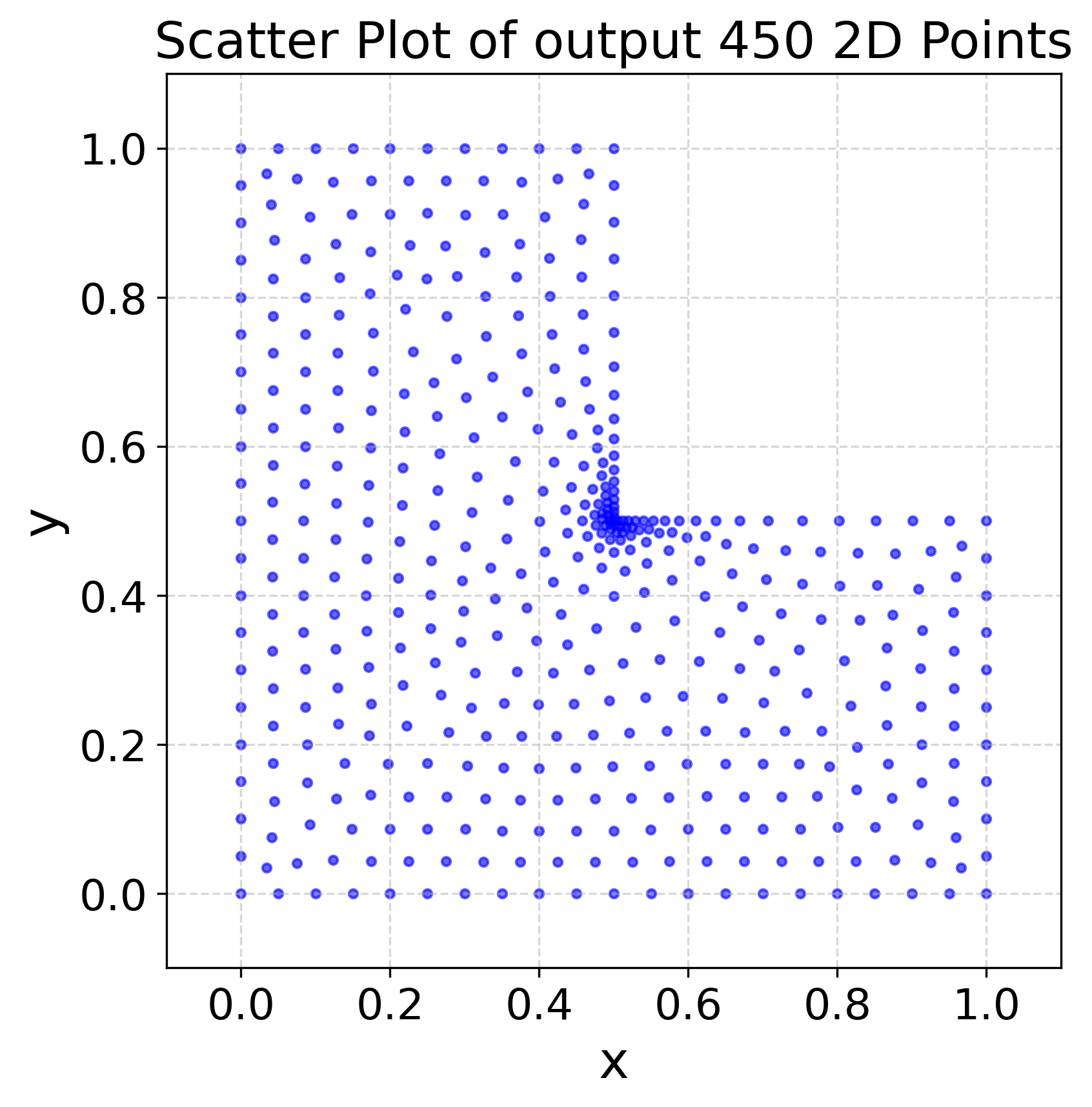}
        \caption*{(b) Output sample points.}
    \end{minipage}
    \caption{
    \textbf{L-shaped Darcy flow.}
    Sample points for the two input fields and the output field.
    }
    \label{fig:DarcyLshapegridpoints}
\end{figure}

\begin{figure}[!ht]
    \centering
    \begin{minipage}[t]{0.33\textwidth}
        \centering
        \includegraphics[width=\linewidth]{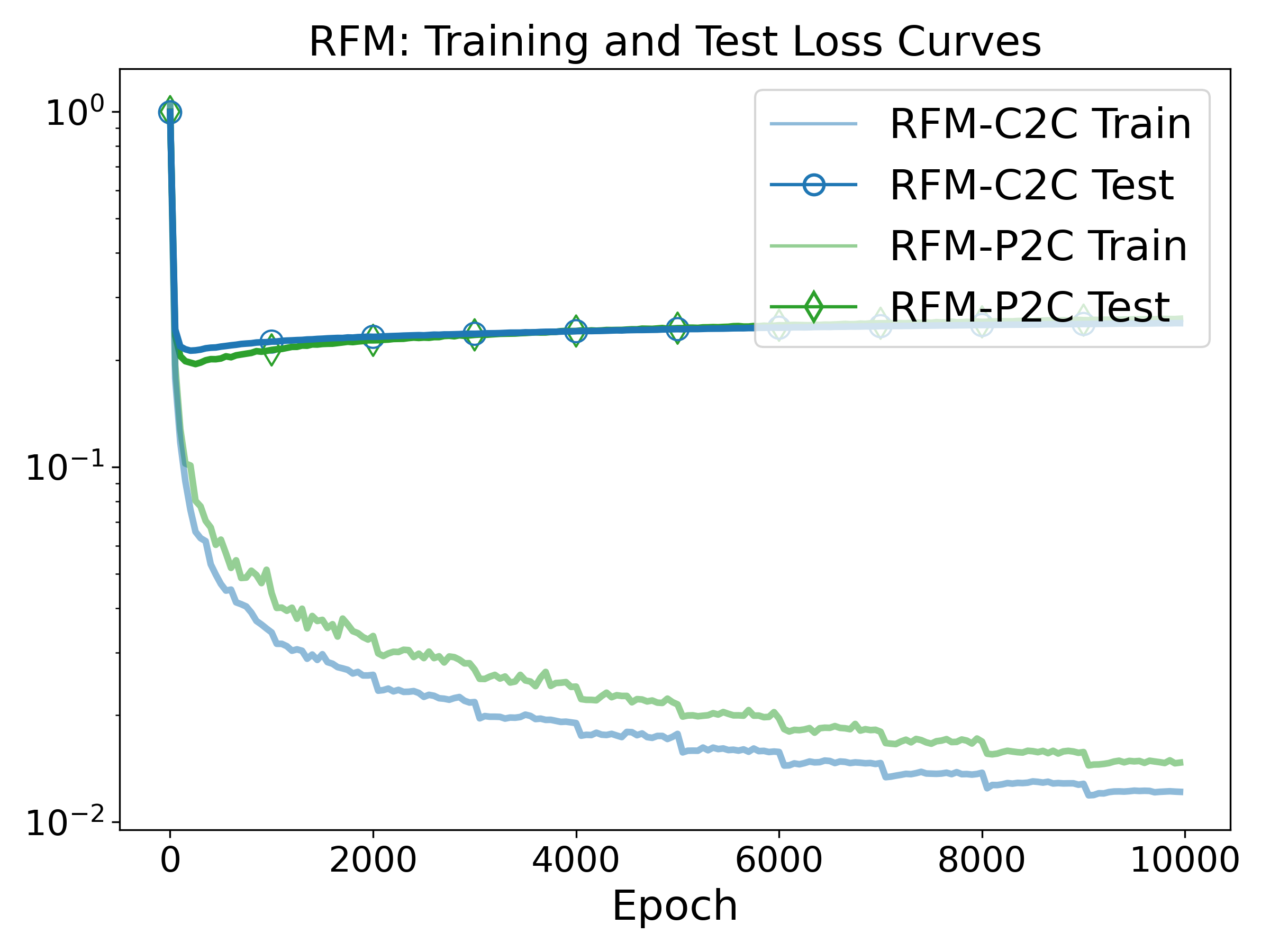}
        \caption*{(a) 1000 training samples.}
    \end{minipage}
    \hfill
    \begin{minipage}[t]{0.33\textwidth}
        \centering
        \includegraphics[width=\linewidth]{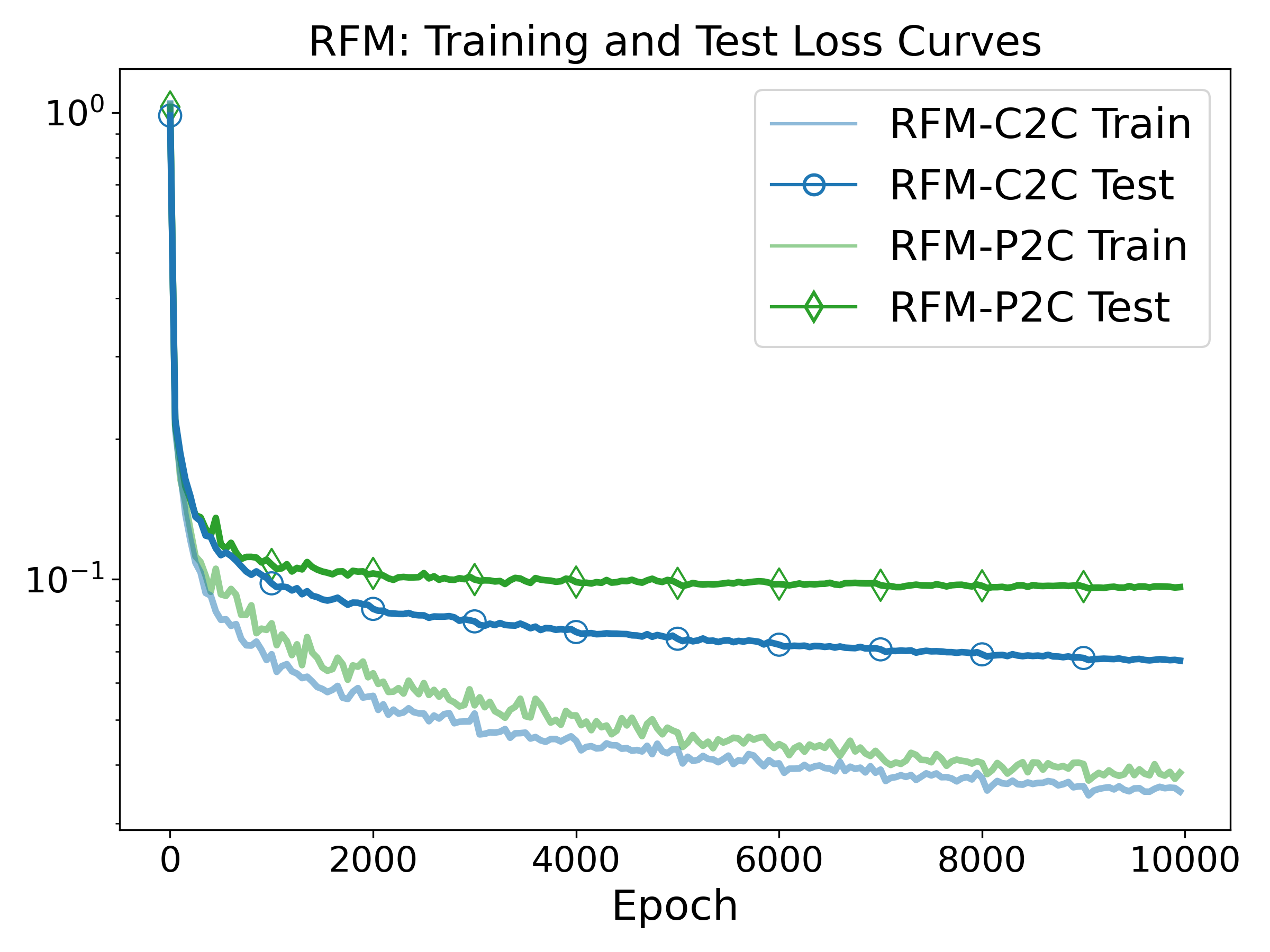}
        \caption*{(b) 5100 training samples.}
    \end{minipage}
    \hfill
    \begin{minipage}[t]{0.33\textwidth}
        \centering
        \includegraphics[width=\linewidth]{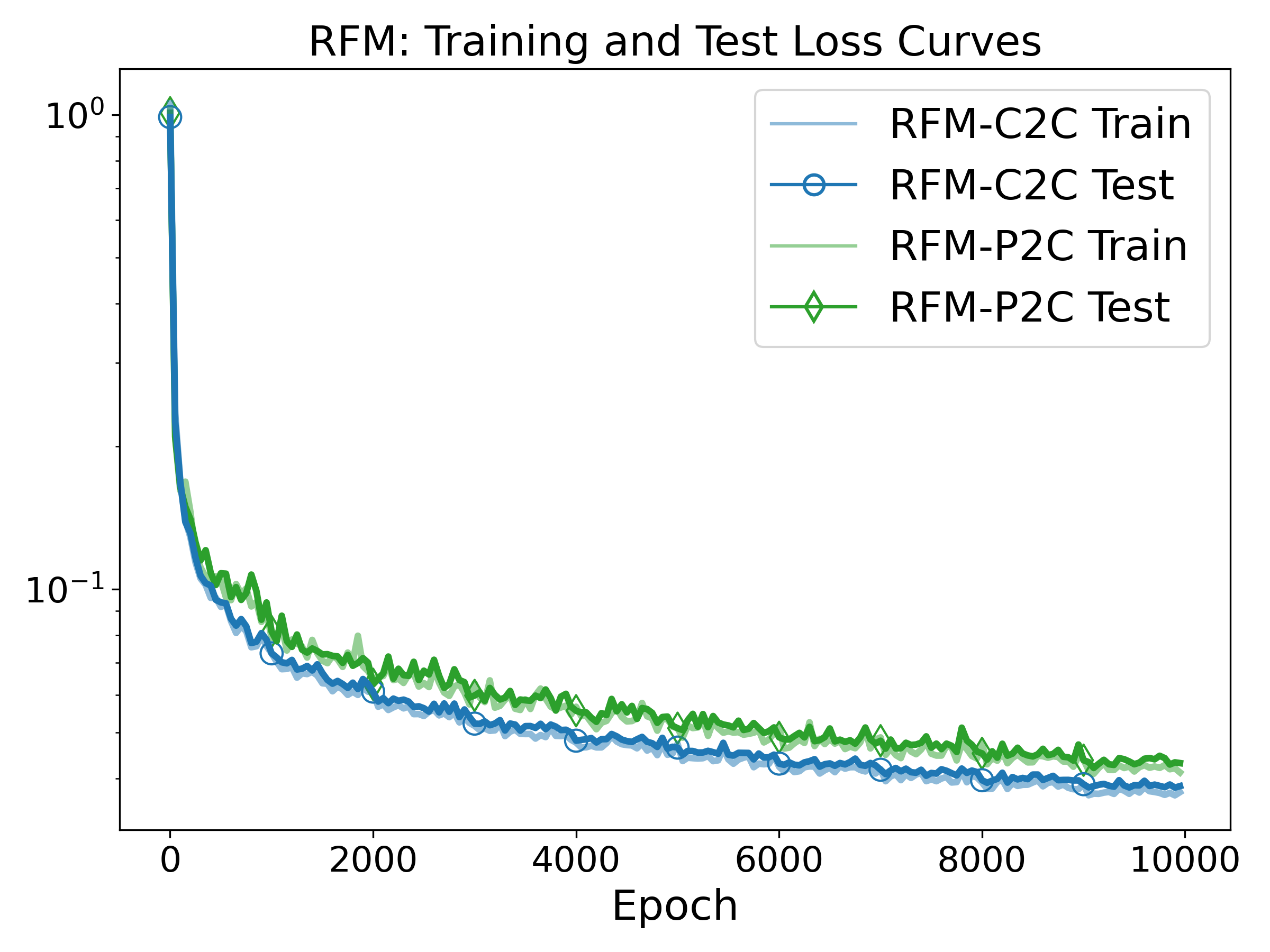}
        \caption*{(c) 51000 training samples.}
    \end{minipage}
    \caption{
    \textbf{L-shaped Darcy flow.}
    Training curve comparison between RFM-C2C and RFM-P2C for different training-set sizes.
    }
    \label{fig:Darcy_Lshape_training_curve}
\end{figure}

\begin{figure}[!ht]
    \centering
    \begin{minipage}[t]{0.33\textwidth}
        \centering
        \includegraphics[width=\linewidth]{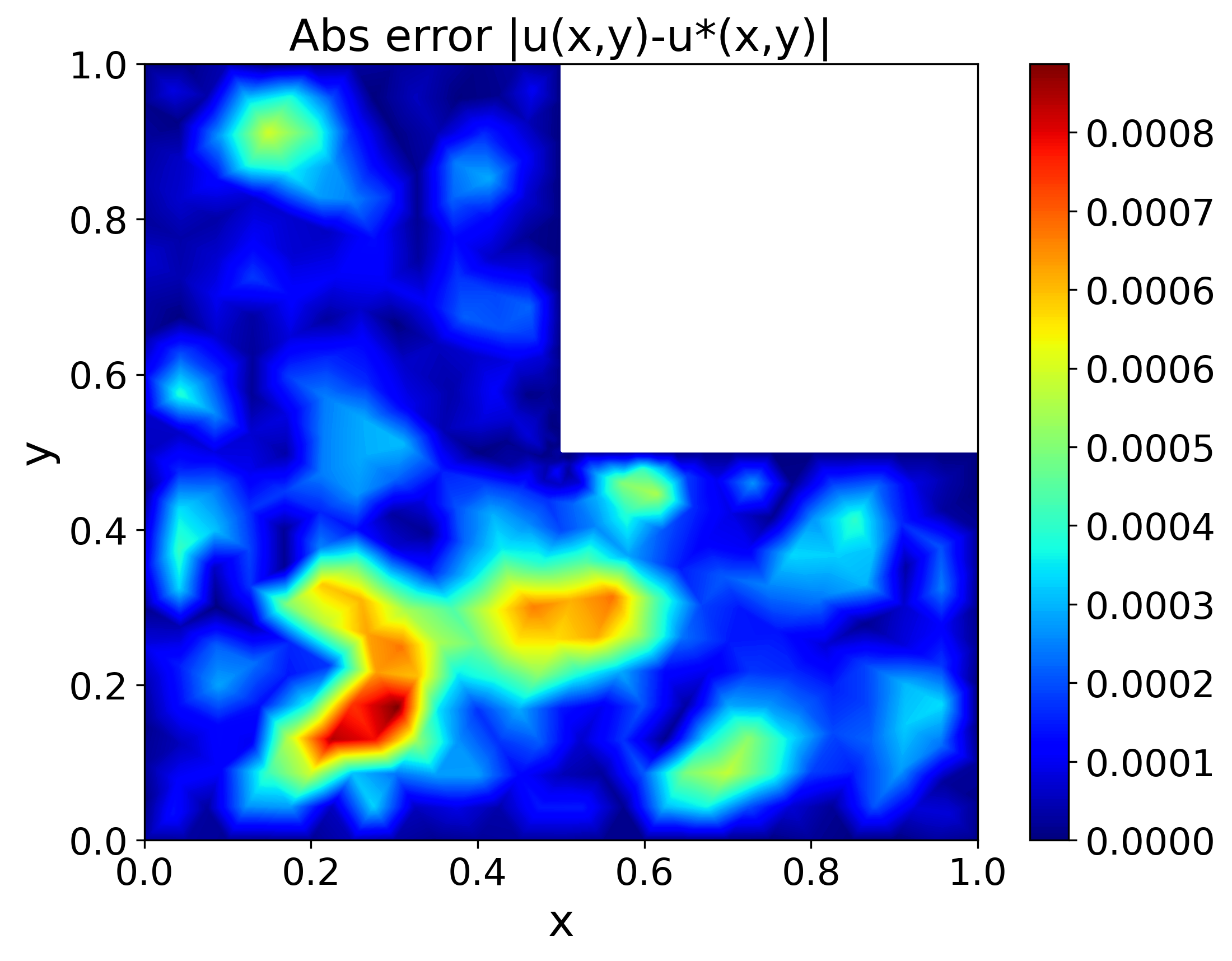}
        \caption*{(a) Pointwise error.}
    \end{minipage}
    \hfill
    \begin{minipage}[t]{0.33\textwidth}
        \centering
        \includegraphics[width=\linewidth]{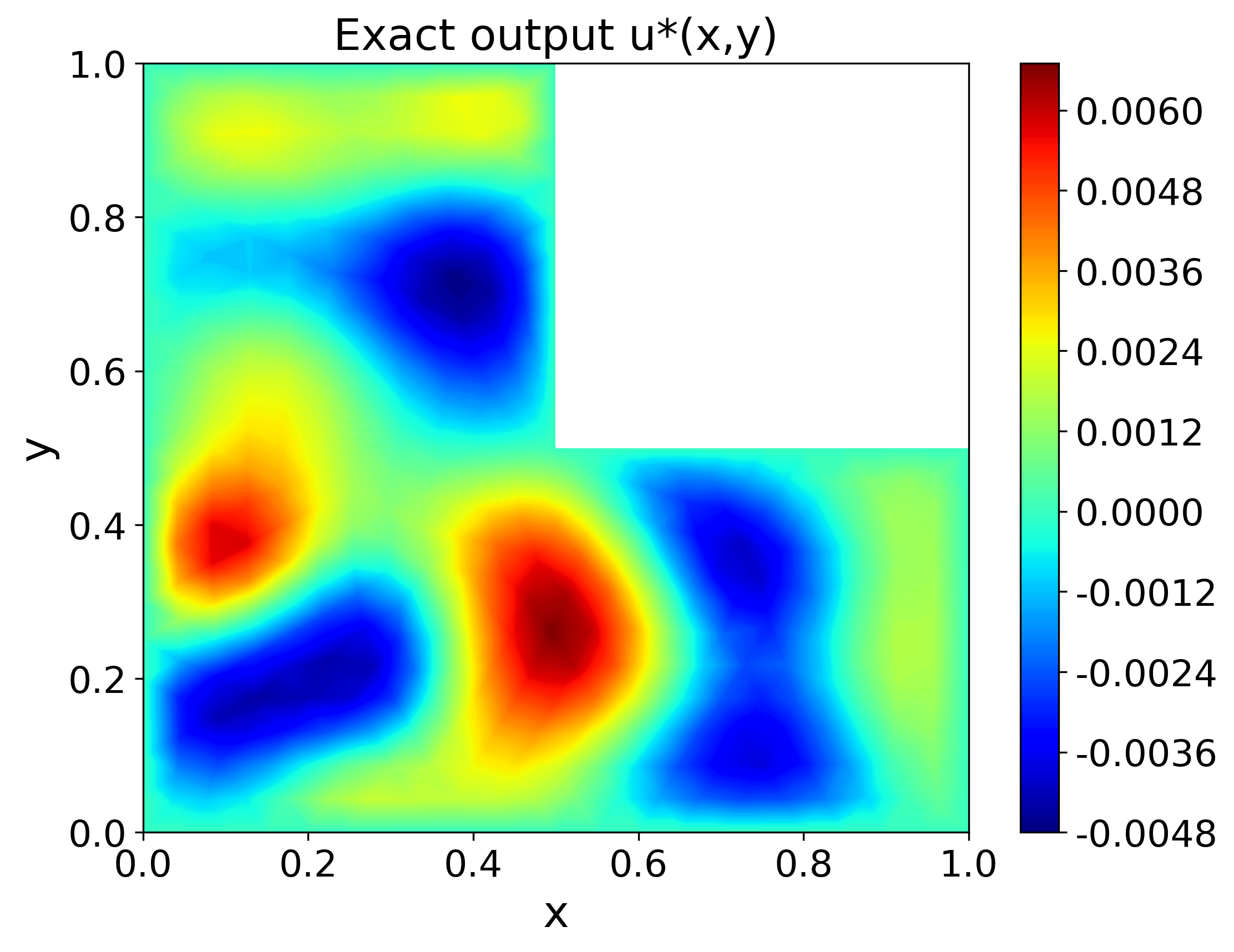}
        \caption*{(b) Exact solution.}
    \end{minipage}
    \hfill
    \begin{minipage}[t]{0.33\textwidth}
        \centering
        \includegraphics[width=\linewidth]{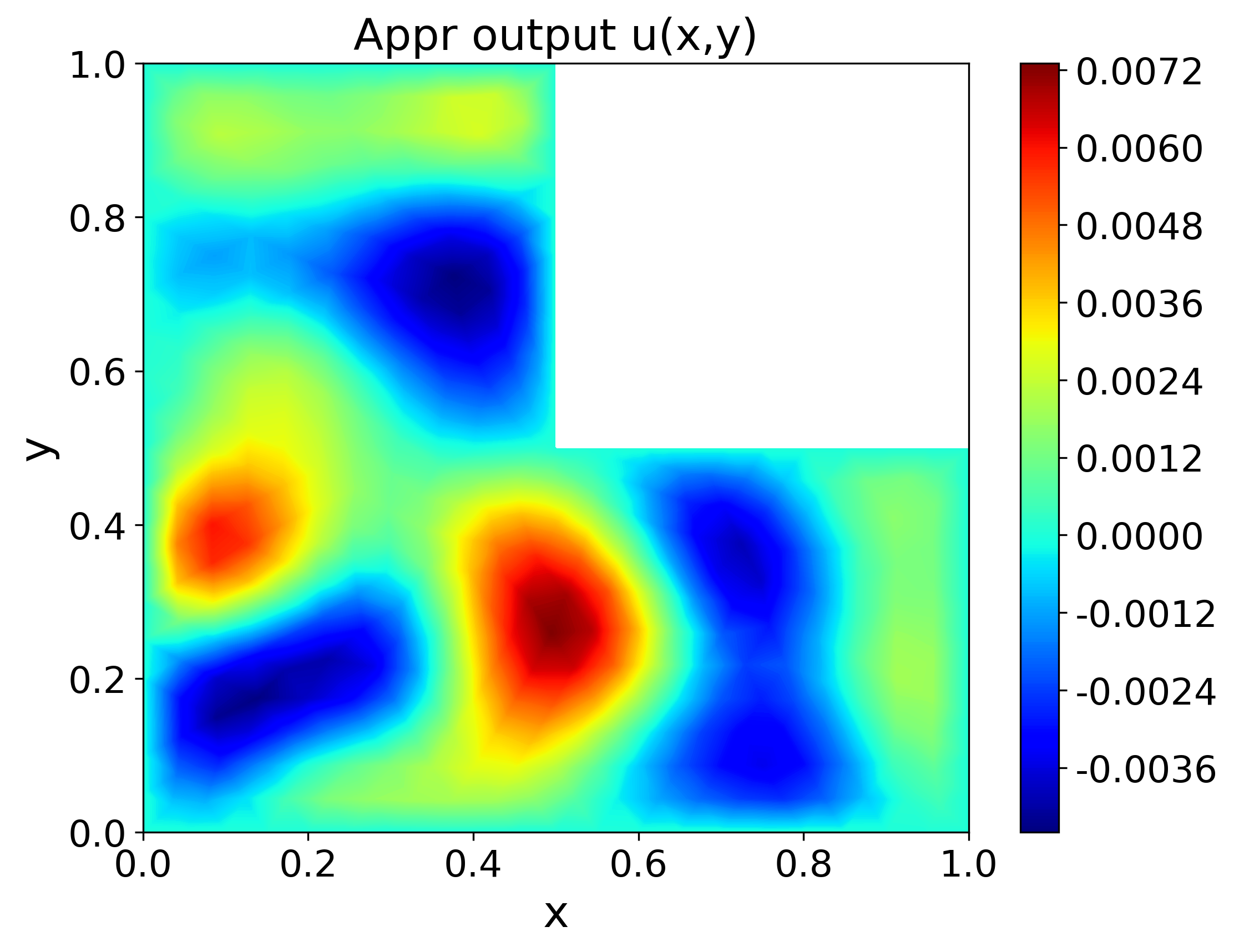}
        \caption*{(c) Predicted solution.}
    \end{minipage}
    \caption{
    \textbf{L-shaped Darcy flow.}
    Worst-case trained on 51000 training samples: (a) pointwise error, (b) exact solution, and (c) predicted solution.
    }
    \label{fig:DarcyLshape_worst_case}
\end{figure}

Table~\ref{tab:2d_darcy_Lshape_time} compares RFM-C2C and RFM-P2C using training subsets of size 1000, 5100, and 51000 and the same test subset of size 7000. 
The C2C method consistently achieves lower test error and shorter training time, with a more visible advantage in the smaller-data regimes. 
Figures~\ref{fig:DarcyLshapegridpoints}--\ref{fig:DarcyLshape_worst_case} show the sample grids, training curves, and a worst-case prediction.

\begin{table}[!htpb]
    \centering
    \caption{
    L-shaped Darcy flow.
    Vector-valued RFM-C2C versus RFM-P2C with different training-set sizes.
    The C2C model uses \(m_1=512\) input coefficients and \(m_2=1024\) output coefficients, while the P2C model uses \(736\times 2\) pointwise input values.
    Within each training-set size, the lower C2C/P2C training time and test error are shown in bold.
    }
    \label{tab:2d_darcy_Lshape_time}
    \footnotesize
    \renewcommand{\arraystretch}{1.12}
    \setlength{\tabcolsep}{3.5pt}
    \begin{tabular}{@{}lcccccc@{}}
        \toprule
        & \multicolumn{2}{c}{Training time}
        & \multicolumn{2}{c}{Test RL2E}
        & \multicolumn{2}{c}{Network architecture}\\
        \cmidrule(lr){2-3}\cmidrule(lr){4-5}\cmidrule(l){6-7}
        Training samples & C2C & P2C & C2C & P2C & C2C & P2C \\ 
        \midrule
        1000 
        & \textbf{19.11 s}
        & 20.59 s
        & \textbf{2.54e-1}
        & 2.62e-1 
        & [512,256\(\times\)2,1024]
        & [1472,256\(\times\)2,1024]\\ 
        5100 
        & \textbf{20.96 s}
        & 24.01 s
        & \textbf{6.70e-2}
        & 9.62e-2 
        & [512,256\(\times\)2,1024]
        & [1472,256\(\times\)2,1024]\\ 
        51000 
        & \textbf{156.48 s}
        & 186.67 s
        & \textbf{3.86e-2} 
        & 4.21e-2 
        & [512,256\(\times\)2,1024]
        & [1472,256\(\times\)2,1024]\\ 
        \bottomrule
    \end{tabular}
\end{table}

\subsubsection{Nonlinear time-dependent and weak-solution problems}
\label{subsubsec:time-dependent-weak}

\paragraph{KdV--Burgers equation}

To test nonlinear time-dependent solution operators, we consider the KdV--Burgers equation with periodic boundary conditions:
\begin{equation}
\begin{aligned}
\frac{\partial u}{\partial t}+u\frac{\partial u}{\partial x}
&=
\nu\frac{\partial^2u}{\partial x^2}
-
\mathcal D\frac{\partial^3u}{\partial x^3},
&& x\in(0,1),\ t\in(0,1),\\
u(0,t)&=u(1,t),
&& t\in[0,1],\\
\frac{\partial u}{\partial x}(0,t)&=\frac{\partial u}{\partial x}(1,t),
&& t\in[0,1],\\
u(x,0)&=f(x),
&& x\in[0,1].
\end{aligned}
\label{eq:kdv-burgers}
\end{equation}

\begin{figure}[!ht]
    \centering
    \begin{minipage}[t]{0.33\textwidth}
        \centering
        \includegraphics[width=\linewidth]{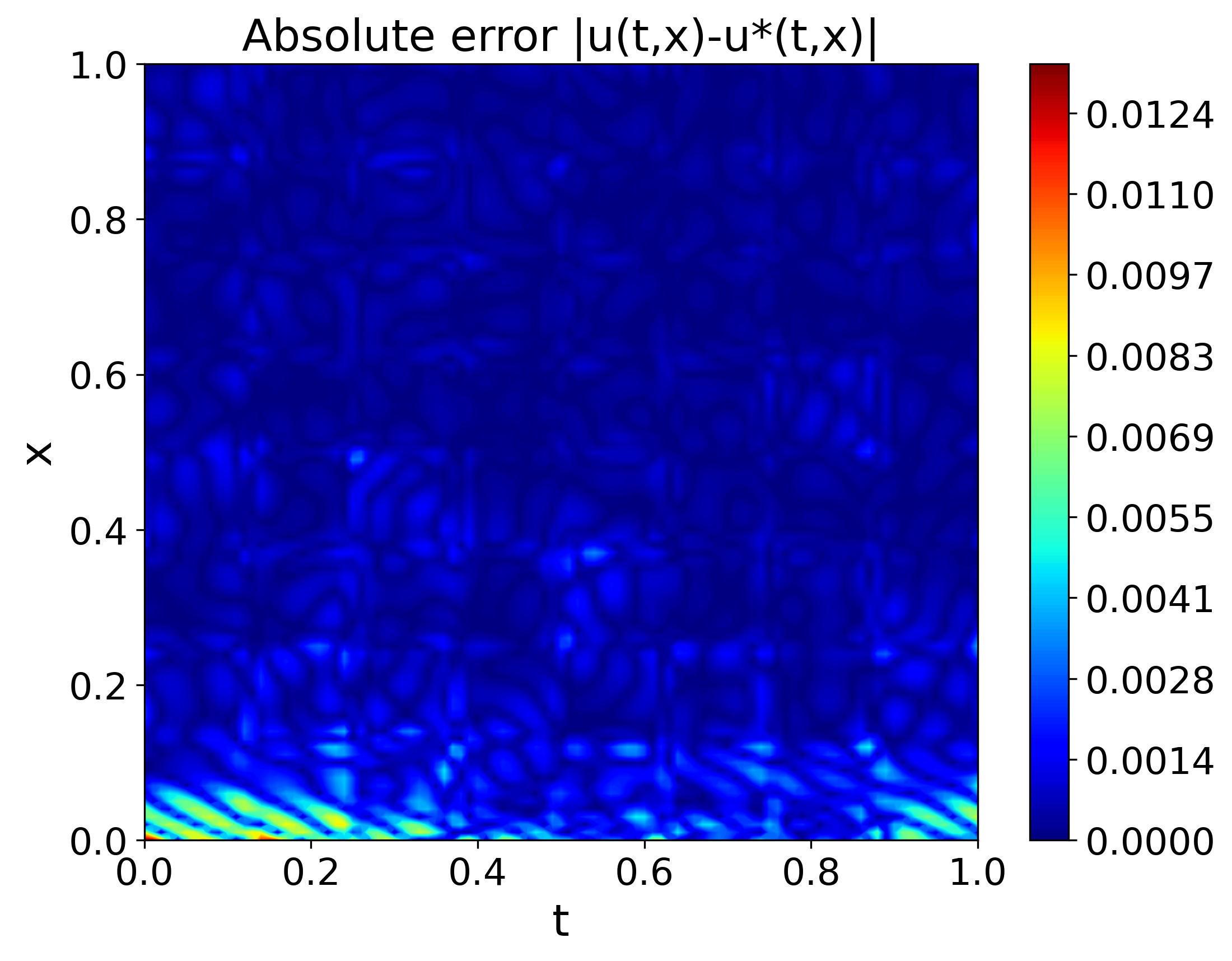}
        \caption*{(a) Pointwise error.}
    \end{minipage}
    \hfill
    \begin{minipage}[t]{0.33\textwidth}
        \centering
        \includegraphics[width=\linewidth]{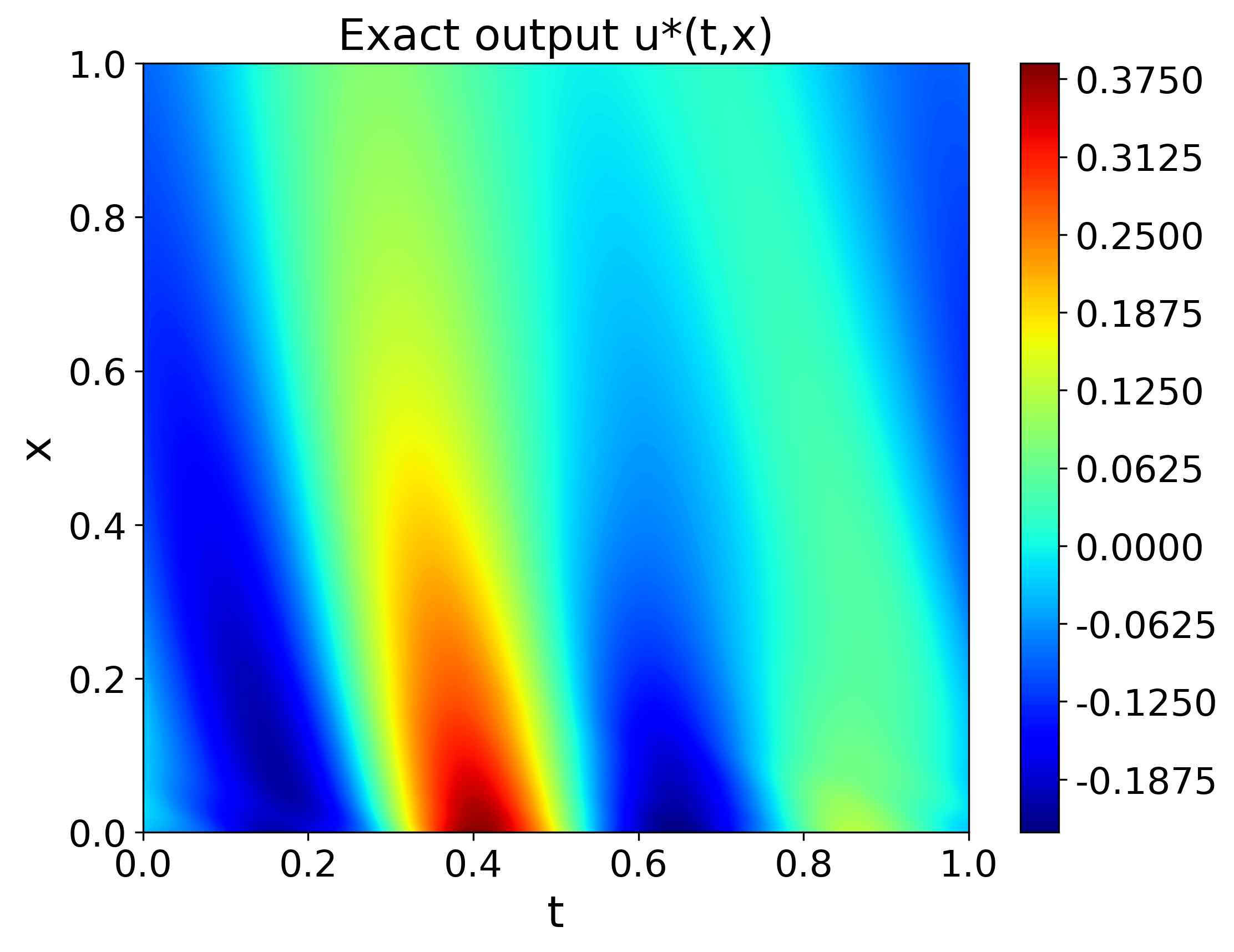}
        \caption*{(b) Exact solution.}
    \end{minipage}
    \hfill
    \begin{minipage}[t]{0.33\textwidth}
        \centering
        \includegraphics[width=\linewidth]{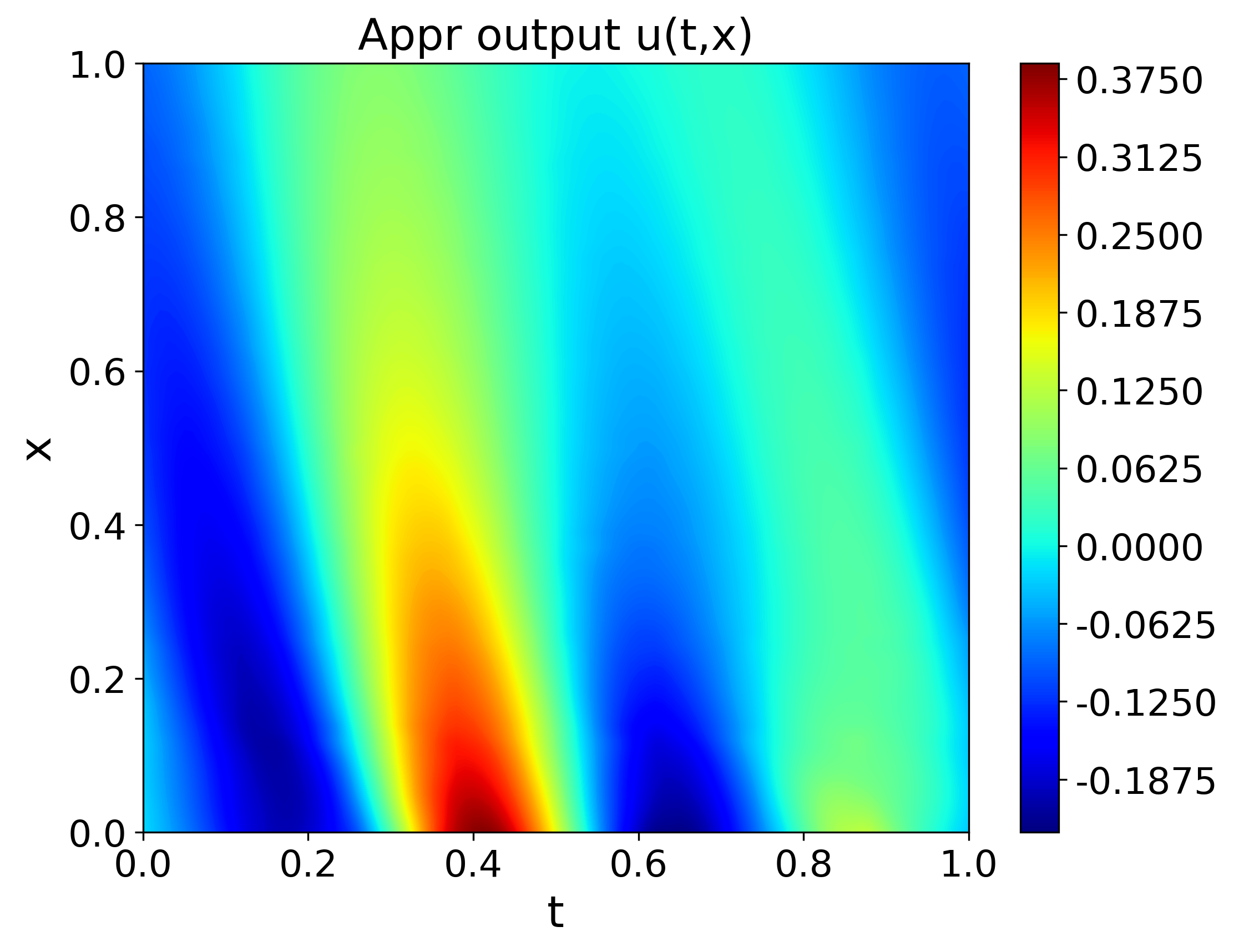}
        \caption*{(c) Predicted solution.}
    \end{minipage}
    \caption{
    \textbf{KdV--Burgers equation.}
    Worst-case test sample solved by RFM-C2C: (a) pointwise error, (b) exact solution, and (c) predicted solution.
    }
    \label{fig:KdV_Burgworst-case}
\end{figure}

Here \(\nu=0.01\) and \(\mathcal D=0.001\). 
The initial condition \(f(x)\) is drawn from \(\mathcal N(0,25^2(-\Delta+25I)^{-4})\). 
The reference solutions are generated using a pseudo-spectral method with a fine temporal discretization, following~\cite{hua2023basis}. 
The input resolution is 101, and the output is a space-time field on a \(101\times101\) grid. 
The target operator is $
G:f(x)\mapsto u(x,t).
$ We use non-separable spatiotemporal RFM bases for the output. 
The RFM-C2C model achieves a test RL2E of \(3.29e-3\). 
The worst-case prediction is shown in Figure~\ref{fig:KdV_Burgworst-case}.

\paragraph{Helmholtz equation with weak solutions}

\begin{figure}[!ht]
    \centering
    \begin{minipage}[t]{0.33\textwidth}
        \centering
        \includegraphics[width=\linewidth]{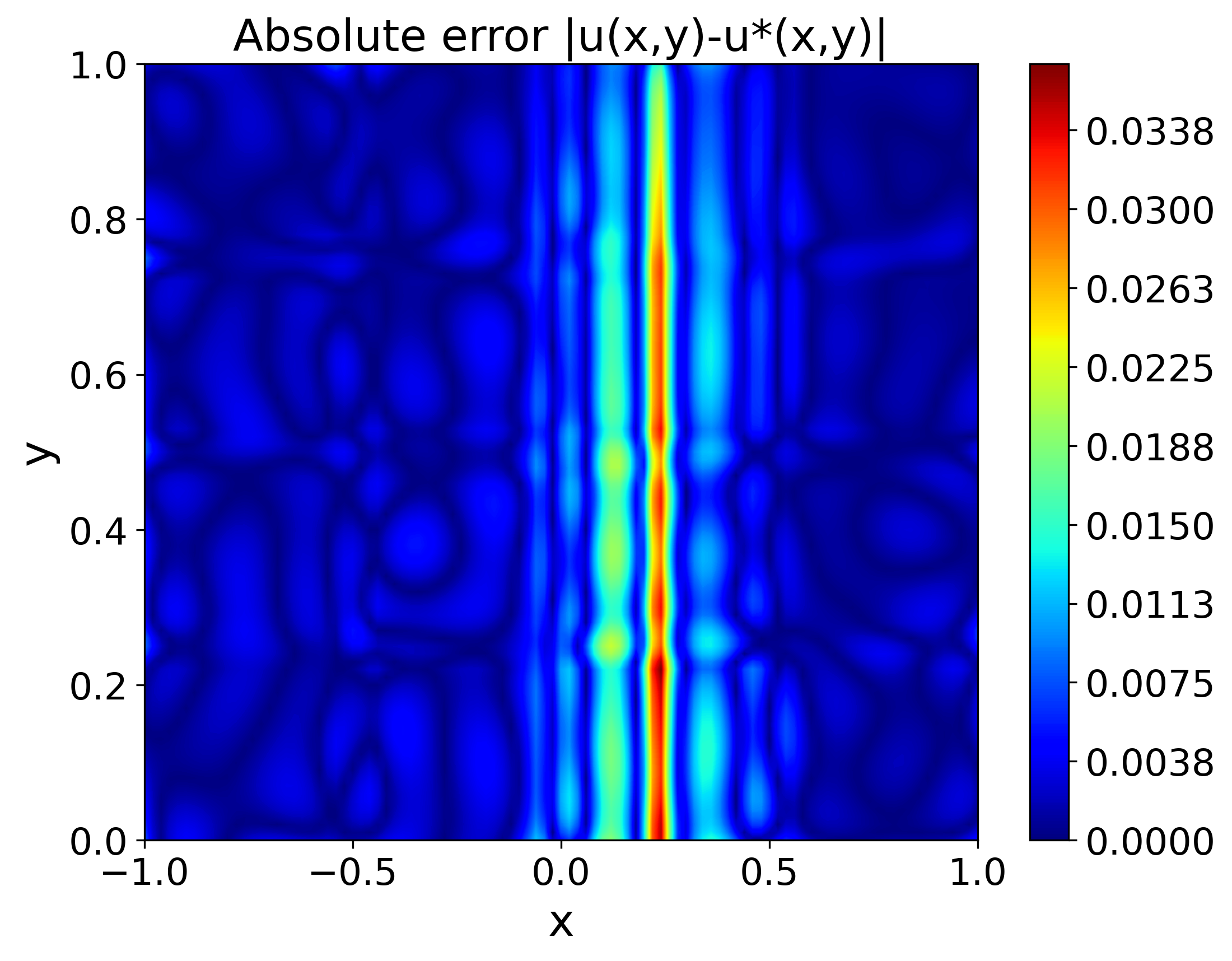}
        \caption*{(a) Pointwise error (RFM-C2C).}
    \end{minipage}
    \hfill
    \begin{minipage}[t]{0.33\textwidth}
        \centering
        \includegraphics[width=\linewidth]{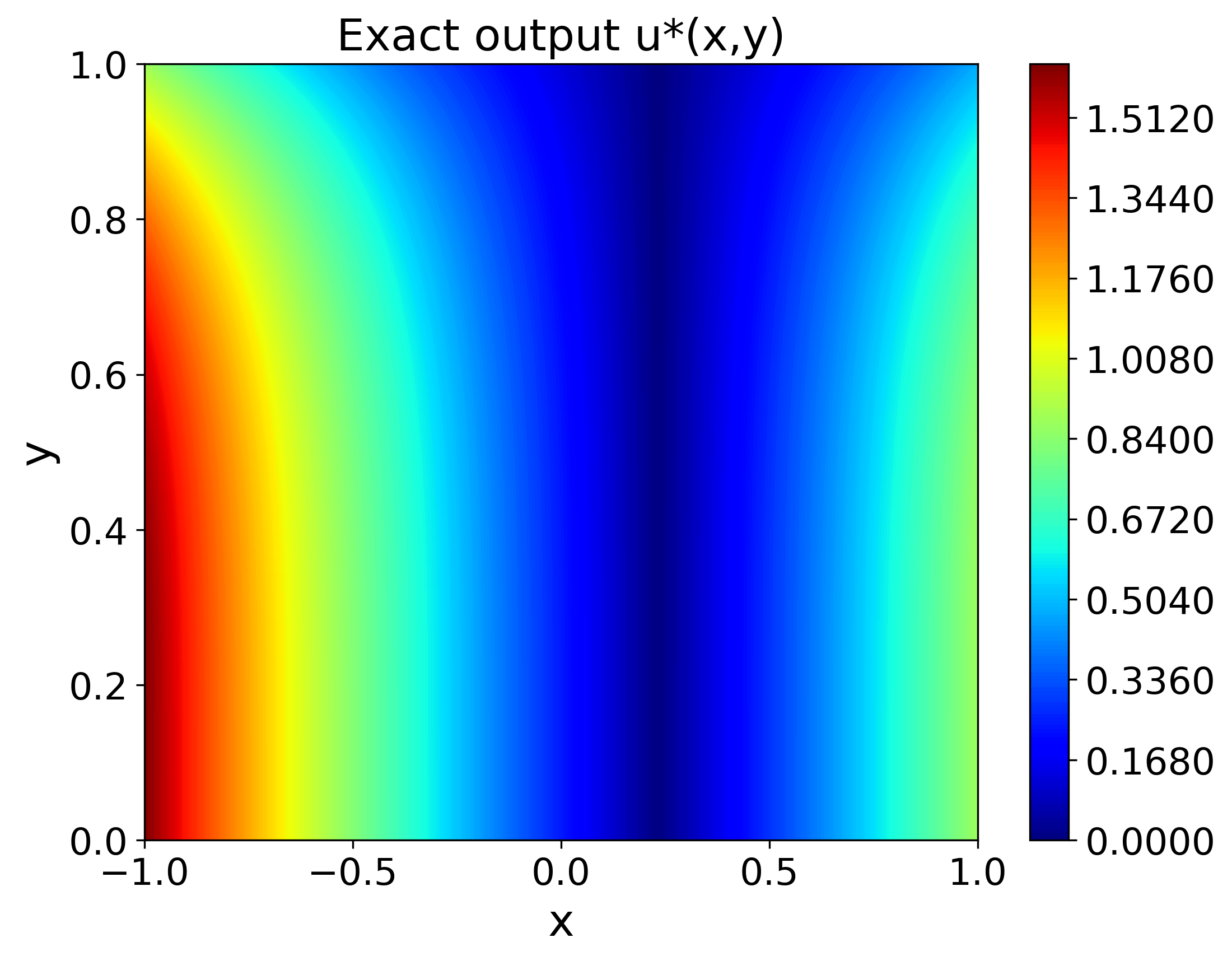}
        \caption*{(b) Exact solution (RFM-C2C).}
    \end{minipage}
    \hfill
    \begin{minipage}[t]{0.33\textwidth}
        \centering
        \includegraphics[width=\linewidth]{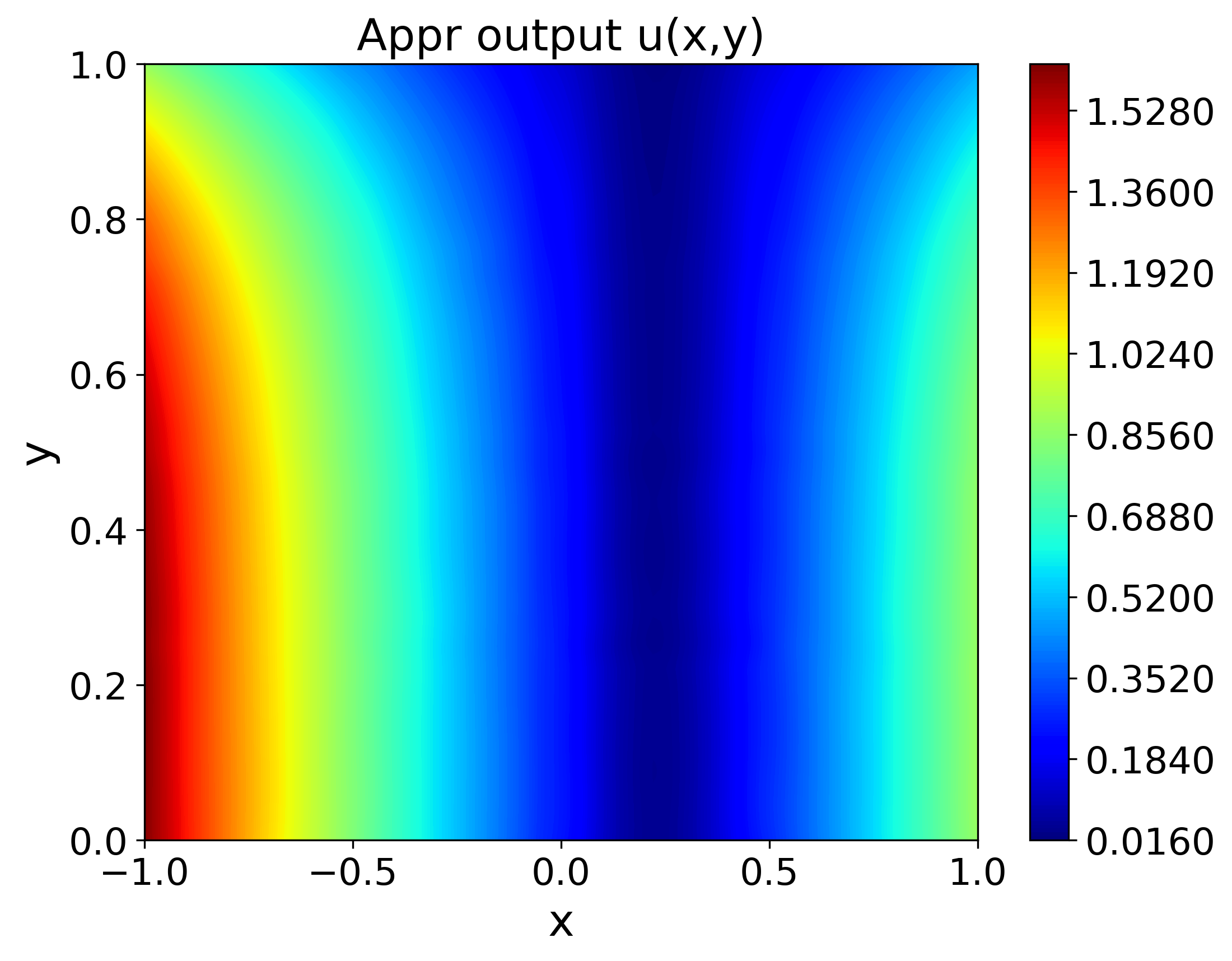}
        \caption*{(c) Predicted solution (RFM-C2C).}
    \end{minipage}
    \begin{minipage}[t]{0.33\textwidth}
        \centering
        \includegraphics[width=\linewidth]{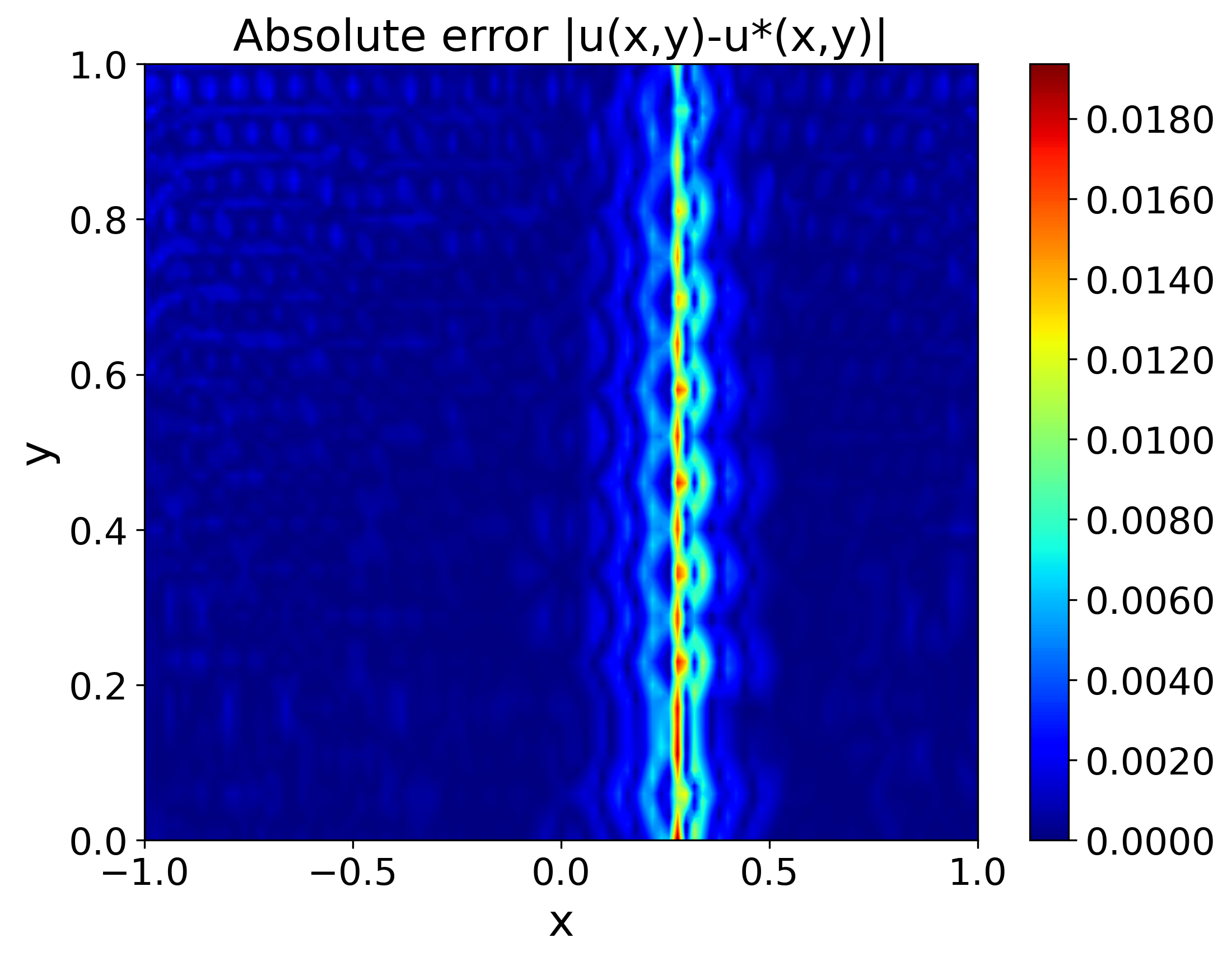}
        \caption*{(d) Pointwise error (FEM-C2C).}
    \end{minipage}
    \hfill
    \begin{minipage}[t]{0.33\textwidth}
        \centering
        \includegraphics[width=\linewidth]{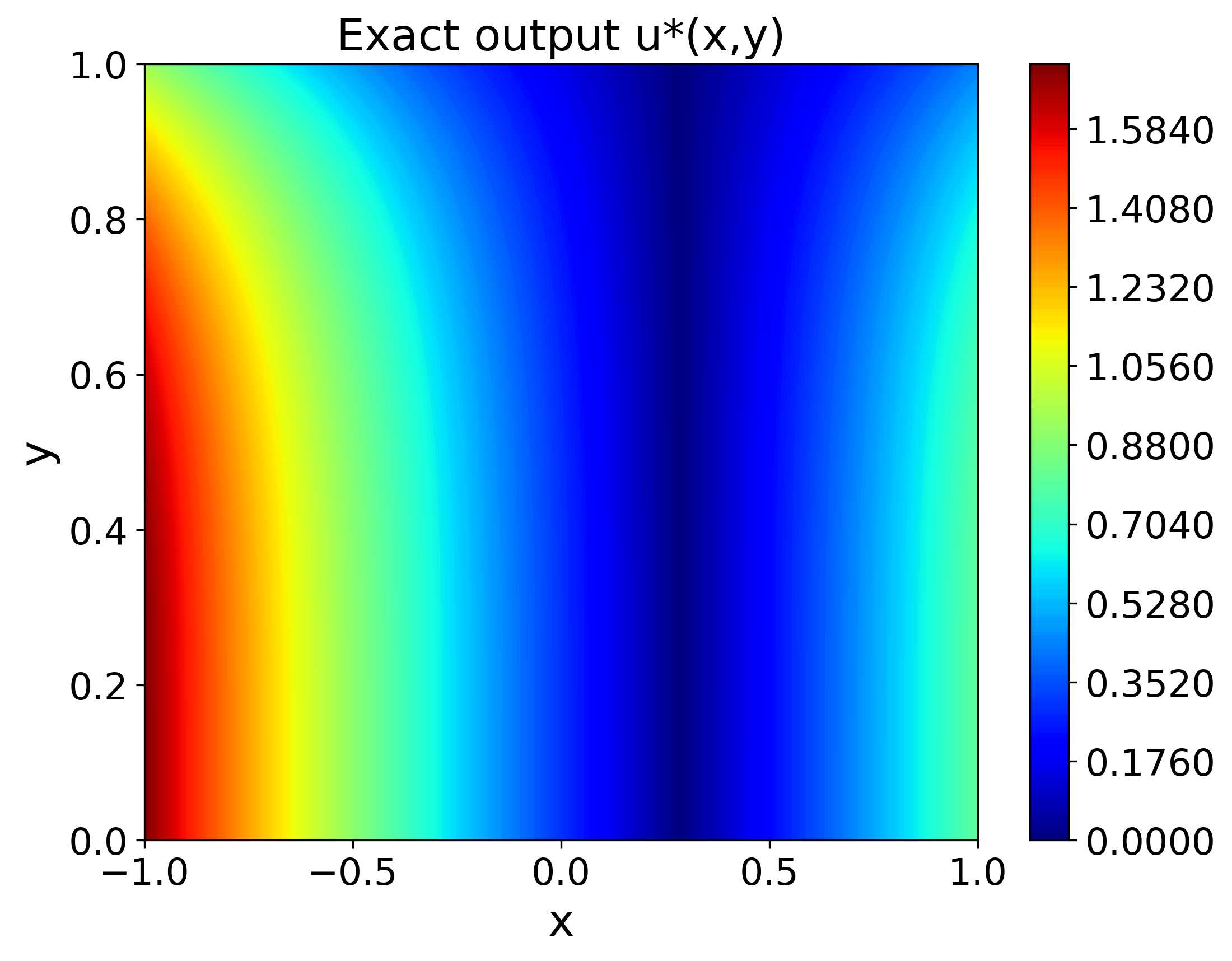}
        \caption*{(e) Exact solution (FEM-C2C).}
    \end{minipage}
    \hfill
    \begin{minipage}[t]{0.33\textwidth}
        \centering
        \includegraphics[width=\linewidth]{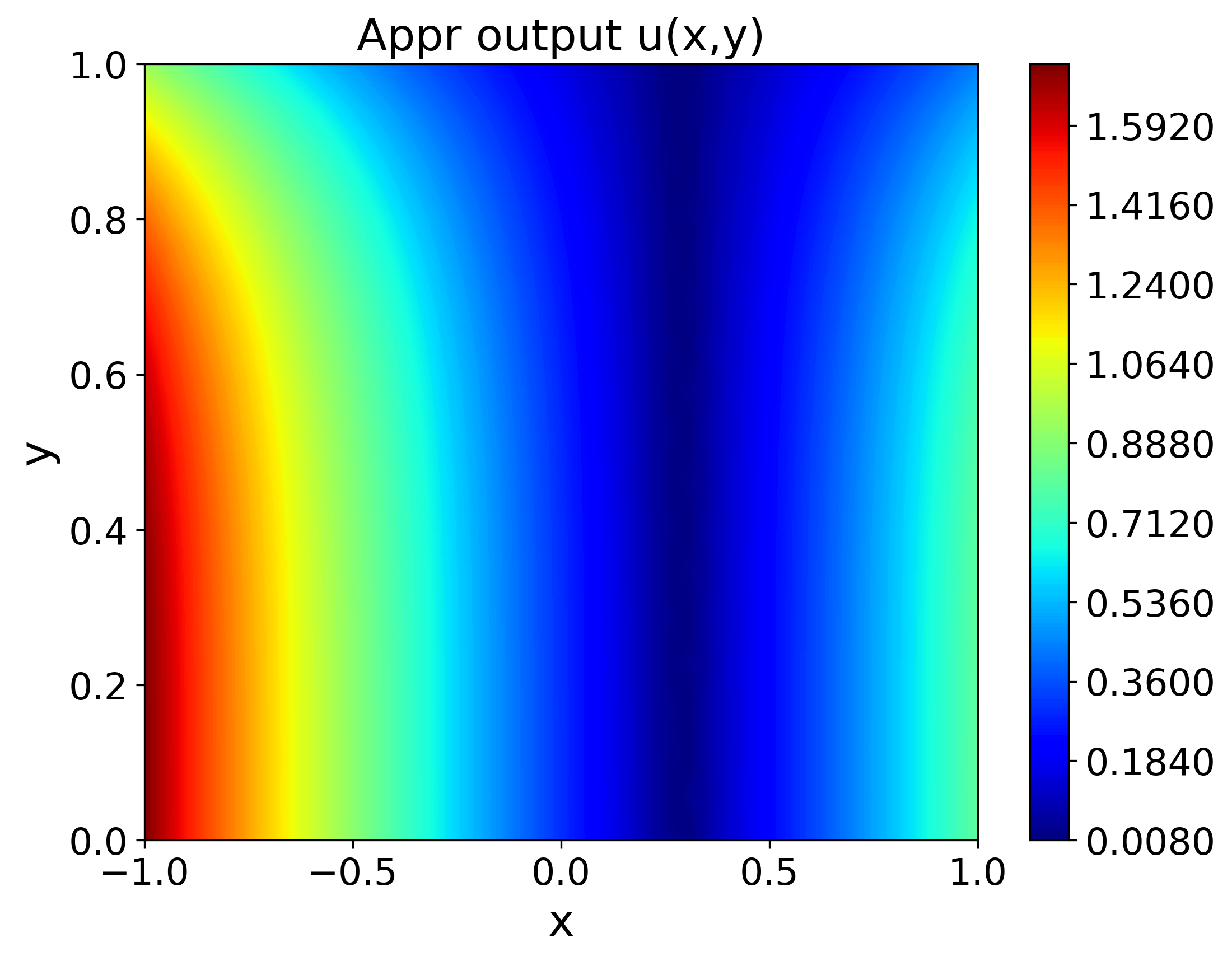}
        \caption*{(f) Predicted solution (FEM-C2C).}
    \end{minipage}
    % \begin{minipage}[t]{0.24\textwidth}
    %     \centering
    %     \includegraphics[width=\linewidth]{2D_Helmholtz_y0.png}
    %     \caption*{(g) Boundary \(y=0\) (RFM-C2C).}
    % \end{minipage}
    % \hfill
    % \begin{minipage}[t]{0.24\textwidth}
    %     \centering
    %     \includegraphics[width=\linewidth]{2D_Helmholtz_y1.png}
    %     \caption*{(h) Boundary \(y=1\) (RFM-C2C).}
    % \end{minipage}
    % \hfill
    % \begin{minipage}[t]{0.24\textwidth}
    %     \centering
    %     \includegraphics[width=\linewidth]{2D_Helmholtz_FEM_y0.png}
    %     \caption*{(i) Boundary \(y=0\) (FEM-C2C).}
    % \end{minipage}
    % \hfill
    % \begin{minipage}[t]{0.24\textwidth}
    %     \centering
    %     \includegraphics[width=\linewidth]{2D_Helmholtz_FEM_y1.png}
    %     \caption*{(j) Boundary \(y=1\) (FEM-C2C).}
    % \end{minipage}
    \caption{
    \textbf{2D Helmholtz equation.}
    (a)--(c) show the worst-case RFM-C2C prediction, and (d)--(f) show the worst-case FEM-C2C prediction.
    }
    \label{fig:Helmholtz_worst_case}
\end{figure}

\begin{figure}[!ht]
    \centering
    % \begin{minipage}[t]{0.3\textwidth}
    %     \centering
    %     \includegraphics[width=\linewidth]{2D_Helmholtz_error.png}
    %     \caption*{(a) Pointwise error (RFM-C2C).}
    % \end{minipage}
    % \hfill
    % \begin{minipage}[t]{0.3\textwidth}
    %     \centering
    %     \includegraphics[width=\linewidth]{2D_Helmholtz_exact.png}
    %     \caption*{(b) Exact solution (RFM-C2C).}
    % \end{minipage}
    % \hfill
    % \begin{minipage}[t]{0.3\textwidth}
    %     \centering
    %     \includegraphics[width=\linewidth]{2D_Helmholtz_pred.png}
    %     \caption*{(c) Predicted solution (RFM-C2C).}
    % \end{minipage}
    % \begin{minipage}[t]{0.3\textwidth}
    %     \centering
    %     \includegraphics[width=\linewidth]{2D_Helmholtz_FEM_error.png}
    %     \caption*{(d) Pointwise error (FEM-C2C).}
    % \end{minipage}
    % \hfill
    % \begin{minipage}[t]{0.3\textwidth}
    %     \centering
    %     \includegraphics[width=\linewidth]{2D_Helmholtz_FEM_exact.png}
    %     \caption*{(e) Exact solution (FEM-C2C).}
    % \end{minipage}
    % \hfill
    % \begin{minipage}[t]{0.3\textwidth}
    %     \centering
    %     \includegraphics[width=\linewidth]{2D_Helmholtz_FEM_pred.png}
    %     \caption*{(f) Predicted solution (FEM-C2C).}
    % \end{minipage}
    \begin{minipage}[t]{0.24\textwidth}
        \centering
        \includegraphics[width=\linewidth]{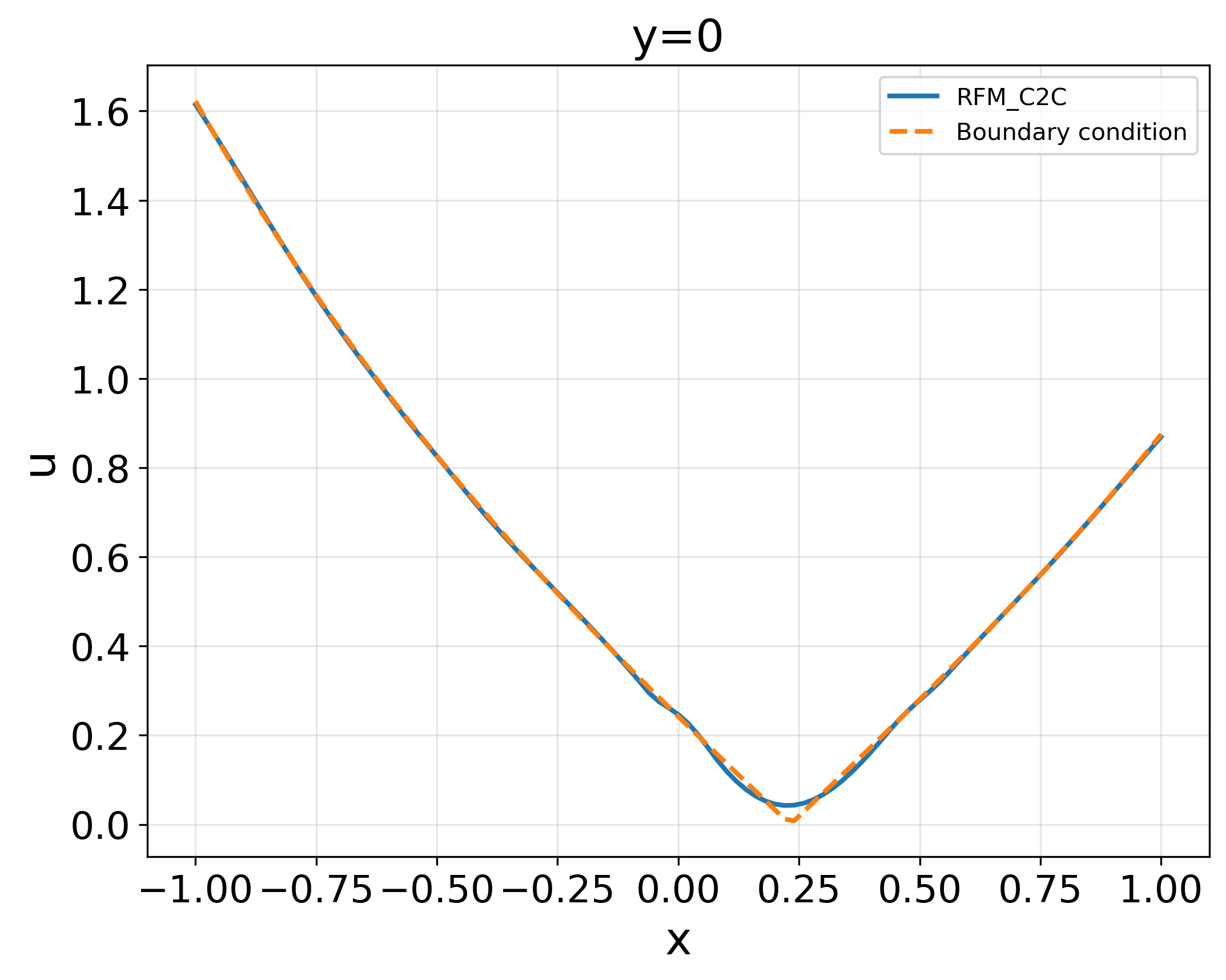}
        \caption*{(a) Boundary \(y=0\) (RFM-C2C).}
    \end{minipage}
    \hfill
    \begin{minipage}[t]{0.24\textwidth}
        \centering
        \includegraphics[width=\linewidth]{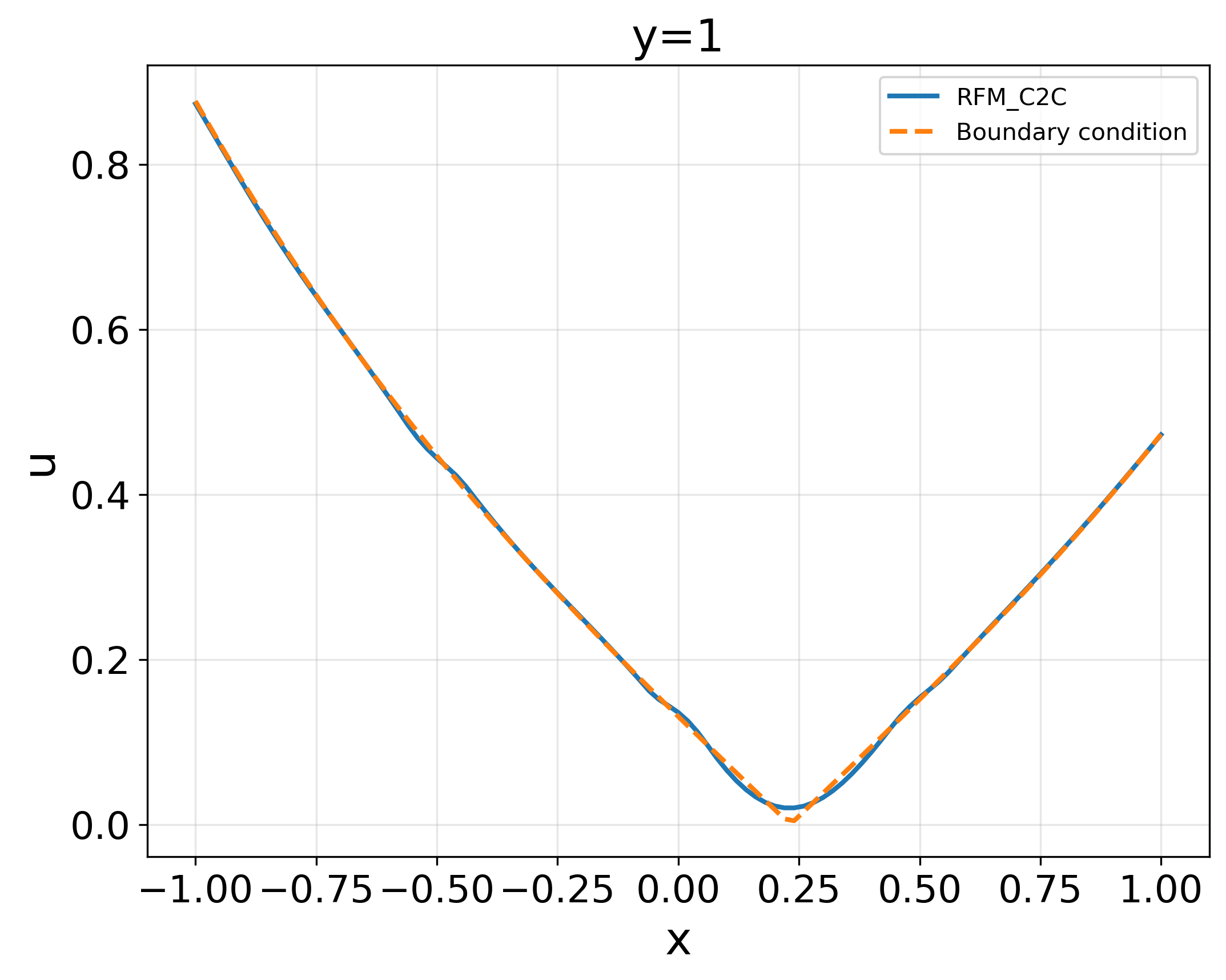}
        \caption*{(b) Boundary \(y=1\) (RFM-C2C).}
    \end{minipage}
    \hfill
    \begin{minipage}[t]{0.24\textwidth}
        \centering
        \includegraphics[width=\linewidth]{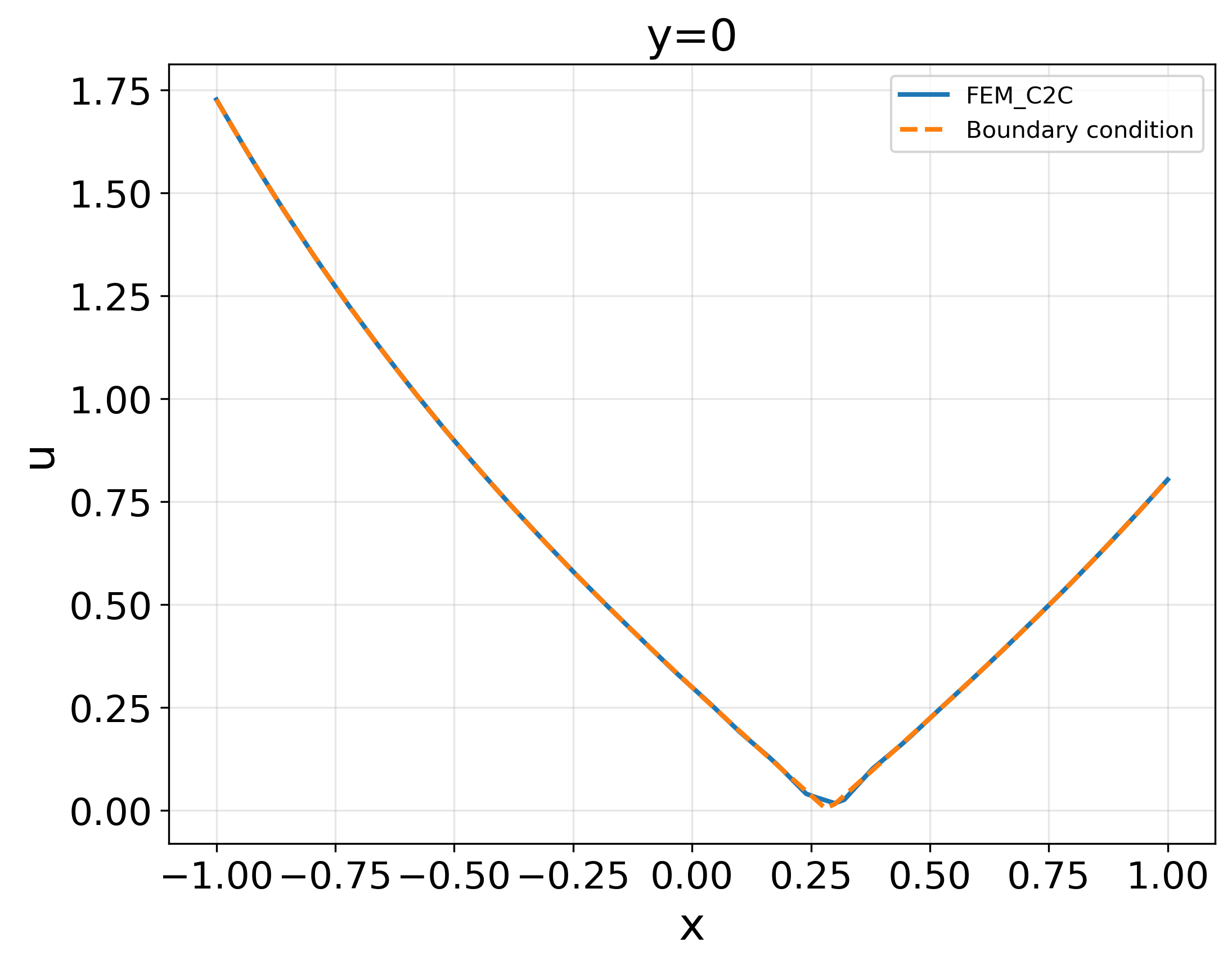}
        \caption*{(c) Boundary \(y=0\) (FEM-C2C).}
    \end{minipage}
    \hfill
    \begin{minipage}[t]{0.24\textwidth}
        \centering
        \includegraphics[width=\linewidth]{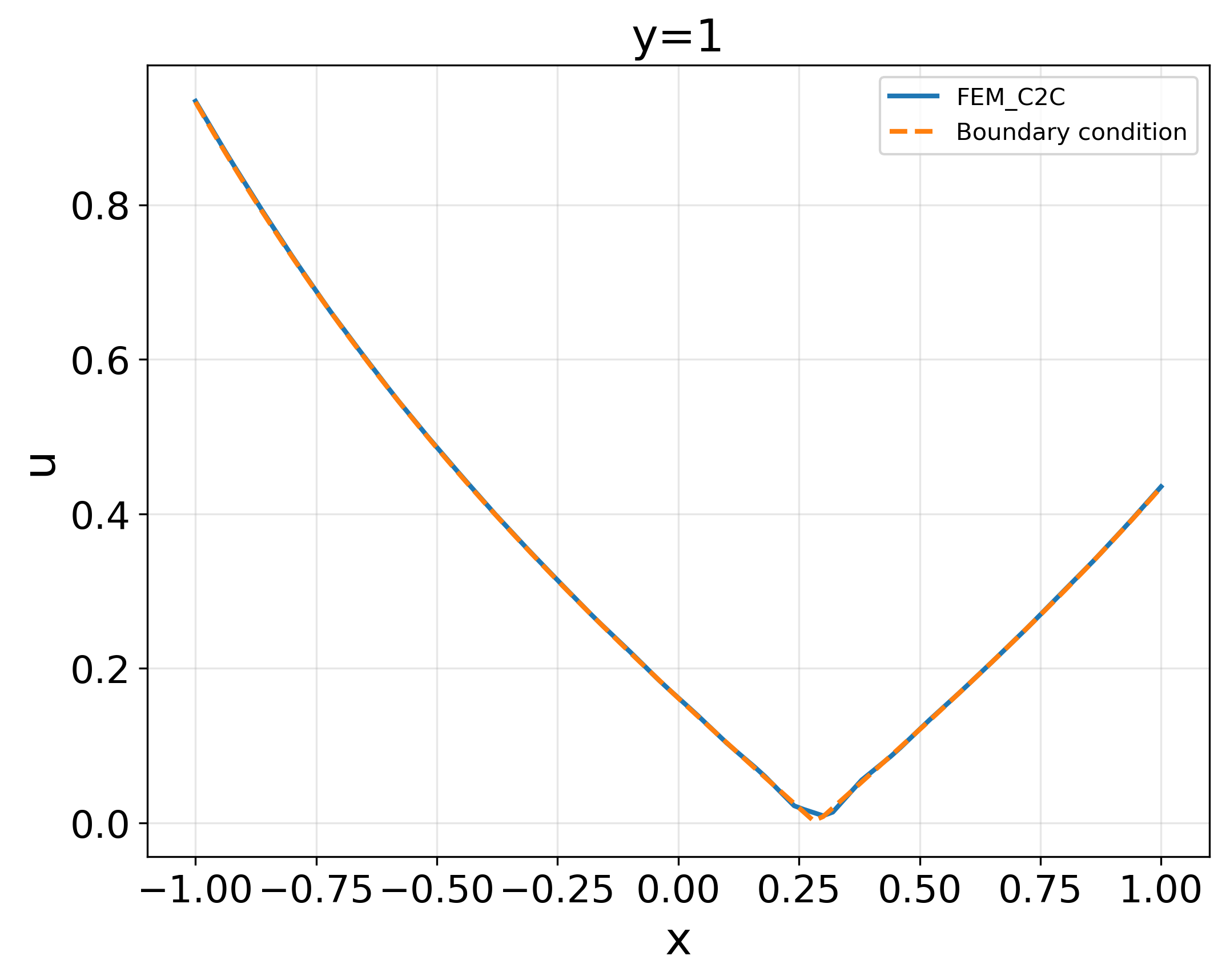}
        \caption*{(d) Boundary \(y=1\) (FEM-C2C).}
    \end{minipage}
    \caption{
    \textbf{2D Helmholtz equation.}
    (a)--(b) show boundary comparisons for RFM-C2C, and (c)--(d) show boundary comparisons for FEM-C2C.
    }
    \label{fig:Helmholtz_boundary}
\end{figure}

We next consider a Helmholtz equation on \(\Omega=(-1,1)\times(0,1)\). For each \(\beta\), define the interface
\[
\Gamma_\beta=\{(\beta,y):0<y<1\}.
\]
In the sense of distributions, the problem is
\begin{equation}
\begin{cases}
\Delta u+\lambda u
=f_{\mathrm{reg}}+2\alpha\cos(y^2)\delta_{\Gamma_\beta}, & \text{in }\mathcal D'(\Omega),\\
u|_{\partial\Omega}=\alpha\sinh(|x-\beta|)\cos(y^2),
\end{cases}
\label{eq:helmholtz-weak}
\end{equation}
where \(\lambda=1\), \(\delta_{\Gamma_\beta}\) is the line Dirac measure supported on \(\Gamma_\beta\), and
\[
f_{\mathrm{reg}}(x,y)
=
\alpha\sinh(|x-\beta|)
\left[(1+\lambda-4y^2)\cos(y^2)-2\sin(y^2)\right].
\]
The corresponding weak solution, which is also available analytically, is
\[
u(x,y)=\alpha\sinh(|x-\beta|)\cos(y^2).
\]
Indeed, the jump in \(\partial_xu\) across \(\Gamma_\beta\) is \(2\alpha\cos(y^2)\), producing the singular interface term above. 
We sample 50 values of \(\alpha\sim\mathcal U([1,2])\) and 50 values of \(\beta\sim\mathcal U([-0.5,0.5])\), giving 2500 samples in total, with 2000 for training and 500 for testing. 
In the numerical experiment, the neural input is the sampled regular part \(f_{\mathrm{reg}}\). Within this manufactured two-parameter family, \(f_{\mathrm{reg}}\) determines \(\alpha\), \(\beta\), and hence the associated interface term. The target operator is therefore
$
G:f_{\mathrm{reg}}(x,y)\mapsto u(x,y).
$
Both input and output functions are evaluated on a \(101\times101\) grid. 
RFM-C2C uses \((m_1,m_2)=(256,512)\), and FEM-C2C uses \((m_1,m_2)=(273,558)\). 
The RL2E values are \(6.51e-3\) for RFM-C2C and \(2.13e-3\) for FEM-C2C. 
Figures~\ref{fig:Helmholtz_worst_case}-\ref{fig:Helmholtz_boundary} show worst-case predictions and boundary comparisons.

% \begin{figure}[!ht]
%     \centering
%     \begin{minipage}[t]{0.3\textwidth}
%         \centering
%         \includegraphics[width=\linewidth]{2D_Helmholtz_FEM_error.png}
%         \caption*{(a) Pointwise error.}
%     \end{minipage}
%     \hfill
%     \begin{minipage}[t]{0.3\textwidth}
%         \centering
%         \includegraphics[width=\linewidth]{2D_Helmholtz_FEM_exact.png}
%         \caption*{(b) Exact solution.}
%     \end{minipage}
%     \hfill
%     \begin{minipage}[t]{0.3\textwidth}
%         \centering
%         \includegraphics[width=\linewidth]{2D_Helmholtz_FEM_pred.png}
%         \caption*{(c) Predicted solution.}
%     \end{minipage}
%     \begin{minipage}[t]{0.3\textwidth}
%         \centering
%         \includegraphics[width=\linewidth]{2D_Helmholtz_FEM_y0.png}
%         \caption*{(d) Boundary \(y=0\).}
%     \end{minipage}
%     \hfill
%     \begin{minipage}[t]{0.3\textwidth}
%         \centering
%         \includegraphics[width=\linewidth]{2D_Helmholtz_FEM_y1.png}
%         \caption*{(e) Boundary \(y=1\).}
%     \end{minipage}
%     \caption{
%     \textbf{2D Helmholtz equation solved by FEM-C2C.}
%     Worst-case test sample: (a) pointwise error, (b) exact solution, (c) predicted solution, and (d)--(e) boundary comparisons.
%     }
%     \label{fig:Helmholtz_FEM}
% \end{figure}
\subsubsection{High-dimensional proof-of-concept}
\label{subsubsec:high-dimensional-proof}

We include a high-dimensional proof-of-concept to illustrate the use of RFM bases in a high-dimensional physical domain. 
Consider the \(d=10\) Poisson equation on \(\Omega=(-1,1)^d\):
\begin{equation}
\begin{aligned}
-\Delta u &= f \quad \text{in }\Omega,\\
u &= \alpha\sum_{i=1}^d\sin\!\left(\frac{\pi}{2}x_i\right) \quad \text{on }\partial\Omega,
\end{aligned}
\label{eq:high-d-poisson}
\end{equation}
where
$
f(\mathbf x)
=
\alpha\frac{\pi^2}{4}
\sum_{i=1}^d\sin\!\left(\frac{\pi}{2}x_i\right)$,$
u(\mathbf x)
=
\alpha
\sum_{i=1}^d\sin\!\left(\frac{\pi}{2}x_i\right).
$ The target map is
$
G:f(\mathbf x)\mapsto u(\mathbf x).
$
We use RFM bases for both input and output spaces and compare \(m_1=m_2=128\) with \(m_1=m_2=8192\).

The training data are generated by choosing eight values of \(\alpha\sim\mathcal U([-1,1])\) in the training range. 
Two unseen values of \(\alpha\) in the same range are used for interpolation tests, denoted by Gen1, and ten values of \(\alpha\sim\mathcal U([1,2])\) are used for extrapolation tests, denoted by Gen2. 
Table~\ref{tab:high_Poi_result} shows that increasing the number of RFM basis functions improves accuracy but increases training time. 
This example is intended as a proof-of-concept for high-dimensional function evaluation with prescribed mesh-free bases, rather than as a full high-dimensional operator-learning benchmark.

\begin{table}[!htpb]
    \centering
    \caption{
    High-dimensional Poisson equation.
    Mean relative \(L^2\) errors for different numbers of RFM basis functions.
    }
    \label{tab:high_Poi_result}
    \small
    \renewcommand{\arraystretch}{1.12}
    \begin{tabular}{@{}lcccccc@{}}
        \toprule
        & \multicolumn{3}{c}{\(m_1=m_2=128\)}
        & \multicolumn{3}{c}{\(m_1=m_2=8192\)} \\
        \cmidrule(lr){2-4}\cmidrule(l){5-7}
        Dimension \(d\) & Gen1 & Gen2 & Training time & Gen1 & Gen2 & Training time \\ 
        \midrule
        10 
        & 1.69e-2 
        & 3.00e-2 
        & 1,208.5 s
        & 6.47e-3 
        & 1.27e-2 
        & 76,663.5 s \\ 
        \bottomrule
    \end{tabular}
\end{table}

\begin{figure}[!ht]
    \centering
    \includegraphics[width=0.4\textwidth]{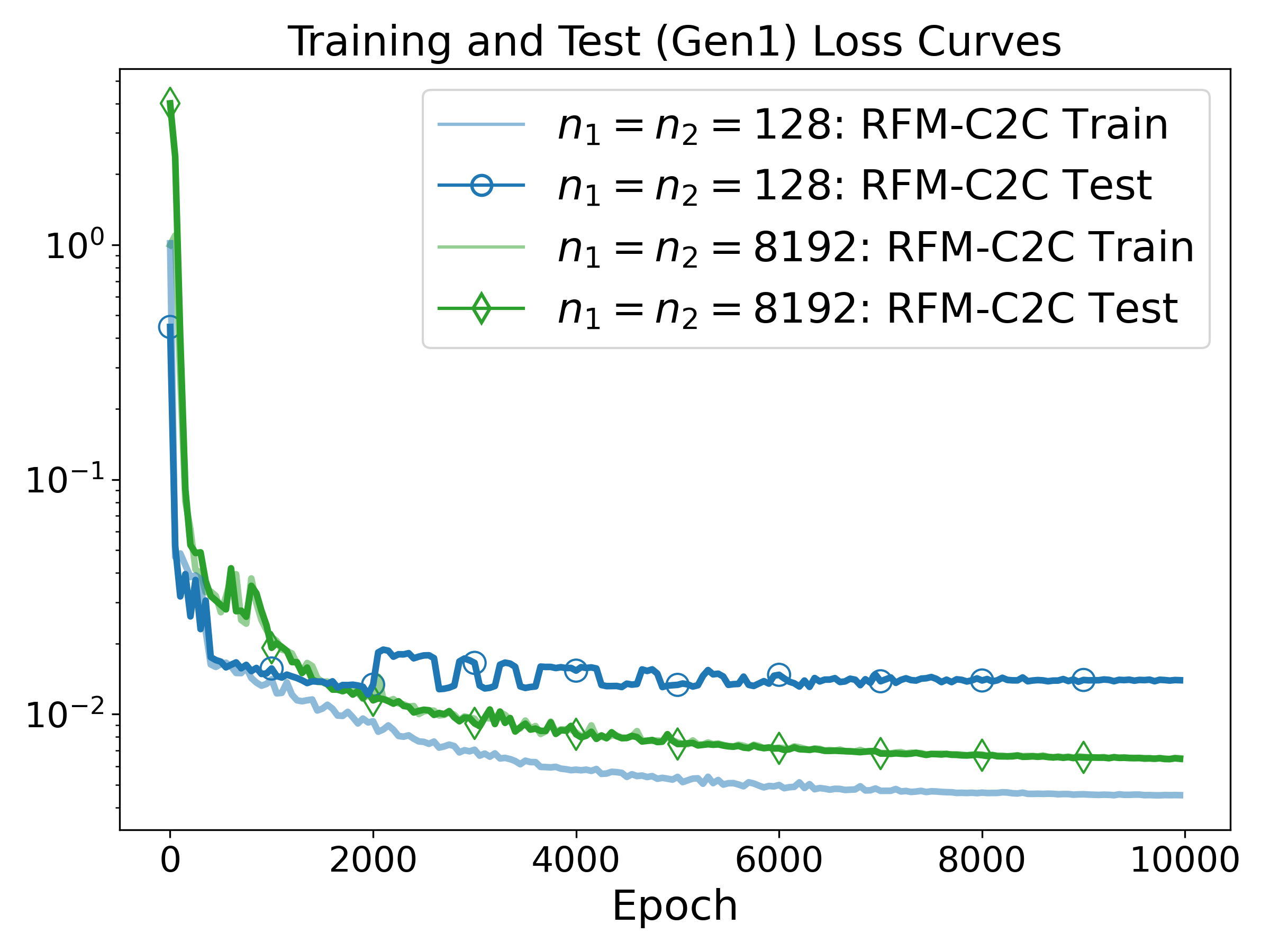}
    \caption{
    \textbf{High-dimensional Poisson equation.}
    Training curves for different numbers of RFM basis functions.
    }
    \label{fig:high_d_poi_curve}
\end{figure}

\begin{figure}[!ht]
    \centering
    \begin{minipage}[t]{0.33\textwidth}
        \centering
        \includegraphics[width=\linewidth]{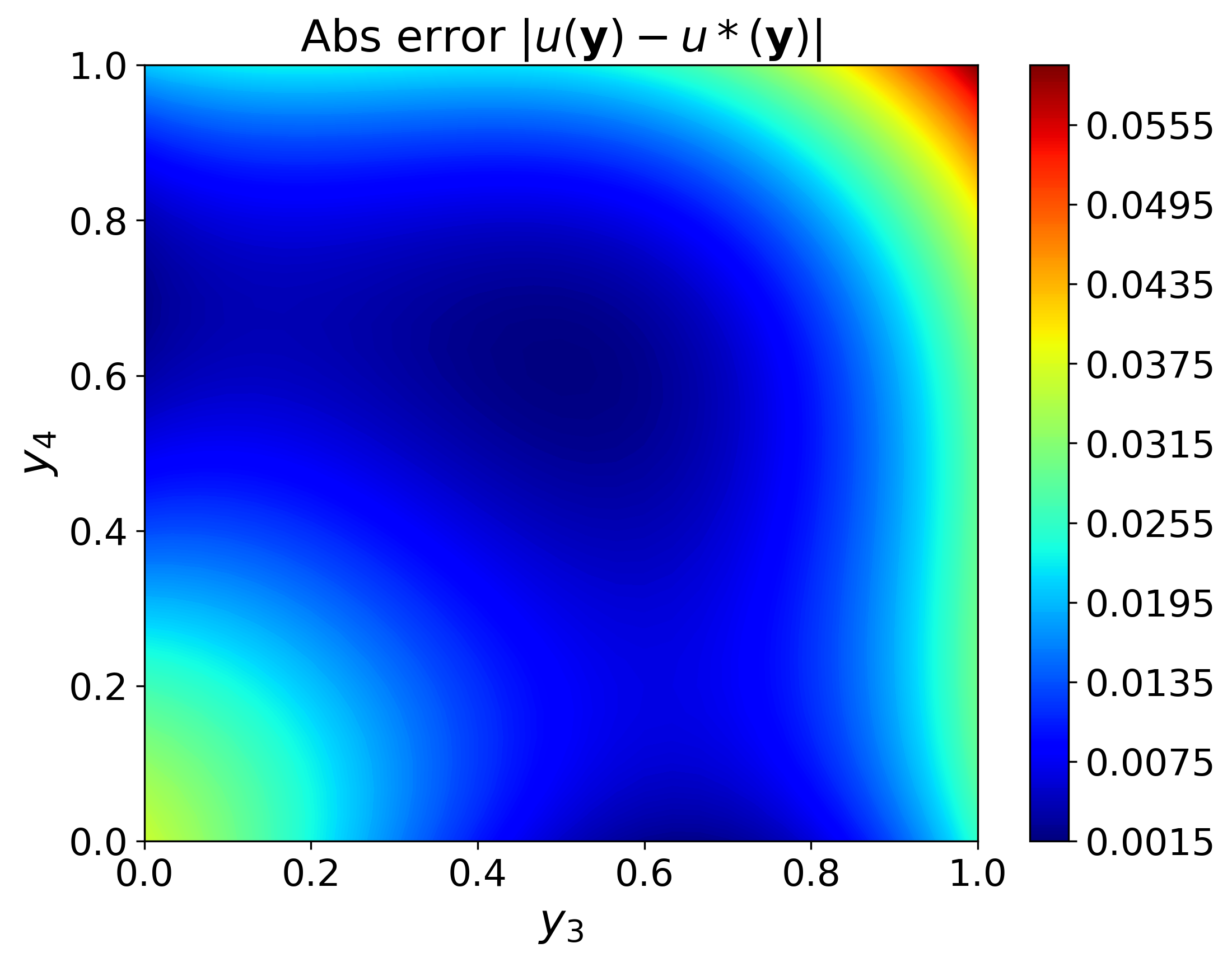}
        \caption*{(a) Pointwise error.}
    \end{minipage}
    \hfill
    \begin{minipage}[t]{0.33\textwidth}
        \centering
        \includegraphics[width=\linewidth]{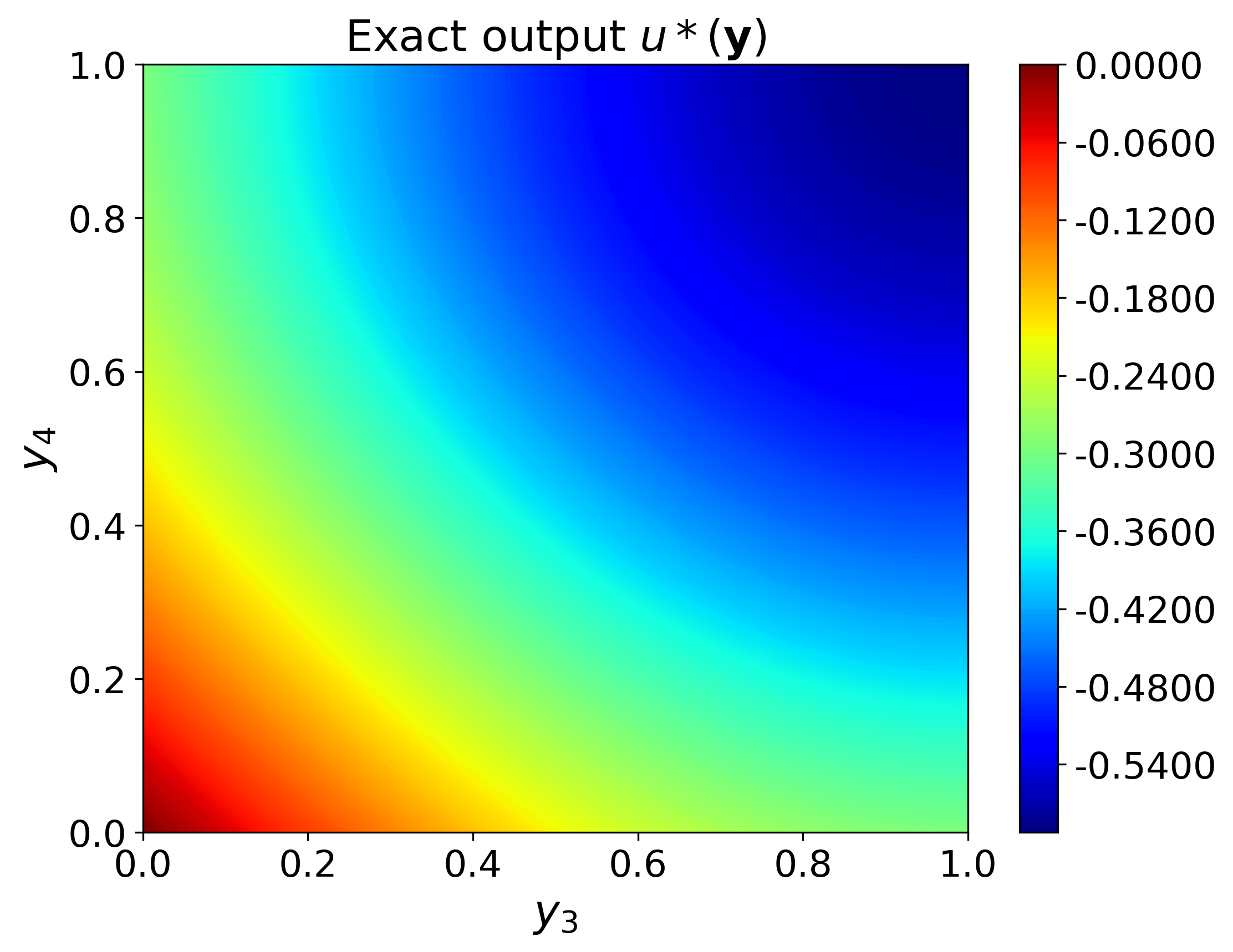}
        \caption*{(b) Exact solution.}
    \end{minipage}
    \hfill
    \begin{minipage}[t]{0.33\textwidth}
        \centering
        \includegraphics[width=\linewidth]{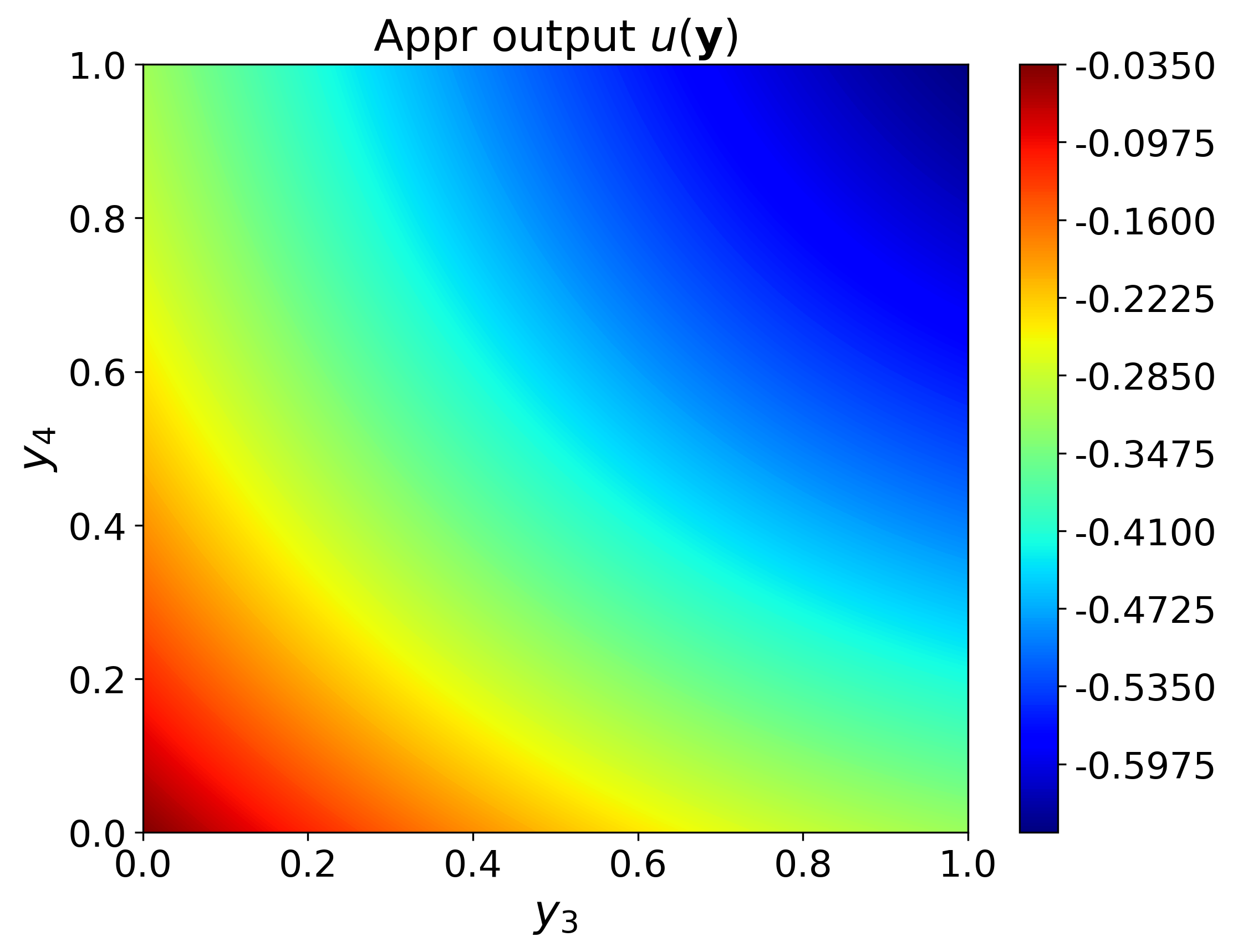}
        \caption*{(c) Predicted solution.}
    \end{minipage}
    \caption{
    \textbf{High-dimensional Poisson equation.}
    Worst-case Gen1 sample for \(m_1=m_2=8192\), visualized along the cross-section \((y_3,y_4)\) with all other coordinates fixed at zero.
    }
    \label{fig:high_d_poi_worst}
\end{figure}

\newpage
\subsection{Inverse Stokes boundary recovery from sample-dependent sparse sensors}
\label{subsec:inverse-stokes}

We finally consider an inverse boundary problem for the steady Stokes system, motivated by flow-based obstacle detection. 
The obstacle boundary is not directly accessible, and only sparse velocity measurements are available in a downstream sensing region. 
Unlike fixed-grid observations, the sensor locations are independently sampled for each realization and therefore are not aligned across different samples. 
The objective is to recover the obstacle geometry, represented by a signed distance function (SDF), from these sample-dependent sparse flow observations.

Let \(\Omega_0\subset\mathbb R^2\) be a fixed computational domain, and let \(\mathcal B\subset\Omega_0\) denote an unknown obstacle with boundary \(\Gamma=\partial\mathcal B\). 
The fluid domain is
\[
\Omega=\Omega_0\setminus\mathcal B.
\]
For each obstacle geometry, the velocity field \(\mathbf u=(u_x,u_y)^\top\) and pressure \(p\) satisfy
\begin{equation}
\label{eq:stokes_forward}
\begin{cases}
-\mu\Delta\mathbf u+\nabla p=\mathbf 0, & \mathbf x\in\Omega,\\
\nabla\cdot\mathbf u=0, & \mathbf x\in\Omega,\\
\mathbf u=\mathbf 0, & \mathbf x\in\Gamma,\\
\mathbf u=\mathbf U_\infty, & \mathbf x\in\partial\Omega_0,
\end{cases}
\end{equation}
where \(\mu>0\) is the dynamic viscosity and \(\mathbf U_\infty\) is a prescribed far-field velocity. 
The obstacle geometry is represented by the SDF
\begin{equation}
\label{eq:sdf_def}
\phi(\mathbf x)=
\begin{cases}
-\operatorname{dist}(\mathbf x,\Gamma), & \mathbf x\in\mathcal B,\\
\phantom{-}\operatorname{dist}(\mathbf x,\Gamma), & \mathbf x\in\Omega,
\end{cases}
\qquad
\Gamma=\{\mathbf x\in\Omega_0:\phi(\mathbf x)=0\}.
\end{equation}

\paragraph{Dataset and sparse observation setting.}

We take
\[
\Omega_0=[-3,3]\times[-3,3],
\]
and restrict the obstacles to the central region
\[
\Omega_{\rm obj}=[-1.5,1.5]\times[-1.5,1.5].
\]
Each obstacle consists of \(K_{\rm obs}\in\{1,2,3\}\) connected components. 
The individual components are sampled from circular, elliptical, and star-shaped families, with randomly generated centers, sizes, orientations, and shape parameters. 
Samples with excessive component overlap or components too close to the boundary of \(\Omega_{\rm obj}\) are rejected. 
The dataset contains 2500 training samples and 500 testing samples.

For each sample, \(n_s=100\) sensor locations are independently drawn from the downstream observation region
\begin{equation}
\label{eq:obs_region}
\begin{aligned}
\Omega_{\rm obs}={}&
\{(x,y):1.5\le x\le3,\ -3\le y\le3\}\\
&\cup\{(x,y):0.5\le x\le3,\ -3\le y\le-2\}\\
&\cup\{(x,y):0.5\le x\le3,\ 2\le y\le3\}.
\end{aligned}
\end{equation}
The corresponding observation set is
\begin{equation}
\label{eq:inverse_obs}
\mathcal Y
=
\left\{
\bigl(\mathbf x_j,\mathbf u(\mathbf x_j;\phi)\bigr)
\right\}_{j=1}^{n_s},
\qquad
\mathbf x_j\in\Omega_{\rm obs}.
\end{equation}
Because the sensor locations are resampled independently for each realization, different samples do not share a common sensor grid. 
The inverse operator to be learned is
\[
G_{\rm inv}:\mathcal Y\mapsto\phi.
\]
This is also a nonlocal inverse problem, since the unknown obstacle is located near the center of the domain, whereas all velocity observations are collected downstream.

\paragraph{Coefficient representation.}

The two velocity components are encoded separately using prescribed RFM bases on
\[
D_v=[0,3]\times[-3,3].
\]
We use a \(2\times2\) partition and \(M_v=64\) basis functions for each component. 
For every sample, the basis evaluation matrices are assembled at its actual sensor locations, allowing the scattered observations to be converted into coefficient vectors
\[
\mathbf w_x,\mathbf w_y\in\mathbb R^{M_v}.
\]
The neural-network input is
\begin{equation}
\label{eq:input_coeff_stokes}
\mathbf w_{\rm obs}
=
\begin{bmatrix}
\mathbf w_x\\
\mathbf w_y
\end{bmatrix}
\in\mathbb R^{2M_v}.
\end{equation}

The output SDF is represented on a \(128\times128\) grid using \(M_u=1024\) prescribed Gaussian RBFs:
\begin{equation}
\label{eq:stokes_output_rbf}
\widehat\phi(\mathbf x)
=
\sum_{m=1}^{M_u}a_m\psi_m(\mathbf x),
\qquad
\psi_m(\mathbf x)
=
\exp\left(
-\frac{\|\mathbf x-\mathbf c_m\|_2^2}{2\sigma_m^2}
\right).
\end{equation}
The resulting coefficient-to-coefficient inverse map is
\begin{equation}
\label{eq:stokes_c2c_map}
\widehat{\mathbf a}
=
F_\theta(\mathbf w_{\rm obs}),
\qquad
\widehat{\boldsymbol\phi}
=
A_u\widehat{\mathbf a},
\end{equation}
where \(A_u\) is the prescribed RBF decoder. 
Further details on obstacle generation, coefficient fitting, RBF construction, and the training objective are provided in Appendix~\ref{app:stokes_inverse_details}.

\paragraph{Numerical results}

\begin{figure}[!ht]
    \centering
    \begin{minipage}[t]{0.32\textwidth}
        \centering
        \includegraphics[width=\linewidth]{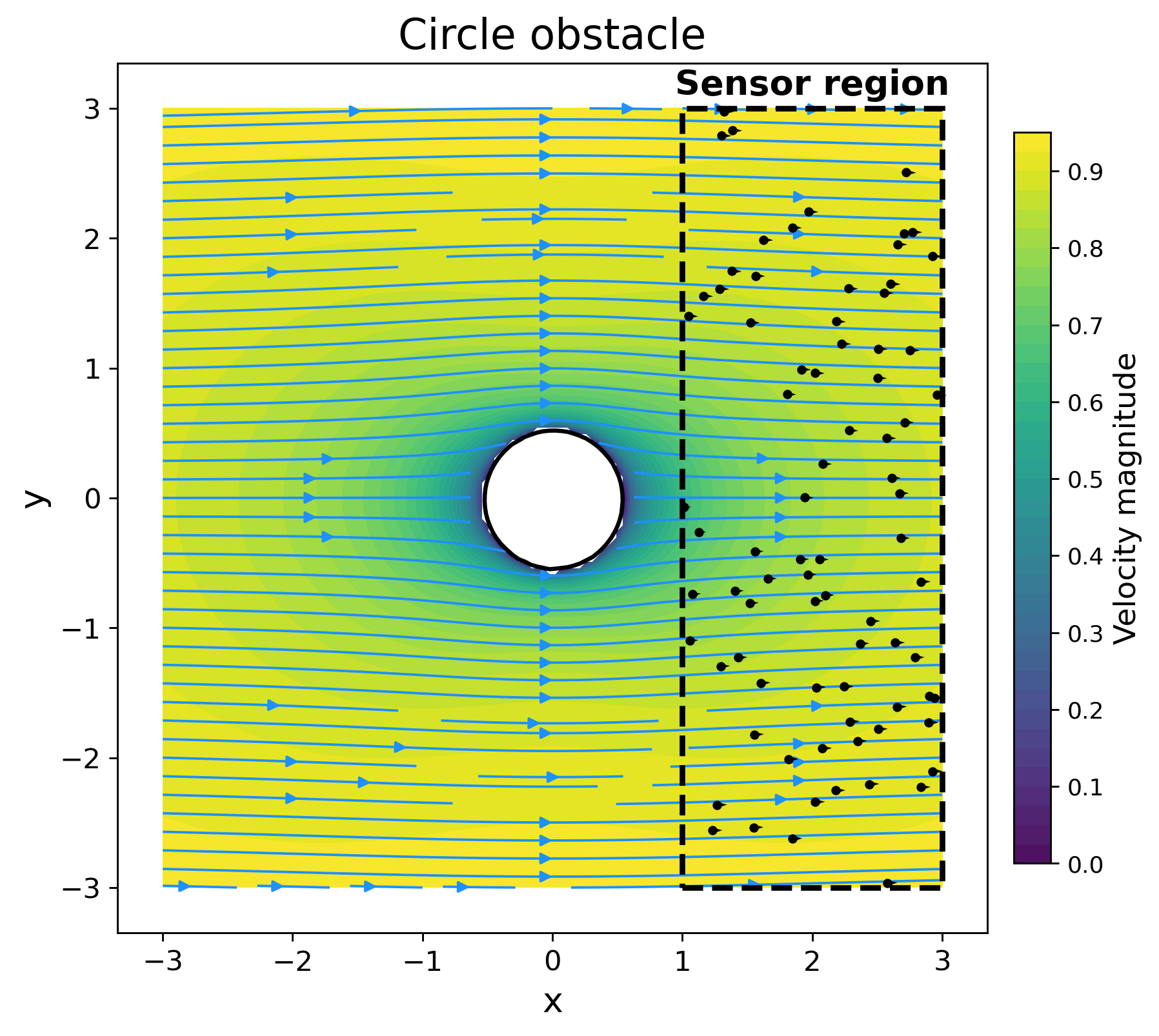}
        \caption*{(a) Circle obstacle.}
    \end{minipage}
    \hfill
    \begin{minipage}[t]{0.32\textwidth}
        \centering
        \includegraphics[width=\linewidth]{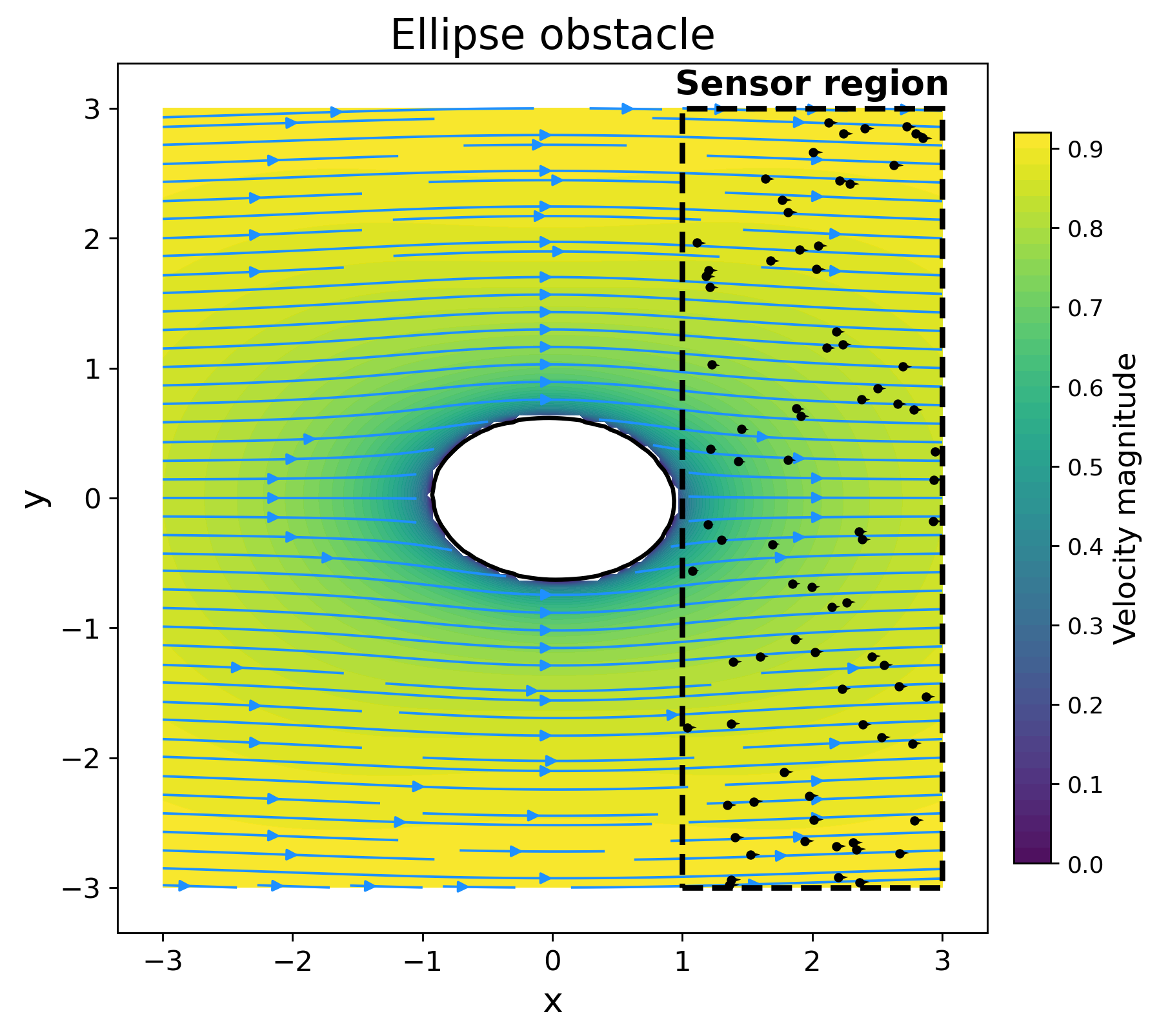}
        \caption*{(b) Ellipse obstacle.}
    \end{minipage}
    \hfill
    \begin{minipage}[t]{0.32\textwidth}
        \centering
        \includegraphics[width=\linewidth]{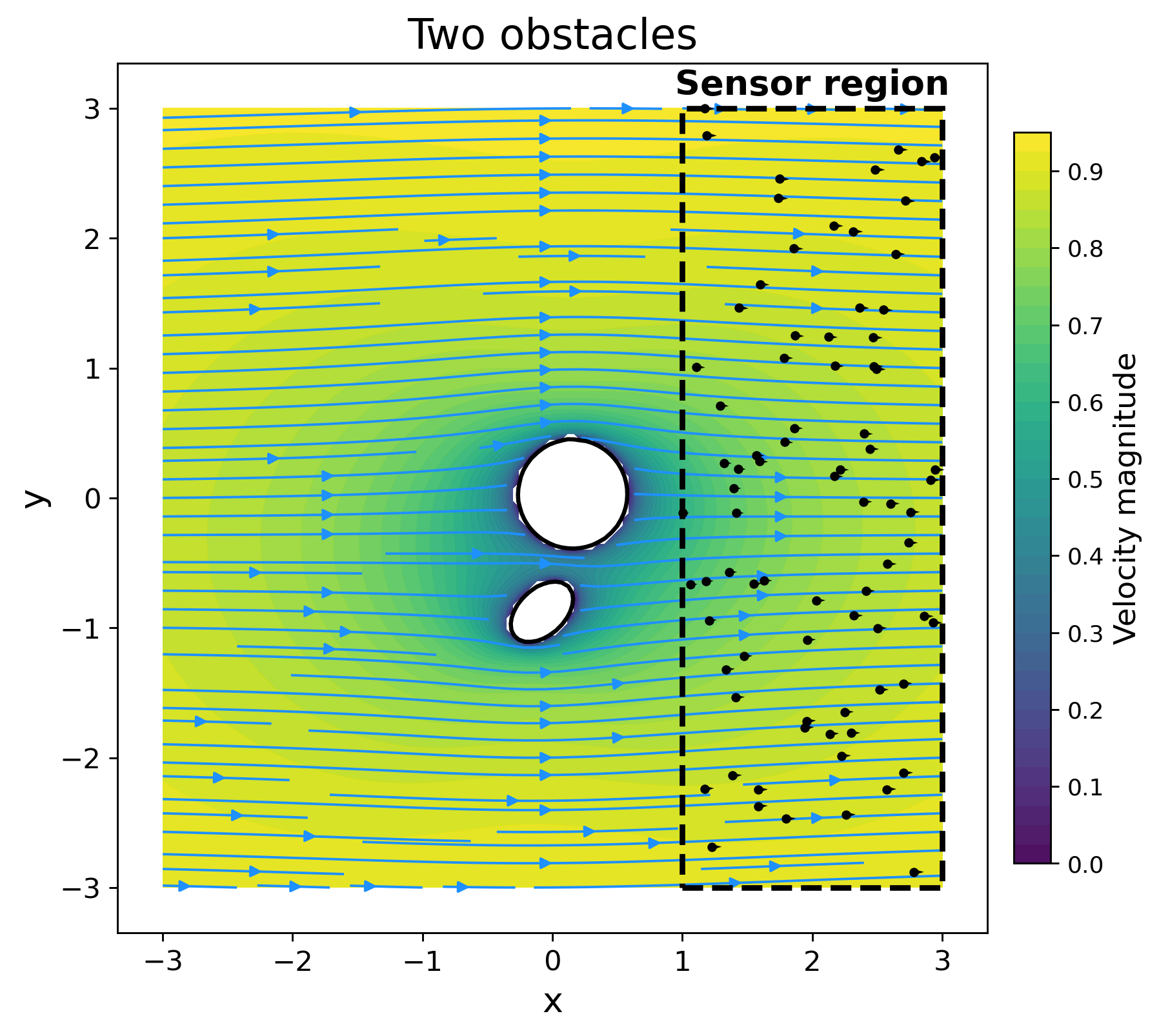}
        \caption*{(c) Two obstacles.}
    \end{minipage}
    \caption{\textbf{Forward Stokes flow and downstream sensing region.} The background color denotes velocity magnitude, streamlines indicate flow direction, and black dots denote sampled sensor locations.}
    \label{fig:stokes_forward_sensor_setting}
\end{figure}

\begin{figure}[!ht]
    \centering
    \begin{minipage}[t]{0.32\textwidth}
        \centering
        \includegraphics[width=\linewidth]{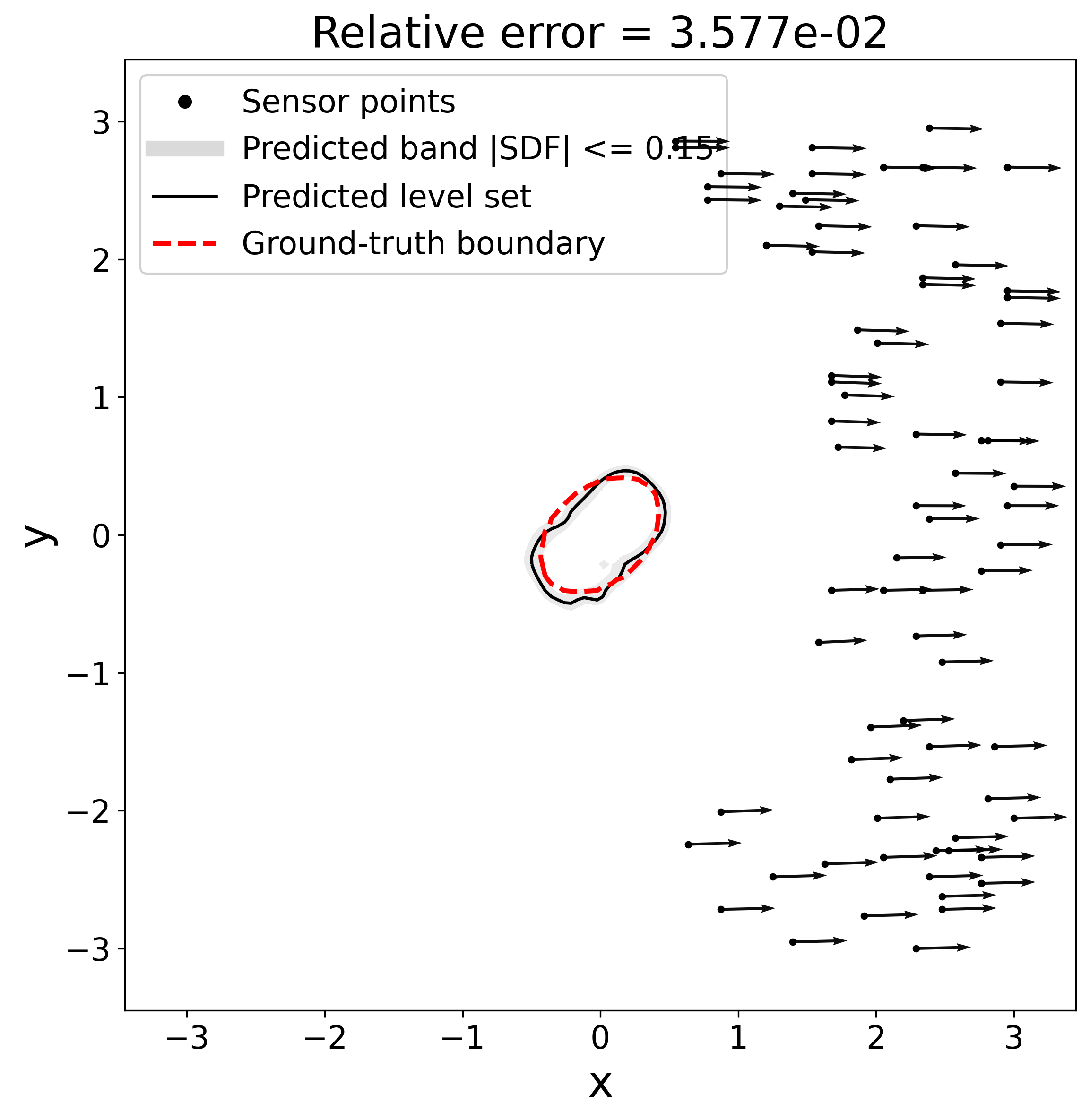}
        \caption*{(a)}
    \end{minipage}
    \hfill
    \begin{minipage}[t]{0.32\textwidth}
        \centering
        \includegraphics[width=\linewidth]{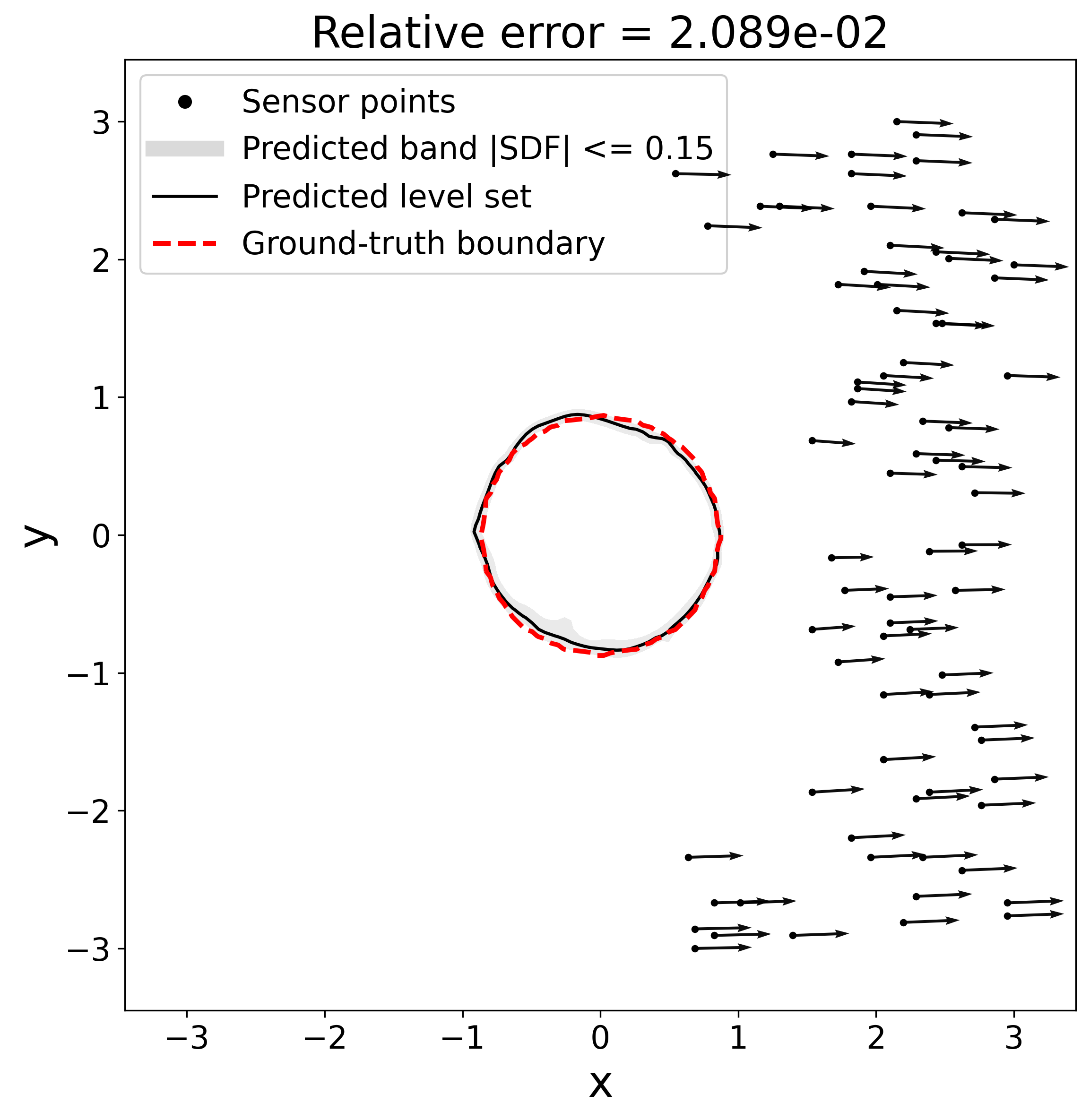}
        \caption*{(b)}
    \end{minipage}
    \hfill
    \begin{minipage}[t]{0.32\textwidth}
        \centering
        \includegraphics[width=\linewidth]{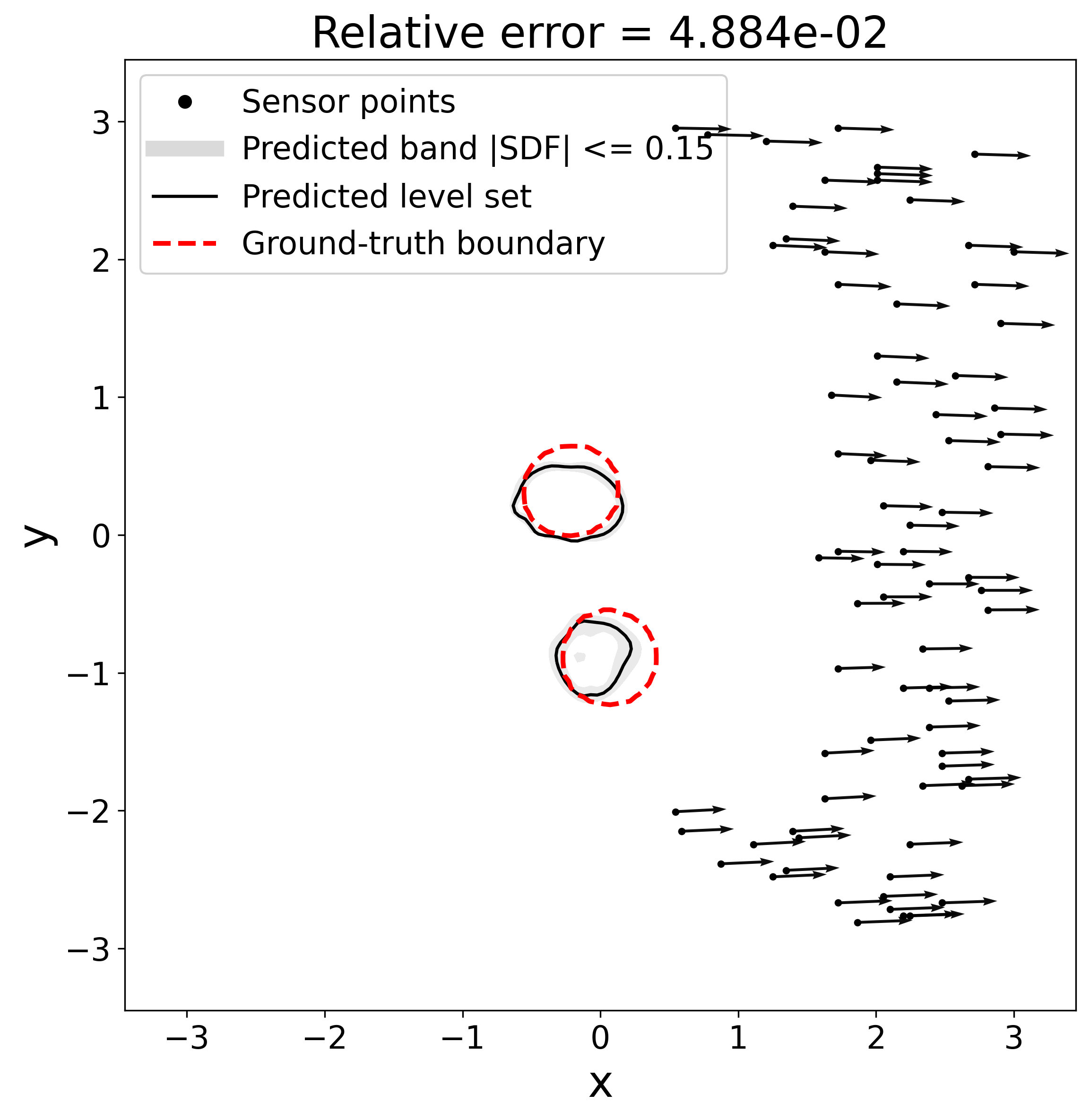}
        \caption*{(c)}
    \end{minipage}
    \caption{\textbf{Inverse obstacle recovery from sparse velocity observations.} Black dots denote sensor locations, arrows denote observed local velocity vectors, red dashed curves are ground-truth boundaries, and black solid curves are recovered zero level sets. Gray regions show the predicted near-boundary band $|\widehat\phi|\le0.15$.}
    \label{fig:stokes_inverse_reconstruction}
\end{figure}

Figure~\ref{fig:stokes_forward_sensor_setting} shows representative forward velocity fields and independently sampled sensor configurations. 
Figure~\ref{fig:stokes_inverse_reconstruction} shows representative reconstructions for single- and multi-component obstacles. 
The recovered boundary is extracted as the zero level set of the predicted SDF.

\newpage
\section{Conclusion}

We introduced the \emph{Fixed-Basis Coefficient-to-Coefficient Network} (FB-C2CNet), a neural operator framework built on prescribed, data-independent input and output approximation spaces. 
Input observations are encoded into basis coefficients through a non-trainable coefficient-recovery procedure, the neural network learns only the finite-dimensional coefficient-to-coefficient map, and the predicted field is reconstructed by a prescribed output basis. 
By separating numerical representation from neural-network optimization, FB-C2CNet reduces the dimension and complexity of the trainable model and provides a natural interface with classical approximation spaces, including random-feature, finite element, and radial-basis-function spaces.

We analyzed the main error sources in this framework. 
Regularized coefficient recovery introduces a stability--bias trade-off: stronger regularization suppresses unstable directions but may also remove informative components of the input. 
We further showed that the approximation error associated with the prescribed output space forms an intrinsic accuracy floor for any FB-C2CNet prediction. 
This decomposition clarifies the respective roles of basis approximation, coefficient recovery, and neural-network learning, and provides practical guidance for selecting basis spaces and regularization parameters.

Numerical experiments across regular and irregular domains, scalar and multi-component operators, nonlinear time-dependent equations, weak-solution problems, a high-dimensional proof-of-concept, and an inverse Stokes boundary-recovery problem demonstrate the flexibility of the proposed framework. 
The comparisons with pointwise-input variants show that coefficient-space learning can substantially reduce the neural input dimension and make the optimization problem empirically easier. 
Compared with representative neural operator architectures, FB-C2CNet achieves competitive accuracy with a markedly shorter total offline cost, particularly when only the coefficient map needs to be trained. 
The non-aligned sampling and inverse Stokes experiments further illustrate that, with suitable prescribed bases, the method can accommodate scattered and sample-dependent sensor locations without requiring a common observation grid.

The performance of FB-C2CNet ultimately depends on the approximation capability and numerical conditioning of the prescribed basis spaces. 
A poorly chosen or excessively large basis may lead to inaccurate or unstable coefficient representations and can offset the computational benefit of coefficient-space learning. 
Future work will therefore investigate adaptive and problem-dependent basis construction, hybrid prescribed-and-learned representations, noise-robust coefficient recovery, and extensions to larger-scale inverse and multiphysics problems.

\section*{Declaration of generative AI and AI-assisted technologies in the manuscript preparation process}

During the preparation of this work, the authors used OpenAI Codex for language editing and manuscript formatting. After using this tool, the authors reviewed and edited the content as needed and take full responsibility for the content of the published article.

\section*{Data and code availability}

The data and code supporting the findings of this study are available from the corresponding author upon reasonable request.

\section*{Acknowledgments}

The work of Y.X. and W.H.Z. was supported in part by the National Key Research and Development Program
of China (No. 2025YFA1016800) and the Project of Hetao Shenzhen-HKUST Innovation Cooperation
Zone HZQB-KCZYB-2020083.

% National Science Foundation.

%%%%%%%%%%%%%%%%%%%%%%%%%
%%%%%%%%%%%%%%%%%%%%%%%%%
\bibliographystyle{elsarticle-num}
\bibliography{references}
\newpage
\appendix
% \appendix

\section{Experimental Details}
\subsection{Partition-of-unity window used for the RFM basis}
\label{app:pou-window}

In one dimension, the localized partition function used in the RFM experiments is
\begin{equation}
\omega_n(x)
=
\begin{cases}
\dfrac{1+\sin(2\pi\widetilde x)}{2},
&
-\dfrac54\leq\widetilde x<-\dfrac34,\\[5pt]
1,
&
-\dfrac34\leq\widetilde x<\dfrac34,\\[5pt]
\dfrac{1-\sin(2\pi\widetilde x)}{2},
&
\dfrac34\leq\widetilde x<\dfrac54,\\[5pt]
0,
&
\text{otherwise},
\end{cases}
\qquad
\widetilde x=\frac{x-x_n}{r_n}.
\label{eq:pou-window-1d}
\end{equation}
For multidimensional domains, tensor-product partitions may be constructed as
$\omega_n(\mathbf x)=\prod_{k=1}^d\omega_{n,k}(x_k)$.

\begin{figure}[!ht]
    \centering
    \includegraphics[width=0.6\textwidth]{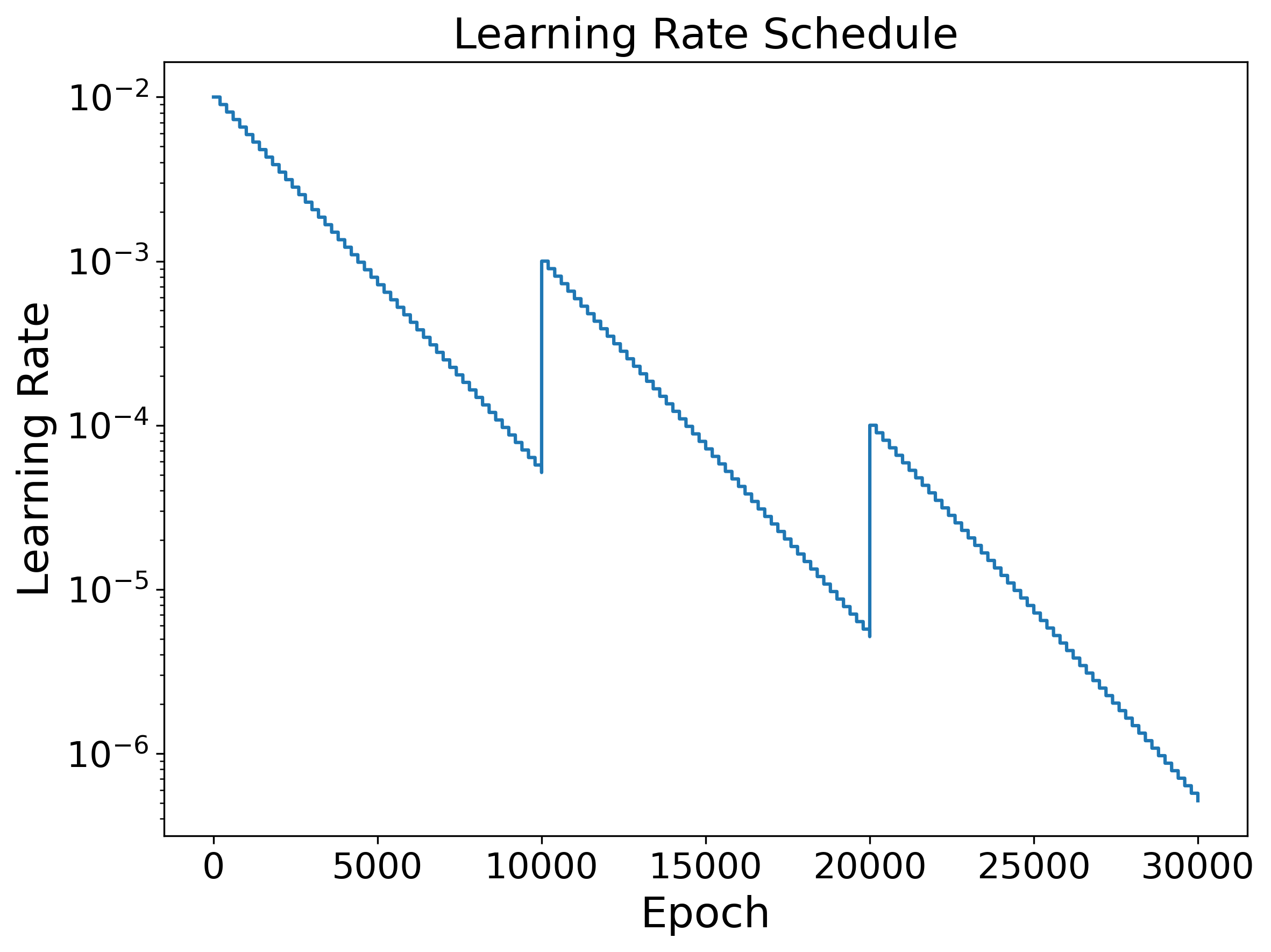}
\caption{Learning-rate schedule. The initial rates are $10^{-2}$, $10^{-3}$, and $10^{-4}$
for three consecutive 10,000-step segments,
with a decay factor of 0.9 applied every 200 steps within each segment.}

    \label{fig:learning_rate}
\end{figure}

Unless otherwise stated, the learning rate follows the schedule illustrated in Figure~\ref{fig:learning_rate}.

In our method, the quantities shown in Table~\ref{tab:key-quantities} are primarily considered. Dataset sizes and RFM basis-system parameters are summarized in Tables~\ref{tab:dataset_parameters} and~\ref{tab:rfm_parameters}, respectively.
\begin{table}[!ht]
\centering
\caption{Key quantities and notation.}
\label{tab:key-quantities}
\small
\renewcommand{\arraystretch}{1.12}
\begin{tabular}{@{}llc@{}}
\toprule
\textbf{Category} & \textbf{Quantity} & \textbf{Notation} \\
\midrule
\multirow{4}{*}{Datasets}
& Number of training samples & $N_{\text{train}}$ \\
& Number of testing samples & $N_{\text{test}}$ \\
& Number of input observation points & $n_1$ \\
& Number of output evaluation points & $n_2$ \\
\midrule
\multirow{2}{*}{Basis systems}
& Number of input basis functions & $m_1$ \\
& Number of output basis functions & $m_2$ \\
\midrule
\multirow{3}{*}{Random feature model}
& Initial-variance parameter & $R_m$ \\
& Number of partition-of-unity blocks & $M_p$ \\
& SVD cut-off & $\mathrm{cut}$ \\
\bottomrule
\end{tabular}
\end{table}

\begin{table}[!ht]
\centering
\caption{Dataset sizes and discrete input/output dimensions used in the experiments.}
\label{tab:dataset_parameters}
\small
\renewcommand{\arraystretch}{1.12}
\begin{tabular}{@{}lcccc@{}}
\toprule
Dataset & $N_{\mathrm{train}}$ & $N_{\mathrm{test}}$ & $n_1$ & $n_2$ \\
\midrule
1D Darcy flow ($n=2000$) & 800 & 200 & 2000 & 2000 \\
1D Darcy flow ($n=100$) & 800 & 200 & 100 & 100 \\
2D Darcy flow (regular) & 2400 & 600 & 19881 & 19881 \\
2D Poisson equation & 1600 & 400 & 1780 & 1780 \\
Elastic plate equation & 1850 & 100 & 101 & $1048\times2$ \\
L-shaped Darcy flow ($N_{\mathrm{train}}=1000$) & 1000 & 7000 & $736\times2$ & 450 \\
L-shaped Darcy flow ($N_{\mathrm{train}}=5100$) & 5100 & 7000 & $736\times2$ & 450 \\
L-shaped Darcy flow ($N_{\mathrm{train}}=51000$) & 51000 & 7000 & $736\times2$ & 450 \\
KdV--Burgers equation & 1600 & 400 & 101 & 19881 \\
2D Helmholtz equation & 2000 & 500 & 19881 & 19881 \\
10D Poisson equation ($m_1=m_2=128$) & 8 & 2 & $31^{10}$ & $31^{10}$ \\
10D Poisson equation ($m_1=m_2=8192$) & 8 & 2 & $31^{10}$ & $31^{10}$ \\
Inverse Stokes system & 2500 & 500 & 100 & $128^2$ \\
\bottomrule
\end{tabular}
\end{table}

\begin{table}[!ht]
\centering
\caption{Prescribed input and output basis-system parameters. The scale columns report \(R_m\) for RFM bases. For the inverse Stokes problem, the input basis is RFM, whereas the output basis is Gaussian RBF; hence, \([0.03,2.4]\) gives the range of output RBF widths. A dash indicates that the corresponding parameter is not used.}
\label{tab:rfm_parameters}
\footnotesize
\renewcommand{\arraystretch}{1.12}
\setlength{\tabcolsep}{2.5pt}
\begin{tabularx}{\textwidth}{@{}lcccc>{\raggedright\arraybackslash}X>{\raggedright\arraybackslash}Xc@{}}
\toprule
Dataset & $m_1$ & $m_2$ & Input scale & Output scale & $M_{p,\mathrm{in}}$ & $M_{p,\mathrm{out}}$ & $\mathrm{cut}$ \\
\midrule
1D Darcy ($n=2000$) & 128 & 128 & 3 & 3 & [16] & [16] & $10^{-1}$ \\
1D Darcy ($n=100$) & 128 & 128 & 3 & 3 & [8] & [4] & $10^{-2}$ \\
2D Darcy (regular) & 512 & 2048 & 0.3 & 0.3 & [8,8] & [8,8] & $10^{-1}$ \\
2D Poisson & 512 & 1600 & 3 & 3 & [2,2] & [4,4] & $10^{-2}$ \\
Elastic plate & 64 & 512 & 3 & 3 & [4] & [2,2] & $10^{-3}$ \\
L-shaped Darcy ($N_{\mathrm{train}}=1000$) & 512 & 1024 & 1 & 3 & [8,8] & [8,8] & $10^{-2}$ \\
L-shaped Darcy ($N_{\mathrm{train}}=5100$) & 512 & 1024 & 1 & 3 & [8,8] & [8,8] & $10^{-2}$ \\
L-shaped Darcy ($N_{\mathrm{train}}=51000$) & 512 & 1024 & 1 & 3 & [8,8] & [8,8] & $10^{-2}$ \\
KdV--Burgers & 32 & 1024 & 3 & 1 & [4] & [8,8] & $10^{-2}$ \\
2D Helmholtz & 256 & 512 & 3 & 1 & [4,4] & [4,4] & $10^{-8}$ \\
10D Poisson ($m=128$) & 128 & 128 & 0.3 & 0.3 & [1,1,1,1,1,1,1,1,1,1] & [1,1,1,1,1,1,1,1,1,1] & $10^{-1}$ \\
10D Poisson ($m=8192$) & 8192 & 8192 & 0.3 & 0.3 & [1,2,1,1,2,1,2,2,1,2] & [1,2,1,1,2,1,2,2,1,2] & $10^{-1}$ \\
Inverse Stokes & 128 & 1024 & 1 & [0.03,2.4] & [2,2] & -- & $10^{-2}$ \\
\bottomrule
\end{tabularx}
\end{table}

% \appendix

\subsection*{Experimental Results}

The performance of our method is listed in Tables~\ref{tab:rl2e-summary} and~\ref{tab:performance-summary}. An en dash in Table~\ref{tab:rl2e-summary} indicates that no experiment was carried out. The multiplication sign ($\times$) indicates layer repetition: $512\times5$ denotes five consecutive layers, each having 512 neurons.
% \begin{table}[!ht]
% \centering
% \caption{Mean Relative Error $L^2$ of the testing data.}
% \begin{tabular}{lcc}
% \hline
% &\textbf{FEM-C2C} & \textbf{RFM-C2C} \\
% \hline
% % DeepCRFM_Darcy1D.ipynb
% % DeepCFEM_Darcy1D.ipynb
% \textcolor{red}{1D Darcy Flow} & 4.26e-2 & \textbf{1.20e-2}\\
% % 3000_RL2E_cutoff4_exp_141*141_FEM_test.ipynb
% % 2048_RL2E_cutoff4_exp_NN_test.ipynb
% 2D Darcy Flow (Regular)&7.04e-3 &\textbf{6.07e-3} \\
% % 2000_RL2E_82_cutoff4_FEM.ipynb
% % 2000_RL2E_82_cutoff4_NN.ipynb
% 2D Poisson Equation& 8.41e-3& \textbf{6.40e-3}  \\
% % 1950_RL2E_uxuy_*1e3_cutoff3_vector2_FEM.ipynb
% % 1950_RL2E_uxuy_*1e3_cutoff3_vector2_NN.ipynb
% Elastic Plate Equation& 6.26e-1 & \textbf{8.92e-3}\\
% % Darcy_Lshape_FEM_test11.py
% % 58000_RL2E_cutoff2_vector_2_NN_1024_1024.ipynb
% 2D Darcy Flow (Lshaped)& 2.94e-1 &\textbf{3.71e-2} \\
% 10D Poisson Equation&- &\textbf{6.48e-3} \\
% \hline
% \end{tabular}
% \end{table}

\begin{table}[!ht]
\centering
\small
\caption{Mean relative $L^2$ error of the testing data. The inverse Stokes entry uses an RFM input basis and a Gaussian RBF output basis.}
\renewcommand{\arraystretch}{1.12}
\begin{tabular}{@{}lcccc@{}}
\toprule
Dataset & Scalar FEM-C2C & Scalar RFM-C2C & Vector FEM-C2C & Vector RFM-C2C\\
\midrule
% DeepCRFM_Darcy1D.ipynb
% DeepCFEM_Darcy1D.ipynb
1D Darcy flow ($n_1=n_2=2000$) & 1.367e-2 & \textbf{3.87e-3} & -- & --\\
% 3000_RL2E_cutoff4_exp_141*141_FEM_test.ipynb
% 2048_RL2E_cutoff4_exp_NN_test.ipynb
2D Darcy flow (regular)&6.92e-3 &\textbf{6.47e-3} & -- & --\\
% 2000_RL2E_82_cutoff4_FEM.ipynb
% 2000_RL2E_82_cutoff4_NN.ipynb
2D Poisson equation& \textbf{4.13e-3}& 7.69e-3  & -- & --\\
% 1950_RL2E_uxuy_*1e3_cutoff3_vector2_FEM.ipynb
% 1950_RL2E_uxuy_*1e3_cutoff3_vector2_NN.ipynb
Elastic plate equation& \textbf{1.68e-3} & 8.47e-3 & 6.2443e-1 & 7.81e-3\\
% Darcy_Lshape_FEM_test11.py
% 58000_RL2E_cutoff2_vector_2_NN_1024_1024.ipynb
L-shaped Darcy flow& -- & -- & -- & \textbf{3.855e-2}\\
KdV--Burgers equation& -- &3.29e-3 & -- & --\\
2D Helmholtz equation & \textbf{2.13e-3} &6.51e-3 & -- & --\\
10D Poisson equation& -- &\textbf{6.48e-3} & -- & --\\
Inverse Stokes system & -- & \textbf{5.56e-2} & -- & -- \\
\bottomrule
\end{tabular}
\label{tab:rl2e-summary}
\end{table}

% \begin{table}[h!]
% \centering
% \caption{Mean Relative Error $L^2$.}
% \begin{tabular}{lll}
% \hline
% &\textbf{FEM-C2C} & \textbf{RFM-C2C} \\
% \hline
% % DeepCRFM_Darcy1D.ipynb
% % DeepCFEM_Darcy1D.ipynb
% 1D Darcy Flow& 0.042612 & 0.012040\\
% % 3000_RL2E_cutoff4_exp_141*141_FEM_test.ipynb
% % 2048_RL2E_cutoff4_exp_NN_test.ipynb
% 2D Darcy Flow in Regular domain&0.007044 &0.006067 \\
% % 2000_RL2E_82_cutoff4_FEM.ipynb
% % 2000_RL2E_82_cutoff4_NN.ipynb
% 2D Poisson in Complex domain& 0.008412& 0.006404  \\
% % 1950_RL2E_uxuy_*1e3_cutoff3_vector2_FEM.ipynb
% % 1950_RL2E_uxuy_*1e3_cutoff3_vector2_NN.ipynb
% Elastic Plate& 0.626444 & 0.008916\\
% % Darcy_Lshape_FEM_test11.py
% % 58000_RL2E_cutoff2_vector_2_NN_1024_1024.ipynb
% 2D Darcy Flow in Complex domain& 0.294052 &0.037057 \\
% High Dimension Poisson equation&- &0.006481 \\
% \hline
% \end{tabular}
% \end{table}
\begin{table}[!htbp]
\centering
\footnotesize
\caption{Performance comparison across benchmarks. Performance is measured using relative \(L^2\) error (RL2E) and mean squared error (MSE). The ``Network architecture'' column includes all trainable subnetworks. For BasisONet and B2B, these are the input basis generator, transformation network, and output basis generator. For DeepONet, these are the branch network, trunk network, and bias term in the dot-product output.}
\label{tab:IR_darcy_comparison}
\renewcommand{\arraystretch}{1.1}
\setlength{\tabcolsep}{3pt}
\begin{tabularx}{\textwidth}{@{}lcccc>{\raggedright\arraybackslash}X@{}}
\toprule
% \tiny
Method & Time & Parameters & Test RL2E & MSE & Network architecture\\
\midrule
\multicolumn{6}{@{}l}{\textbf{2D Darcy flow}}\\
% BO & 87m38.8s & 2,712,174 & 0.703e-2& 0.9e-7 &[2,512$\times$ 5,50]+[50,512$\times$ 3,30]+[2,512$\times$ 5,30]\\
% BO & 90m28.0s & 2,712,174 & 0.6379e-2& 0.73899e-7 &[2,512$\times$ 5,50]+[50,512$\times$ 3,30]+[2,512$\times$ 5,30]\\
BasisONet & 88 min 0.7 s & 2,712,174 & 7.02e-3& 9.12e-8 &[2,512$\times$5,50] + [50,512$\times$3,30] + [2,512$\times$5,30]\\
B2B & 94 min 2 s & 499,758 & 1.17e-2& 2.73e-7 & [2,256$\times$3,100] + [100,256$\times$3,101] + [2,256$\times$3,101]\\
% DeepOnet & 12m24s & 504,926 & 13.013e-2 & 3335.17e-7 & [19881,25$\times$ 3,100]+[2,25$\times$ 3,100] + 1 \\
% DeepOnet & 35s & 504,926 & 12.803e-2 & 3197.64e-8 & [19881,25$\times$ 3,100]+[2,25$\times$ 3,100] + 1 \\
DeepONet & 4 min 10 s & 504,926 & 1.15e-1 & 2.59e-5 & [19881,25$\times$3,100] + [2,25$\times$3,100] + 1 \\
% FNO2d & 762m18.8s & 1,188,353 & \textbf{0.86e-3}& 0.52e-8 & modes = 12, width = 32\\
% FNO2d & 73m21.5s & 1,188,353 & \textbf{3.46e-3}& 2.35e-08 & modes = 12, width = 32\\
FNO2D & 31 min 22 s & 1,188,353 & \textbf{5.43e-3}& \textbf{5.19e-8} & modes = 12, width = 32, epochs = 500\\
% C2C-RFM & 6.14s+17m53.7s & \textbf{492,864} & 0.6577e-2 &0.825e-7& 0 + [512,180$\times$ 2,2048] + 0\\
RFM-C2C & \textbf{5.32 s + 3 min 17 s} & \textbf{492,864} & 6.47e-3 &8.68e-8& [512,180$\times$2,2048]\\
FEM-C2C & 7 min 27 s + 4 min 51 s & 504,979 & 6.92e-3 &9.03e-8& [553,180$\times$2,2059]\\
\midrule
\multicolumn{6}{@{}l}{\textbf{2D Poisson equation}}\\
BasisONet & 17 min 52.09 s & 2,712,174 &\textbf{1.49e-3} & \textbf{7.00e-10} &[2,512$\times$5,50] + [50,512$\times$3,30] + [2,512$\times$5,30]\\
% B2B & 88m1s & 499,758 & \textcolor{red}{3e-2}& 28.64e-8 & [2,256$\times$ 3,100]+[100,256$\times$ 3,101]+[2,256$\times$ 3,101]\\
B2B & 88 min 1 s & 499,758 & 3.03e-2& 2.62e-7 & [2,256$\times$3,100] + [100,256$\times$3,101] + [2,256$\times$3,101]\\
% B2B & 12m38s & 499,758 & 4.70289e-2& 52.695e-8 & [2,256$\times$ 3,100]+[100,256$\times$ 3,101]+[2,256$\times$ 3,101]\\
% DeepOnet& 10m12s & 501,417 & 5.568e-2 & 76.80e-8 & [1780,184$\times$ 3,100]+[2,184$\times$ 3,100] + 1 \\
% DeepOnet& 36s & 501,417 & 13.87027e-2 & 471.0875e-8 & [1780,184$\times$ 3,100]+[2,184$\times$ 3,100] + 1 \\
DeepONet& 3 min 57 s & 501,417 & 6.80e-2 & 1.08e-6 & [1780,184$\times$3,100] + [2,184$\times$3,100] + 1 \\
RFM-C2C & \textbf{2.61 s + 51.93 s} & \textbf{489,640} & 7.69e-3 &1.48e-8& [512,210$\times$2,1600]\\
FEM-C2C & 34.74 s + 53.66 s & 491,096 & 4.11e-3 &4.90e-9& [533,210$\times$2,1586]\\
\midrule
\multicolumn{6}{@{}l}{\textbf{Elastic plate}}\\
% \textcolor{red}{63m11s}
% BO & 30m16.04s & 2,716,832  &0.1751e-2 &  &[1,512$\times$ 5,10]+[10,512$\times$ 3,50]+[2,512$\times$ 5,100]\\
BasisONet & 30 min 16.04 s & 2,716,832  &1.70e-3 & 1.85e-14 &[1,512$\times$5,10] + [10,512$\times$3,50] + [2,512$\times$5,100]\\
B2B & 38 min 20 s & 500,644 & 4.24e-3 & 3.36e-14  & [2,249$\times$3,100] + [100,249$\times$3,202] + [2,249$\times$3,101]\\
DeepONet& 4 min 15 s & 498,657 & 8.78e-3 & 1.10e-13 & [101,295$\times$3,200] + [2,295$\times$3,200] + 2\\
RFM-C2C (Vector) & \textbf{0.82 s + 42.11 s} & \textbf{295,936} & 7.81e-3 &4.11e-14& [64,512$\times$1,512]\\
FEM-C2C (Scalar) & 11.75 s + 47.24 s & 587,320 & \textbf{1.68e-3} &\textbf{3.80e-15}& [64,512$\times$1,1080]\\
\midrule
\multicolumn{6}{@{}l}{\textbf{L-shaped Darcy flow}}\\
% BO & 87m38.8s & 2,712,174 & 0.703e-2& 0.9e-7 &[2,512$\times$ 5,50]+[50,512$\times$ 3,30]+[2,512$\times$ 5,30]\\
B2B & 41 min 27 s & 500,643 & 7.49e-2&3.34e-7  & [2,249$\times$3,200] + [200,249$\times$3,101] + [2,249$\times$3,101]\\
DeepONet& 4 min 09 s & 498,985 & 1.69e-1 & 1.59e-6 & [1922,176$\times$3,100] + [2,176$\times$3,100] + 1 \\
RFM-C2C & 1.09 s + 2 min 36.48 s & \textbf{460,288} & \textbf{3.86e-2} &\textbf{8.55e-8}& [512,256$\times$2,1024]\\
\midrule
\multicolumn{6}{@{}l}{\textbf{KdV--Burgers equation}}\\
BasisONet & 28 min 33.3 s & 2,691,182 &3.63e-3 & 1.43e-6 &[1,512$\times$5,10] + [10,512$\times$3,50] + [2,512$\times$5,50]\\
% B2B & 42m57s & 499,502 & 1.21e-2& 2.57e-6 & [1,256$\times$ 3,100]+[100,256$\times$ 3,101]+[2,256$\times$ 3,101]\\
B2B & 42 min 56 s & 499,502 & 1.21e-2& 2.57e-6 & [1,256$\times$3,100] + [100,256$\times$3,101] + [2,256$\times$3,101]\\
% DeepOnet& 13m13s & 500110 & 5.77e-2 & 63.98e-6 & [101,317$\times$ 3,100]+[2,317$\times$ 3,100]\\
DeepONet& 10 min 45 s & 500,110 & 5.77e-2& 6.40e-5 & [101,317$\times$3,100] + [2,317$\times$3,100] + 1\\
RFM-C2C & \textbf{2.40 s + 3 min 40.38 s} & \textbf{498,724} & \textbf{3.29e-3} &\textbf{5.20e-7}& [32,300$\times$3,1024]\\
\midrule
\multicolumn{6}{@{}l}{\textbf{2D Helmholtz equation}}\\
RFM-C2C & \textbf{1.50 s + 50.64 s} & \textbf{492,287} & 6.51e-3 &4.26e-5& [256,415$\times$2,512]\\
FEM-C2C & 69.54 s + 55.01 s & 493,758 & \textbf{2.13e-3} &\textbf{4.20e-6}& [273,400$\times$2,558]\\
\midrule
\multicolumn{6}{@{}l}{\textbf{Inverse Stokes system}}\\
RFM--RBF C2C & 5.60 s + 605 s & 854,016 & 5.56e-2 & 6.24e-1  & [128,512$\times$2,1024]\\
\bottomrule
\end{tabularx}
\label{tab:performance-summary}
% \vspace{0.5cm}
% \footnotesize
% \begin{itemize}
% \item[—] MWT (Gupta et al., 2021) only supports resolution with powers of two.
% \item[—] FNO2D, U-NO and MWT's performance are further improved from originally reported because of the usage of \( H^1 \) loss and scheduler.
% \item[—] The runtime and number of parameters count are using Darcy rough as the example case.
% \end{itemize}
\end{table}

\clearpage
\section{Implementation Details for the Stokes Inverse Problem}
\label{app:stokes_inverse_details}

This appendix provides the implementation details for the inverse Stokes example in
Section~\ref{subsec:inverse-stokes}. 
We describe the sampling of obstacle geometries, sparse sensor locations, fixed input and output bases, and the training objective.

\subsection{Dataset generation}

All samples are generated on the computational domain
\[
    \Omega_0=[-3,3]\times[-3,3].
\]
For each sample, we first generate an obstacle geometry $\mathcal B$ contained in the central region
\[
    \Omega_{\rm obj}
    =
    [-1.5,1.5]\times[-1.5,1.5].
\]
The obstacle may consist of one or several connected components. 
In the experiments, we include circular, elliptical, and star-shaped components. 
Each component is generated by sampling its center, size, orientation, and shape parameters subject to the constraint that the resulting obstacle remains inside $\Omega_{\rm obj}$. 
When multiple components are present, we reject samples with excessive overlap or components that are too close to the boundary of $\Omega_{\rm obj}$.

More precisely, for each sample, we draw the number of components
\[
    K_{\rm obs}\in\{1,2,3\}.
\]
For each component, the center $\mathbf q_k$ is sampled from a central box
\[
    \mathbf q_k\sim {\rm Unif}(\Omega_{\rm obj}^{\rm center}),
\]
where $\Omega_{\rm obj}^{\rm center}\subset \Omega_{\rm obj}$ is chosen so that the whole component remains inside $\Omega_{\rm obj}$. 
The component shape is then sampled from one of the following families.

\paragraph{Circular obstacles.}
A circular component is specified by a center $\mathbf q_k$ and a radius $r_k$,
\[
    \mathcal B_k
    =
    \{\mathbf x:\|\mathbf x-\mathbf q_k\|_2\leq r_k\},
    \qquad
    r_k\sim {\rm Unif}(r_{\min},r_{\max}).
\]

\paragraph{Elliptical obstacles.}
An elliptical component is specified by a center $\mathbf q_k$, semi-axes $a_k,b_k$, and an orientation angle $\theta_k$,
\[
    a_k\sim {\rm Unif}(a_{\min},a_{\max}),
    \qquad
    b_k\sim {\rm Unif}(b_{\min},b_{\max}),
    \qquad
    \theta_k\sim {\rm Unif}(0,\pi).
\]
Let $R_{\theta_k}$ be the two-dimensional rotation matrix. 
The component is defined by
\[
    \mathcal B_k
    =
    \left\{
    \mathbf x:
    \left(
    \frac{[(R_{\theta_k}^{\top}(\mathbf x-\mathbf q_k))]_1}{a_k}
    \right)^2
    +
    \left(
    \frac{[(R_{\theta_k}^{\top}(\mathbf x-\mathbf q_k))]_2}{b_k}
    \right)^2
    \leq 1
    \right\}.
\]

\paragraph{Star-shaped obstacles.}
A star-shaped component is represented in polar coordinates around its center $\mathbf q_k$ by
\[
    \rho_k(\vartheta)
    =
    r_{0,k}
    \left(
        1+\beta_k \cos(m_k\vartheta+\varphi_k)
    \right),
\]
where
\[
    r_{0,k}\sim {\rm Unif}(r_{0,\min},r_{0,\max}),
    \qquad
    \beta_k\sim {\rm Unif}(\beta_{\min},\beta_{\max}),
    \qquad
    \varphi_k\sim {\rm Unif}(0,2\pi),
\]
and $m_k$ is sampled from a prescribed finite set of modes. 
The corresponding component is
\[
    \mathcal B_k
    =
    \left\{
    \mathbf q_k+\rho(\cos\vartheta,\sin\vartheta):
    0\leq \rho\leq \rho_k(\vartheta),
    \ 0\leq \vartheta<2\pi
    \right\}.
\]

The full obstacle is the union
\[
    \mathcal B
    =
    \bigcup_{k=1}^{K_{\rm obs}}\mathcal B_k.
\]
After generating $\mathcal B$, we compute its signed distance function
\[
    \phi(\mathbf x)
    =
    \begin{cases}
    -\operatorname{dist}(\mathbf x,\partial\mathcal B), & \mathbf x\in \mathcal B,\\
    \phantom{-}\operatorname{dist}(\mathbf x,\partial\mathcal B), & \mathbf x\in \Omega_0\setminus \mathcal B.
    \end{cases}
\]
The forward Stokes solver is then applied to obtain the corresponding velocity field
\[
    \mathbf u(\mathbf x;\phi)
    =
    (u_x(\mathbf x;\phi),u_y(\mathbf x;\phi))^\top.
\]

\subsection{Sampling of sparse velocity observations}

For each sample, sparse sensor locations are drawn from the downstream observation region
\[
\begin{aligned}
    \Omega_{\rm obs}
    ={}&
    \{(x,y): 1.5\leq x\leq 3,\ -3\leq y\leq 3\} \\
    &\cup
    \{(x,y): 0.5\leq x\leq 3,\ -3\leq y\leq -2\} \\
    &\cup
    \{(x,y): 0.5\leq x\leq 3,\ 2\leq y\leq 3\}.
\end{aligned}
\]
The number of sensors is fixed to $n_s$ for each sample. 
Unless otherwise stated, the sensor points are sampled uniformly from $\Omega_{\rm obs}$. 
Equivalently, one may first select one of the three rectangular subregions with probability proportional to its area and then sample uniformly inside the selected rectangle. 
The sparse observation set is
\[
    \mathcal Y
    =
    \left\{
    \bigl(\mathbf x_j,\mathbf u(\mathbf x_j;\phi)\bigr)
    \right\}_{j=1}^{n_s}.
\]
Here
\[
    \mathbf u(\mathbf x_j;\phi)
    =
    \bigl(
    u_x(\mathbf x_j;\phi),
    u_y(\mathbf x_j;\phi)
    \bigr)^\top.
\]
In the reported visualizations, we also plot the observed velocity vectors at the sampled sensor locations.

\subsection{Fixed output RBF basis for the SDF}

The output SDF is represented on the full grid by a fixed Gaussian RBF expansion. 
Let $\{\mathbf z_i\}_{i=1}^{N_g}$ be the output grid points. 
We use
\[
    M_u=1024
\]
Gaussian RBFs
\[
    \psi_m(\mathbf x)
    =
    \exp\left(
        -\frac{\|\mathbf x-\mathbf c_m\|_2^2}{2\sigma_m^2}
    \right),
    \qquad
    m=1,\dots,M_u.
\]
The predicted SDF takes the form
\[
    \widehat\phi(\mathbf x)
    =
    \sum_{m=1}^{M_u}a_m\psi_m(\mathbf x).
\]
On the output grid, this gives
\[
    \widehat{\boldsymbol\phi}
    =
    A_u\mathbf a,
    \qquad
    (A_u)_{im}=\psi_m(\mathbf z_i),
    \qquad
    A_u\in\mathbb R^{N_g\times M_u}.
\]

Since the obstacle is known to lie in the central part of the domain, we use a prior-aware sampling strategy for the RBF centers. 
Specifically, let
\[
    M_{\rm center}=\lfloor 0.75 M_u\rfloor,
    \qquad
    M_{\rm global}=M_u-M_{\rm center}.
\]
We sample
\[
    \mathbf c_m\sim {\rm Unif}\bigl([-1.75,1.75]\times[-1.75,1.75]\bigr),
    \qquad
    m=1,\dots,M_{\rm center},
\]
and
\[
    \mathbf c_m\sim {\rm Unif}(\Omega_0),
    \qquad
    m=M_{\rm center}+1,\dots,M_u.
\]
The first group gives higher resolution near the possible obstacle region, while the second group preserves global coverage of the SDF field.

The RBF widths are fixed after random initialization. 
In our implementation, they are sampled independently from a prescribed interval
\[
    \sigma_m\sim {\rm Unif}(\sigma_{\min},\sigma_{\max}),
\]
or set according to the average spacing of the sampled centers. 
The same fixed output basis $A_u$ is used for all training and testing samples.

\subsection{Fixed input RFM representation}

The sparse velocity observations are converted into fixed-dimensional input coefficients using separate RFM bases for the two velocity components. 
We construct two independent fixed bases
\[
    \{\eta^{(x)}_k\}_{k=1}^{M_v},
    \qquad
    \{\eta^{(y)}_k\}_{k=1}^{M_v},
\]
for $u_x$ and $u_y$, respectively. 
The RFM bases are built on the velocity observation domain
\[
    D_v=[0,3]\times[-3,3],
\]
using a $2\times 2$ partition and
\[
    M_v=64
\]
basis functions for each component.

Given sensor locations $\{\mathbf x_j\}_{j=1}^{n_s}$, define the RFM design matrices
\[
    (A_x)_{jk}
    =
    \eta^{(x)}_k(\mathbf x_j),
    \qquad
    (A_y)_{jk}
    =
    \eta^{(y)}_k(\mathbf x_j),
    \qquad
    j=1,\dots,n_s,\quad k=1,\dots,M_v.
\]
Let
\[
    \mathbf u_x
    =
    \bigl(u_x(\mathbf x_1),\dots,u_x(\mathbf x_{n_s})\bigr)^\top,
    \qquad
    \mathbf u_y
    =
    \bigl(u_y(\mathbf x_1),\dots,u_y(\mathbf x_{n_s})\bigr)^\top.
\]
The input coefficients are obtained by ridge regression:
\begin{equation}
\label{eq:app_ridge_x}
    \mathbf w_x
    =
    \arg\min_{\mathbf w\in\mathbb R^{M_v}}
    \|A_x\mathbf w-\mathbf u_x\|_2^2
    +
    \lambda_v\|\mathbf w\|_2^2,
\end{equation}
and
\begin{equation}
\label{eq:app_ridge_y}
    \mathbf w_y
    =
    \arg\min_{\mathbf w\in\mathbb R^{M_v}}
    \|A_y\mathbf w-\mathbf u_y\|_2^2
    +
    \lambda_v\|\mathbf w\|_2^2.
\end{equation}
Equivalently,
\[
    \mathbf w_x
    =
    (A_x^\top A_x+\lambda_v I)^{-1}A_x^\top \mathbf u_x,
    \qquad
    \mathbf w_y
    =
    (A_y^\top A_y+\lambda_v I)^{-1}A_y^\top \mathbf u_y.
\]
The final input coefficient vector is
\[
    \mathbf w_{\rm obs}
    =
    \begin{bmatrix}
        \mathbf w_x\\
        \mathbf w_y
    \end{bmatrix}
    \in\mathbb R^{2M_v}.
\]
This split encoding avoids forcing the two velocity components to share the same coefficient vector. 
Before training the neural network, the input coefficients are standardized componentwise using the training-set statistics:
\[
    \widetilde{\mathbf w}_{\rm obs}
    =
    \frac{
        \mathbf w_{\rm obs}-\boldsymbol\mu_w
    }{
        \boldsymbol\sigma_w+\epsilon_{\rm std}
    },
\]
where $\boldsymbol\mu_w$ and $\boldsymbol\sigma_w$ are computed only from the training set.

\subsection{Coefficient-to-coefficient network}

The inverse model maps the encoded velocity coefficients to SDF coefficients:
\[
    F_\theta:\mathbb R^{2M_v}\to\mathbb R^{M_u}.
\]
Given $\widetilde{\mathbf w}_{\rm obs}$, the predicted SDF coefficient vector is
\[
    \widehat{\mathbf a}
    =
    F_\theta(\widetilde{\mathbf w}_{\rm obs}),
\]
and the full SDF field is decoded by
\[
    \widehat{\boldsymbol\phi}
    =
    A_u\widehat{\mathbf a}.
\]
The network is a fully connected multilayer perceptron with two hidden layers, hidden width $512$, and GELU activation functions.

\subsection{Training objective}

The model is first trained using a relative SDF loss over the full grid:
\begin{equation}
\label{eq:app_rel_loss}
    \mathcal L_{\rm rel}
    =
    \frac{
        \|\widehat{\boldsymbol\phi}-\boldsymbol\phi\|_2
    }{
        \|\boldsymbol\phi\|_2+\varepsilon
    }.
\end{equation}
Here $\boldsymbol\phi$ is the ground-truth SDF evaluated on the same grid as $\widehat{\boldsymbol\phi}$.

Since the obstacle boundary is determined by the zero level set of the SDF, we further fine-tune the model with a boundary-focused weighted loss. 
For each grid point $\mathbf z_i$, define
\[
    W_i
    =
    \begin{cases}
        w_{\rm in}, & |\phi_i|\leq \delta,\\
        w_{\rm out}, & |\phi_i|>\delta.
    \end{cases}
\]
The weighted boundary loss is
\begin{equation}
\label{eq:app_band_loss}
    \mathcal L_{\rm band}
    =
    \frac{
        \sum_{i=1}^{N_g}
        W_i(\widehat\phi_i-\phi_i)^2
    }{
        \sum_{i=1}^{N_g}W_i
    }.
\end{equation}
The fine-tuning objective is
\begin{equation}
\label{eq:app_total_loss}
    \mathcal L
    =
    \alpha\mathcal L_{\rm band}
    +
    (1-\alpha)\mathcal L_{\rm rel}.
\end{equation}
In the experiments, we use
\[
    \delta=0.05,
    \qquad
    w_{\rm in}=20,
    \qquad
    w_{\rm out}=1,
    \qquad
    \alpha=0.8.
\]
This loss emphasizes the accuracy of the zero level set while maintaining a globally consistent SDF field.

\subsection{Boundary extraction and visualization}

After prediction, the recovered obstacle boundary is extracted as the zero level set
\[
    \widehat\Gamma
    =
    \{\mathbf x:\widehat\phi(\mathbf x)=0\}.
\]
In the reconstruction plots, we compare $\widehat\Gamma$ with the ground-truth boundary
\[
    \Gamma
    =
    \{\mathbf x:\phi(\mathbf x)=0\}.
\]
We also visualize the predicted near-boundary band
\[
    \mathcal N_{\widehat\phi}
    =
    \{\mathbf x:|\widehat\phi(\mathbf x)|\leq 0.15\}.
\]
This band provides a qualitative indication of the uncertainty and thickness of the predicted boundary region.

\section{Proofs for Section~\ref{sec:basis-coefficient-analysis}}
\label{app:analysis-proofs}

\subsection{Proof of Lemma~\ref{lem:regularized-encoder-stability}}

\begin{proof}
Let $\Phi=U\Sigma V^\top$ be a singular value decomposition. 
From \eqref{eq:ridge-encoder-definition},
\[
E_\lambda
=
V(\Sigma^\top\Sigma+\lambda I)^{-1}\Sigma^\top U^\top.
\]
It follows that
\[
\|E_\lambda\|_2
=
\max_i
\frac{\sigma_i(\Phi)}
{\sigma_i(\Phi)^2+\lambda}.
\]
For $s\geq0$, let $h_\lambda(s)=s/(s^2+\lambda)$. 
Since $h_\lambda'(s)=(\lambda-s^2)/(s^2+\lambda)^2$, the maximum is attained at $s=\sqrt{\lambda}$ and equals $1/(2\sqrt{\lambda})$. 
This proves \eqref{eq:encoder-stability-main}.

Finally,
$\mathbf a_\lambda(f)-\mathbf a_\lambda(g)=E_\lambda(\mathbf f-\mathbf g)$, and therefore
\[
\|\mathbf a_\lambda(f)-\mathbf a_\lambda(g)\|_2
\leq
\|E_\lambda\|_2
\|\mathbf f-\mathbf g\|_2.
\]
Using \eqref{eq:encoder-stability-main} proves \eqref{eq:encoder-lipschitz}.
\end{proof}

\subsection{Proof of Proposition~\ref{prop:encoding-bias-ridge}}

\begin{proof}
Using $\Phi=U\Sigma V^\top$, the reconstruction operator satisfies
\[
H_\lambda
=
U\Sigma(\Sigma^\top\Sigma+\lambda I)^{-1}\Sigma^\top U^\top.
\]
For each left singular vector $u_i$ corresponding to a nonzero singular value $\sigma_i$,
\[
H_\lambda u_i
=
\frac{\sigma_i^2}{\sigma_i^2+\lambda}u_i.
\]
Moreover, $H_\lambda\mathbf f_\perp=0$ whenever
$\mathbf f_\perp\perp\operatorname{Range}(\Phi)$. 
Thus,
\[
\widetilde{\mathbf f}_\lambda
=
H_\lambda\mathbf f
=
\sum_{i=1}^r
\frac{\sigma_i^2}{\sigma_i^2+\lambda}
\alpha_i u_i,
\]
which proves \eqref{eq:regularized-reconstruction-spectral}.

Furthermore,
\[
\mathbf f-\widetilde{\mathbf f}_\lambda
=
\sum_{i=1}^r
\frac{\lambda}{\sigma_i^2+\lambda}
\alpha_i u_i
+
\mathbf f_\perp.
\]
The terms in this decomposition are mutually orthogonal. 
Taking the squared norm proves \eqref{eq:encoding-bias-spectral}.
\end{proof}

\subsection{Proof of Proposition~\ref{prop:output-space-lower-bound}}

\begin{proof}
Since $P_{\mathrm{out}}$ is the orthogonal projection onto
$V_{\mathrm{out}}^{m_2}$,
\[
G(f)-P_{\mathrm{out}}G(f)
\perp
V_{\mathrm{out}}^{m_2}.
\]
Moreover,
$P_{\mathrm{out}}G(f)-\widehat G_{\theta,\lambda}(f)
\in V_{\mathrm{out}}^{m_2}$. 
Therefore,
\[
G(f)-\widehat G_{\theta,\lambda}(f)
=
\bigl(G(f)-P_{\mathrm{out}}G(f)\bigr)
+
\bigl(P_{\mathrm{out}}G(f)-\widehat G_{\theta,\lambda}(f)\bigr)
\]
is an orthogonal decomposition. 
The Pythagorean identity gives
\[
\|G(f)-\widehat G_{\theta,\lambda}(f)\|_Y^2
=
\|G(f)-P_{\mathrm{out}}G(f)\|_Y^2
+
\|P_{\mathrm{out}}G(f)-\widehat G_{\theta,\lambda}(f)\|_Y^2.
\]
Using
$\|G(f)-P_{\mathrm{out}}G(f)\|_Y
=\varepsilon_{\mathrm{out}}(G(f))$
proves \eqref{eq:output-projection-identity}. 
The lower bound \eqref{eq:output-projection-lower-bound} follows immediately.
\end{proof}

\subsection{Proof of Lemma~\ref{lem:lambda-loss-bound}}

\begin{proof}
By Lemma~\ref{lem:regularized-encoder-stability} and Assumption~\ref{assump:generalization-main},
\[
\|\mathbf a_\lambda(f)\|_2
=
\|E_\lambda\mathbf f\|_2
\leq
\frac{R_f}{2\sqrt{\lambda}}.
\]
The affine-growth assumption gives
\[
\|N_\theta(\mathbf a_\lambda(f))\|_2
\leq
B_N+\frac{L_NR_f}{2\sqrt{\lambda}}.
\]
Using the decoder bound and the triangle inequality,
\begin{align*}
\left\|
D_{\mathrm{out}}
\left(
N_\theta(\mathbf a_\lambda(f))-\mathbf b(f)
\right)
\right\|_Y
&\leq
\|D_{\mathrm{out}}\|_{\mathrm{op}}
\left(
\|N_\theta(\mathbf a_\lambda(f))\|_2
+
\|\mathbf b(f)\|_2
\right)\\
&\leq
C_D
\left(
B_b+B_N+\frac{L_NR_f}{2\sqrt{\lambda}}
\right)
=
B_\lambda.
\end{align*}
Squaring both sides proves
$0\leq\ell_{\theta,\lambda}(f)\leq B_\lambda^2$.
\end{proof}

\subsection{Proof of Theorem~\ref{thm:lambda-dependent-generalization}}

\begin{proof}
By Lemma~\ref{lem:lambda-loss-bound}, every function in
$\mathcal F_\lambda=\{\ell_{\theta,\lambda}:\theta\in\Theta\}$
takes values in $[0,B_\lambda^2]$. 
By Assumption~\ref{assump:generalization-main},
$\operatorname{Pdim}(\mathcal F_\lambda)\leq\Pi$. 
Applying the standard pseudo-dimension generalization estimate for bounded real-valued function classes~\cite{anthony1999neural,shalev2014understanding} gives
\[
\mathbb E
\sup_{\theta\in\Theta}
\left|
\mathcal R_\lambda(\theta)
-
\widehat{\mathcal R}_{M,\lambda}(\theta)
\right|
\leq
C
B_\lambda^2
\sqrt{\frac{\Pi}{M}}
\log(eM).
\]
This proves \eqref{eq:lambda-dependent-generalization-bound}.
\end{proof}

\subsection{Proof of Corollary~\ref{cor:erm-excess-risk-lambda}}

\begin{proof}
Define
\[
\Delta_{M,\lambda}
=
\sup_{\theta\in\Theta}
\left|
\mathcal R_\lambda(\theta)
-
\widehat{\mathcal R}_{M,\lambda}(\theta)
\right|.
\]
For every $\theta\in\Theta$,
\[
\mathcal R_\lambda(\theta)
\leq
\widehat{\mathcal R}_{M,\lambda}(\theta)
+
\Delta_{M,\lambda},
\qquad
\widehat{\mathcal R}_{M,\lambda}(\theta)
\leq
\mathcal R_\lambda(\theta)
+
\Delta_{M,\lambda}.
\]
Using \eqref{eq:approximate-erm-condition},
\begin{align*}
\mathcal R_\lambda(\widehat\theta_\lambda)
&\leq
\widehat{\mathcal R}_{M,\lambda}(\widehat\theta_\lambda)
+
\Delta_{M,\lambda}\\
&\leq
\inf_{\theta\in\Theta}
\widehat{\mathcal R}_{M,\lambda}(\theta)
+
\eta_{\mathrm{opt}}
+
\Delta_{M,\lambda}\\
&\leq
\inf_{\theta\in\Theta}
\mathcal R_\lambda(\theta)
+
\eta_{\mathrm{opt}}
+
2\Delta_{M,\lambda}.
\end{align*}
Taking expectations and applying
Theorem~\ref{thm:lambda-dependent-generalization} gives
\[
\mathbb E\mathcal R_\lambda(\widehat\theta_\lambda)
\leq
\inf_{\theta\in\Theta}
\mathcal R_\lambda(\theta)
+
\eta_{\mathrm{opt}}
+
2C
B_\lambda^2
\sqrt{\frac{\Pi}{M}}
\log(eM).
\]
Absorbing the factor $2$ into the universal constant proves
\eqref{eq:erm-excess-risk-lambda}.
\end{proof}

\end{document}